\documentclass[reqno, 10pt, centertags,draft]{amsart}
\usepackage{amsmath,amsthm,amscd,amssymb,latexsym,upref}

\makeatletter
\def\theequation{\@arabic\c@equation}

\newcommand{\bbN}{{\mathbb{N}}}
\newcommand{\bbR}{{\mathbb{R}}}

\newcommand{\bbZ}{{\mathbb{Z}}}
\newcommand{\bbC}{{\mathbb{C}}}

\newcommand{\C}{\mathbb{C}}
\newcommand{\R}{\mathbb{R}}
\newcommand{\Z}{\mathbb{Z}}
\newcommand{\N}{\mathbb{N}}

\newcommand{\cA}{{\mathcal A}}
\newcommand{\cB}{{\mathcal B}}
\newcommand{\cC}{{\mathcal C}}

\newcommand{\cG}{{\mathcal G}}
\newcommand{\cH}{{\mathcal H}}

\newcommand{\cK}{{\mathcal K}}

\newcommand{\cM}{{\mathcal M}}

\newcommand{\cQ}{{\mathcal Q}}

\newcommand{\cS}{{\mathcal S}}

\newcommand{\no}{\nonumber}
\newcommand{\lb}{\label}
\newcommand{\f}{\frac}

\newcommand{\ol}{\overline}

\newcommand{\wti}{\widetilde  }

\newcommand{\loc}{\text{\rm{loc}}}

\newcommand{\ran}{\text{\rm{ran}}}

\newcommand{\dom}{\text{\rm{dom}}}

\newcommand{\slim}{\text{\rm{s-lim}}}
\newcommand{\slimes}{\text{\rm{l.i.m.}}}

\newcommand{\supp}{\text{\rm{supp}}}

\newcommand{\res}{\text{\rm{res}}}

\newcommand{\bi}{\bibitem}

\newcommand{\tr}{\text{\rm{tr}}}

\renewcommand{\Re}{\text{\rm Re}}
\renewcommand{\Im}{\text{\rm Im}}
\renewcommand{\ln}{\text{\rm ln}}

\numberwithin{equation}{section}

\newtheorem{theorem}{Theorem}[section]

\newtheorem{lemma}[theorem]{Lemma}
\newtheorem{corollary}[theorem]{Corollary}
\newtheorem{hypothesis}[theorem]{Hypothesis}
\theoremstyle{definition}
\newtheorem{definition}[theorem]{Definition}
\newtheorem{remark}[theorem]{Remark}

\begin{document}

\title[Hill Operators and Spectral Operators of Scalar Type]{A Criterion
for Hill Operators to be \\ Spectral Operators of Scalar Type}
\author[F.\ Gesztesy and V.\ Tkachenko]{Fritz
Gesztesy and Vadim Tkachenko}
\address{Department of Mathematics,
University of Missouri, Columbia, MO 65211, USA}
\email{fritz@math.missouri.edu}
\urladdr{http://www.math.missouri.edu/personnel/faculty/gesztesyf.html}
\address{Department of Mathematics,
Ben Gurion University of the Negev, Beer--Sheva 84105, Israel}
\email{tkachenk@math.bgu.ac.il}
\date{\today}
\thanks{Based upon work supported by the National Science
Foundation under Grant No.\ DMS-0405526 and  Israel Science Foundation under 
Grant No.\ 186/01.}

\begin{abstract}

We derive necessary and sufficient conditions for a Hill operator (i.e., a 
one-dimensional periodic Schr\"odinger operator) $H=-d^2/dx^2+V$ to be a spectral
operator of scalar type. The conditions show the remarkable fact that the
property of a Hill operator being a spectral operator is independent of
smoothness (or even analyticity) properties of the potential $V$. In the course of our 
analysis we also establish a functional model for periodic Schr\"odinger operators 
that are spectral operators of scalar type and develop the corresponding eigenfunction expansion. 

The problem of deciding which Hill operators are spectral
operators of scalar type appears to have been open for about 40 years. 
\end{abstract}

\maketitle

\section{Introduction} \lb{s1}

The principal aim of this paper is to establish a criterion for deciding
when a Hill operator (i.e., a one-dimensional periodic Schr\"odinger operator)  
\begin{equation}\label{1.00}
H=-\frac{d^2}{dx^2}+V(x) 
\end{equation}
in $L^2(\bbR)$ with $V(x)$ periodic, is a spectral
operator of scalar type in the sense of Dunford. To describe this
longstanding open problem in more detail requires a bit of preparation and
so we present a brief historical introduction of topics closely related to
the material in this paper.

For the notion of spectral operators we refer to the classical monograph of
Dunford and Schwartz \cite{DS88a} (cf. Appendix \ref{A} for a very brief
summary). Generally speaking, a spectral operator $T$ in a Hilbert space $\cH$
possesses a spectral measure $E_T(\cdot)$ (i.e., a homomorphism from the
$\sigma$-algebra of Borel subsets of $\bbC$ into the Boolean algebra of
projection operators on $\cH$ with $E_T$ bounded and
$E_T(\bbC)=I_{\cH}$) such that $E_T(\omega)$ commutes with $T$ and
$\sigma\big(T|_{E_{T}(\omega)\cH}\big)\subseteq \ol\omega$ for all Borel
subsets $\omega$ of $\bbC$ (cf.\ Definition \ref{d2.14}). Here
$\sigma(\cdot)$ denotes the spectrum. Spectral operators of scalar type $S$
are then of the form
$S=\slim_{n\uparrow\infty}\int_{|\lambda|\leq n} \lambda \, dE_S(\lambda)$
(cf.\ Definition \ref{d2.17}).

The literature on self-adjoint Hill operators of the form
$H=-d^2/dx^2+V$ in $L^2(\bbR)$ with $V$ real-valued and periodic, is too
enormous to be reviewed here, so we just note that the basic spectral
properties of such $H$, namely, a countable set of closed intervals which may
degenerate into finitely many closed intervals and a half-line, and an
explicit formula for the spectral expansions generated by $H$, have been
established since 1950 (see, e.g., \cite[Ch.\ 6]{Ea73},
\cite[Sect.\ XIII.16]{RS78}, \cite{Ti50}, \cite[Ch.\ XXI]{Ti58}). In
particular, since all self-adjoint operators are spectral operators in
the sense of \cite{DS88a}, we immediately turn to the case of periodic
but non-self-adjoint operators $H$ with complex-valued and periodic
potentials $V$. In this case the spectrum shows much more complexity
compared to the self-adjoint situation. Indeed, as shown by Serov \cite{Se60} 
(see also \cite{Ro63}, \cite{Tk64}, \cite{Mc65}) in
the first half of the 1960's, the spectrum now consists of a countable
system of piecewise analytic arcs which may in fact exhibit crossings
(again this can degenerate into a finite system of piecewise analytic
arcs and a semi-infinite arc). However, necessary and sufficient
conditions (in fact, even just sufficient conditions) for the existence
of a uniformly bounded family of spectral projections of a
(non-self-adjoint) Hill operator $H$, and especially, the property of
$H$ being a spectral operator of scalar type, remained elusive since the
mid 1960's. Following up on successful applications of the formalism of
spectral operators of scalar type to Schr\"odinger operators on the
half-line $(0,\infty)$ with sufficiently fast decaying potentials as
$|x|\uparrow\infty$ (cf.\ the discussions in \cite[Sect.\ XX.1]{DS88a},
\cite{Na60}, and \cite[Appendix II]{Na68}), McGarvey
\cite{Mc62}--\cite{Mc66} started a systematic study of periodic
differential operators (including higher-order differential operators) in a
series of papers in 1962--1966. Using a combination of direct integral
decompositions and perturbation techniques, McGarvey \cite{Mc65a} proved
in 1965 that certain $n$th-order differential operators in $L^2(\bbR)$
with periodic coefficients and certain second-order periodic differential
operators of the form
\begin{equation}
- \frac{d^2}{dx^2} +p(x)\frac{d}{dx} + q(x), \quad x\in\mathbb{R}, \,
\label{1.1a}
\end{equation}
with $\pi$-periodic functions $p$ and $q$ under the restriction
\begin{equation}
\Im\bigg(\int_0^\pi dx \, p(x)\bigg) \neq 0,
\end{equation}
are spectral operators at infinity. The spectra of operators \eqref{1.1a}
outside a sufficiently large disc  are composed of some separated ovals,
permitting McGarvey to prove that these operators are in some
sense  {\em asymptotically} spectral operators \cite{Mc65a}. Such results
ignore the existence of local spectral singularities and, at any rate, are
not applicable to non-self-adjoint Hill operators $H=-d^2/dx^2+V$ in
$L^2(\bbR)$.

In spite of a flurry of activities in connection with spectral theory for
non-self-adjoint Hill operators since the early eighties, many of which
were inspired by connections to the Korteweg--deVries hierarchy of evolution
equations (we refer, e.g., to \cite{BG06}--\cite{Ch06}, 
\cite{DM02}--\cite{DM03a}, \cite{Ga80}, \cite{Ga80a},
\cite{GW96}--\cite{GU83}, \cite{Ko97}, \cite{PT88},
\cite{PT91}, \cite{ST96}--\cite{ST97}, \cite{Sh03}--\cite{Sh04a},
\cite{Tk92}--\cite{Tk02}, \cite{We98}, \cite{We98a}), to the best of our
knowledge, no progress on the question of which Hill operators are
spectral operators of scalar type was made since McGarvey's investigations
in the early 1960's. In Section 8 we will further discuss results from  
\cite{EFZ05}, \cite{Me77}, \cite{Mi06}, \cite{Ve80}--\cite{Vo63}, \cite{Zh69} on spectral expansions associated with a non-self-adjoint 
Hill operator $H$.

At first sight, a natural expectation with respect to such non-self-adjoint operators would be the following: The ``better'' the properties of its potential $V$ and the ``smaller'' its imaginary part, the better should be its chance to be a scalar spectral operator. From this (admittedly, perhaps a bit naive) point of view, operators with potentials $V$ analytic in a half-plane $\Im (z) > a$, $a\leq 0$, should be the best possible candidates. However, it follows from results obtained by Gasymov \cite{Ga80}, that actually no such operator with non-constant potential is spectral. In other words, no smoothness or analyticity conditions imposed on $V$ can guarantee that $H$ is a spectral operator of scalar type.

In this paper we prove two versions of a criterion for a Hill operator to be a spectral operator of scalar type, one analytic and one geometric. The analytic version is stated in terms of Hill's discriminant $\Delta_+(z)=[\theta(z,\pi)+\phi'(z,\pi)]/2$ and the functions $\Delta_-(z)=[\theta(z,\pi)-\phi'(z,\pi)]/2$
and $\phi(z,\pi)$, where $\{\theta(z,x),\phi(z,x)\}$ is a fundamental system of distributional solutions of $H\psi(z,x)=z\psi(z,x)$, satisfying canonical boundary conditions at $x=0$ recorded in \eqref{3.3}. The triple $\Delta_+(z)$, $\Delta_-(z)$, and $\phi(z,\pi)$ was introduced in \cite{ST96}, \cite{Tk92} as a complete system of independent parameters which uniquely determines the potential $V$. The geometric version of the criterion uses algebraic and geometric properties of spectra of periodic/antiperiodic and Dirichlet boundary value problems generated by $H$ in the space $L^2([0,\pi])$.

Finally we briefly describe the content of each section. Section \ref{s3}
contains basic facts on Floquet theory and 
a summary of the spectral results  of (non-self-adjoint) Hill operators 
 (see \cite{Se60}, \cite{Ro63}, \cite{Tk64}, \cite{Mc65}) and some comments on
direct integral decompositions of $H$ and its Green's function. Section
\ref{s4} then summarizes our principal new results. We provide three
theorems (Theorems \ref{t4.3}--\ref{t4.5} announced in \cite{GT06}) which 
each provide a criterion
for $H$ to be a spectral operator of scalar type. Section \ref{s5} provides a detailed discussion of the reduced operators
$H(t)$, $t\in [0,2\pi]$, in $L^2([0,\pi])$, their spectra, and
their spectral expansion. Here $H$ is the direct integral over $H(t)$,
$H=(2\pi)^{-1} \int_{[0,2\pi]}^{\oplus} dt\, H(t)$. Section \ref{s6} proves
some auxiliary results and Section \ref{s7} proves the necessity of our
conditions for $H$ to be a spectral operator of scalar type. The proof of
sufficiency of our conditions for $H$ to be a spectral operator of scalar
type is then presented in Section \ref{s8}. Our final Section \ref{s9}
presents a series of concluding remarks which put our results in proper
perspective and underscores the subtleties involved when trying to determine 
whether or not a Hill operator is a spectral operator of scalar type.

\section{Preliminaries on Hill Operators and Floquet Theory} \lb{s3}

In this section we briefly recall some standard results on (not necessarily 
self-adjoint) Hill operators and their associated  Floquet theory.

Throughout this paper the $L^2(\Omega)$ and $L^2_{\loc}(\Omega)$ spaces
without specifying the corresponding measure on the set
$\Omega\subseteq\bbR$ refer to Lebesgue measure on $\Omega$. The scalar
product in $L^2(\Omega)$ will be denoted by $(\cdot,\cdot)_{L^2(\Omega)}$
(it is assumed to be linear in the second factor), the corresponding
norm is denoted by $\|\cdot\|_{L^2(\Omega)}$. For simplicity, the identity
operators in $L^2(\Omega)$ will be denoted by $I$. We use the symbol
$\prime$ to denote $x$-derivatives and $\bullet$ to denote
$z$-derivatives (and occasionally, $\zeta$ derivatives, where
$z=\zeta^2$). Moreover, we  denote by $\sigma(\cdot)$,
$\sigma_{\rm p}(\cdot)$, $\sigma_{\rm r}(\cdot)$, $\sigma_{\rm c}(\cdot)$,
and $\rho(\cdot)$ the spectrum, the point spectrum (i.e., the set of eigenvalues), the residual spectrum,  the continuous spectrum, and
resolvent set of a densely defined, closed, linear operator in a Hilbert
space. We also use the abbreviation $\bbN_0=\bbN\cup\{0\}$.

For the remainder of this paper we assume the following hypothesis.
\begin{hypothesis} \lb{h3.1}
Suppose
\begin{equation}
V\in L^2_{\loc}(\bbR), \quad V(x+\pi)=V(x) \, \text{ for a.e.\
$x\in\bbR$.}  \lb{3.1}
\end{equation}
\end{hypothesis}

Without loss of generality we chose the period of $V$ to be $\pi$ for
subsequent notational convenience.

Given Hypothesis \ref{h3.1}, one introduces the differential expression
\begin{equation}\lb{3.1a}
L = - \f{d^2}{dx^2} + V(x), \quad x\in\bbR \,
\text{ (or $x\in [0,\pi]$)} 
\end{equation}
and defines the corresponding Schr\"odinger
operator $H$ in $L^2(\bbR)$ by
\begin{align}
\begin{split}
& (Hf)(x)=(L f)(x),  \quad x\in\bbR, \\
& f\in\dom(H)=\{g\in L^2(\bbR) \,|\, g,g'\in AC_{\loc}(\bbR); \,
L g \in L^2(\bbR)\}.   \lb{3.2}
\end{split}
\end{align}
Then $H$ is known to be a densely defined, closed, linear operator in
$L^2(\bbR)$. It is self-adjoint if and only if $V$ is real-valued.

Associated with $H$ one introduces the fundamental system of
distributional solutions $\theta(z,\cdot)$ and $\phi(z,\cdot)$ of
$L \psi=z\psi$ satisfying
\begin{equation}
\theta(z,0)=\phi'(z,0)=1, \quad
\theta'(z,0)=\phi(z,0)=0, \quad z\in\bbC.  \lb{3.3}
\end{equation}
For each $x\in\bbR$, $\theta(z,x)$ and $\phi(z,x)$ are entire with respect
to $z$. The monodromy matrix
$\cM(z)$ is then given by
\begin{equation}
\cM(z)=\begin{pmatrix} \theta(z,\pi) & \phi(z,\pi) \\
\theta'(z,\pi) & \phi'(z,\pi) \end{pmatrix}, \quad
z\in\bbC  \lb{3.4}
\end{equation}
and its eigenvalues $\rho_\pm(z)$, the Floquet multipliers, satisfy
\begin{equation}
\rho_+(z)\rho_-(z)=1  \lb{3.5}
\end{equation}
since
\begin{equation}
\det(\cM(z))=\theta(z,\pi)\phi'(z,\pi)-\theta'(z,\pi)\phi(z,\pi)=1.
\lb{3.5a}
\end{equation}
The Floquet discriminant $\Delta_+(\cdot)$ is
then defined by
\begin{equation}
\Delta_+(z)=\tr(\cM(z))/2=
[\theta(z,\pi)+\phi'(z,\pi)]/2, \quad z\in\bbC,  \lb{3.6}
\end{equation}
and one obtains
\begin{equation}
\rho_\pm (z)=\Delta_+(z)\pm i\sqrt{1-\Delta_+(z)^2}  \lb{3.7}
\end{equation}
with an appropriate choice of the square root branches. We also note that
\begin{equation}
|\rho_\pm(z)|=1 \, \text{ if and only if } \, \Delta_+(z)\in [-1,1].
\lb{3.8}
\end{equation}

The following theorem describes the well-known fundamental properties of the spectrum of (non-self-adjoint) Hill operators in \eqref{3.2} (cf.\ \cite{Se60}, \cite{Ro63}, \cite{Tk64}, \cite{Mc65}, 
\cite[p.\ 1486--98]{DS88}).

\begin{theorem} \lb{t3.2} Assume Hypothesis \ref{h3.1}. \\
$(i)$ The point spectrum and residual spectrum of $H$ are empty and hence
the spectrum of $H$ is purely continuous,
\begin{align}
\sigma_{\rm p}(H)&=\sigma_{\rm r}(H)=\emptyset,  \lb{3.10} \\
\sigma(H)&=\sigma_{\rm c}(H).   \lb{3.11}
\end{align}
$(ii)$ $\sigma(H)$ is given by
\begin{align}
\sigma(H) &=\big\{\lambda\in\bbC\,\big|\, -1\leq \Delta_+(\lambda)
\leq 1\big\}   \lb{3.12} \\
&=\{\lambda\in\bbC \,|\, \text{there exists at least one non-trivial
bounded}    \no \\
& \hspace*{.65cm} \text{distributional solution $\psi\in L^\infty(\bbR;dx)$ of
$H\psi=\lambda\psi$}\}.  \lb{3.13}
\end{align}
The latter set equals the conditional stability set of $H$.
In addition, $\sigma(H)$ contains no isolated points.  \\
$(iii)$ $\sigma(H)$ is contained in the semi-strip
\begin{equation}
\sigma(H)\subset \{z\in\bbC \,|\, \Im(z)\in [M_1,M_2], \, \Re(z)\geq
M_3\},  \lb{3.14}
\end{equation}
where
\begin{equation}
M_1=\inf_{x\in\bbR}[\Im(V(x))], \quad
M_2=\sup_{x\in\bbR}[\Im(V(x))], \quad M_3=
\inf_{x\in\bbR}[\Re(V(x))] \lb{3.15}
\end{equation}
$(iv)$ Qualitatively, $\sigma(H)$ consists of countably many, simple,
analytic arcs which may degenerate into finitely many analytic arcs and
one simple semi-infinite analytic arc. Asymptotically, these analytic
arcs approach the half-line
\begin{equation}
L_{\langle V\rangle}=\{z\in\bbC \,|\, z=\langle V\rangle +x, \,
x\geq 0\}.
\end{equation}
Moreover, crossings of arcs $($subject to certain
restrictions, see, e.g., item $(v)$$)$ are permitted.
\\
$(v)$ The resolvent set $\bbC\backslash\sigma (H)$ of $H$ is
path-connected.
\end{theorem}

Here the mean value $\langle h \rangle$ of a periodic function $h \in
L^1_{\loc}(\bbR)$ of period $\pi>0$ is given by
\begin{equation}
\langle h \rangle=\f{1}{\pi} \int_{x_0}^{x_0+\pi} dx\, h(x),
\lb{3.16}
\end{equation}
independent of the choice of $x_0\in\bbR$.

\begin{remark} \lb{r3.3}
A set $\sigma\subset\bbC$ is called an {\it arc} if there exists a
parameterization $\gamma\in C([0,1])$ such that
$\sigma=\{\gamma(t)\in\bbC \,|\, t\in [0,1]\}$. The arc $\sigma$ is
called {\it simple} if there exists a parameterization
$\gamma$ such that $\gamma\colon [0,1]\to\bbC$ is injective. The arc
$\sigma$ is called {\it analytic} if there is a parameterization $\gamma$
that is analytic at each $t\in (0,1)$. Finally, $\sigma_\infty$ is called
a {\it semi-infinite} arc if there exists a parameterization $\gamma\in
C([0,\infty))$ such that $\sigma_\infty=\{\gamma(t)\,|\, t\in
[0,\infty)\}$ and $\sigma_\infty$ is an unbounded subset of $\bbC$.
Analytic semi-infinite arcs are defined analogously, and by a simple
semi-infinite arc we mean one that is without self-intersection (i.e.,
corresponds to a injective parameterization) with the additional
restriction that the unbounded part of $\sigma_\infty$ consists of
precisely one branch tending to infinity.
\end{remark}

Next, we take a closer look at spectral arcs of $H$.
Let $\Delta_+$ be defined as in  \eqref{3.6}, 
and let $\lambda_0\in\sigma(H)$. Then there exists $t_0\in[0,\pi]$
such that $\Delta_+(\lambda_0)=\cos (t_0)$. If $\lambda_0\in\sigma(H)$
with $\Delta_+^{\bullet}(\lambda_0)\neq 0$, then there exist closed
intervals
$[\alpha,\beta] \subset[0,\pi]$ and a function  $\lambda(\cdot)$
continuous on $[\alpha,\beta]$, analytic in an open neighborhood of
$[\alpha,\beta]$, such that $\lambda (t_0)=\lambda_0,
\; t_0\in [\alpha,\beta]$ and
\begin{align}
\begin{split}
&\Delta_+(\lambda(t))=\cos(t), \;
t\in[\alpha,\beta],  \lb{3.16a} \\
&\Delta_+^{\bullet}(\lambda(t)) \neq 0,\quad \lambda'(t) \neq 0, \;
t\in(\alpha,\beta).
\end{split}
\end{align}

\begin{definition} \lb{d3.4} A closed spectral arc of $H$
\begin{equation}
\sigma=\{z\in\bbC \,|\, z=\lambda(t), \, t\in[\alpha,\beta]\}  \lb{3.16b}
\end{equation}
with
$\Delta_+^{\bullet}(z)\neq 0$ for all $z\in\sigma$ will
be called a {\it regular spectral arc} of $H$,
and the points  $\lambda(\alpha)$ and $\lambda(\beta)$ will
be called its {\it endpoints}.
\end{definition}

It follows from this definition that all regular spectral arcs of $H$
are compact subsets of $\bbC$.

Finally, we assume that the orientation of a regular spectral arc $\sigma$
of $H$
is induced by the function $\lambda (t)$ as $t$ varies from $0$ to $\pi$
(cf.\ \eqref{3.16a}). In the special self-adjoint case, this has the
effect that for every odd-numbered spectral band (we index them
by $k\in\bbN$, see also \eqref{3.29}, \eqref{3.30}, and Lemma \ref{l5.1}),
$\lambda(0) < \lambda(\pi)$, whereas for every even-numbered spectral band
$\lambda(0) > \lambda(\pi)$. This is consistent with integrals such as
\eqref{defproj} for $P(\sigma)$ since we assume that 
\begin{equation}
\sqrt{1-\Delta_+(\lambda)^2} \geq 0, \quad \lambda \in \sigma(H).
\end{equation}

Next, we denote by
\begin{equation}
\{\mu_k\}_{k\in\bbN}=\{z\in\bbC \,|\, \phi(z,\pi)=0\}   \lb{3.17}
\end{equation}
the set of zeros of the entire function $\phi(\cdot,\pi)$. The set
$\{\mu_k\}_{k\in\bbN}$ represents the Dirichlet spectrum associated with
the restriction of $L$ to the interval $[0,\pi]$. More precisely,
define the densely defined, closed, linear  
operator $H^{D}$ in $L^2([0,\pi])$ with Dirichlet boundary conditions by
\begin{align}
& (H^Df)(x)=(L f)(x),  \quad x \in [0,\pi],  \lb{3.18} \\
& f\in\dom(H^D)=\{g\in L^2([0,\pi]) \,|\, g,g'\in AC([0,\pi]); \,
L g \in L^2([0,\pi]); \no \\
& \hspace*{7.8cm} g(0)=g(\pi)=0\}  \no
\end{align}
($H^D$ is self-adjoint if and only if $V$ is real-valued), then
\begin{equation}
\sigma(H^D)=\{\mu_k\}_{k\in\bbN}.   \lb{3.19}
\end{equation}

In a similar fashion one defines Schr\"odinger operator
$H^N$ in $L^2([0,\pi])$ with the Dirichlet boundary condition in
\eqref{3.18} replaced by the Neumann boundary condition 
\begin{equation}
g'(0)=g'(\pi)=0.  \lb{3.19a}
\end{equation}
Moreover, we also mention the family of Schr\"odinger operators
$H^{\alpha}$, $\alpha\in\bbR$, in $L^2([0,\pi])$, replacing the boundary
condition in \eqref{3.18} by
\begin{equation}
g'(0)+\alpha g(0)=g'(\pi)+\alpha g(\pi)=0, \quad \alpha\in\bbR.
\lb{3.19b}
\end{equation}
Of course, $H^N=H^0$ and formally, $H^D=H^\infty$.

For future purposes it is convenient to introduce
\begin{equation}
\Delta_-(z)=[\theta(z,\pi)-\phi'(z,\pi)]/2, \quad z\in\bbC. \lb{3.19c}
\end{equation}
Floquet solutions $\psi_\pm(z,\cdot)$ of $L \psi=z\psi$ 
normalized at $x=0$ associated with $H$ are then given by 
($z\in\bbC\backslash\{\mu_j\}_{j\in\bbN}$)
\begin{align}
\psi_\pm(z,x)&=\theta(z,x)+[\rho_\pm(z)-\theta(z,\pi)]
\phi(z,\pi)^{-1}\phi(z,x) \no \\
& =\theta(z,x)+m_\pm(z)\phi(z,x),\
\quad m_\pm(z)=\frac
{-\Delta_-(z)\pm i\sqrt{1-\Delta_+(z)^2}}{\phi(z,\pi)}, \lb{3.20} \\
\psi_\pm(z,0)&=1. \no
\end{align}
One then verifies (for $z\in\bbC\backslash\{\mu_j\}_{j\in\bbN}$, $x\in\bbR$),
\begin{align}
&\psi_{\pm}(z,x+\pi)=\rho_\pm(z) \psi_\pm(z,x),  \no \\
& \hspace*{1.9cm} =e^{\pm it}\psi_{\pm}(z,x) \, \text{ with
$\Delta_+(z)=\cos(t)$,}  \lb{l24} \\
& W(\psi_+(z,\cdot),\psi_-(z,\cdot))= m_-(z)- m_+(z)= -
\f{2i\sqrt{1-\Delta_+ (z)^2}}{\phi(z,\pi)},  \lb{3.21} \\
& m_+(z)+m_-(z)=-2\Delta_-(z)/\phi(z,\pi),  \lb{3.22} \\
& m_+(z)m_-(z)=-\theta'(z,\pi)/\phi(z,\pi), \lb{3.23}
\end{align}
where $W(f,g)=fg'-f'g$ denotes the Wronskian of $f$ and $g$.
Moreover, applying Lagrange's formula one computes 
\begin{equation}
\Delta_+^{\bullet}(z)=-\phi(z,\pi)\f{1}{2} \int_{0}^{\pi}
dx \, \psi_+(z,x)\psi_-(z,x), \quad z\in\bbC.  \lb{3.24}
\end{equation}

To describe the Green's function of $H$ (i.e., the integral kernel of the
resolvent of $H$) we assume that the square root in \eqref{3.7} is chosen
such that
\begin{equation}
|\rho_+(z)|<1, \quad |\rho_-(z)|>1, \quad
z\in\bbC\backslash\sigma(H).  \lb{3.24a}
\end{equation}
Then
\begin{equation}
\psi_\pm(z,\cdot) \in L^2([x_0,\pm\infty)) \, \text{ for all } \,
x_0\in\bbR, \;
z\in\bbC\backslash(\sigma(H)\cup\{\mu_j\}_{j\in\bbN})  \lb{3.24b}
\end{equation}
and the Green's function of $H$ is of the form
\begin{align}
&G(z,x,y)=(H-zI)^{-1}(x,y) \no \\
& \quad =\f{1}{W(\psi_+(z),\psi_-(z))}\begin{cases}
\psi_-(z,x)\psi_+(z,y), & x\leq y, \\
\psi_-(z,y)\psi_+(z,x), & x\geq y, \end{cases} \no \\
& \quad =-\f{\phi(z,\pi)}{2i\sqrt{1-\Delta_+ (z)^2}}
\begin{cases}
\psi_-(z,x)\psi_+(z,y), & x\leq y, \\
\psi_-(z,y)\psi_+(z,x), & x\geq y, \end{cases} \quad
z\in\bbC\backslash\sigma(H).   \lb{3.24c}
\end{align}
By \eqref{3.5a} and \eqref{3.23} one infers that the singularities of
$\psi_\pm(z,x)$ at
$z=\mu_j$, $j\in\bbN$, cancel in \eqref{3.24c} and hence the latter
indeed extends to all $z\in\bbC\backslash\sigma(H)$.

For subsequent purposes, we denote the set of critical points of
$\Delta_+$ by
\begin{equation}
\{\delta_k\}_{k\in\bbN}=\{z\in\bbC \,|\, \Delta_+^{\bullet}(z)=0 \}
\lb{3.25}
\end{equation}
and the set of critical values of $\Delta_+$ by
\begin{equation}
\{\gamma_k=\Delta_+(\delta_k)\}_{k\in\bbN}.
\lb{3.26}
\end{equation}

Next, we introduce one more family of densely defined, closed,
linear operators $H(t)$, $t\in [0, 2\pi]$ in $L^2([0,\pi])$ by
\begin{align}
& (H(t)f)(x)=(L f)(x),  \quad t\in[0,2\pi], \; x \in [0,\pi],  \lb{3.28}
\\
& f\in\dom(H(t))=\{g\in L^2([0,\pi]) \,|\, g,g'\in AC([0,\pi]); \,
L g \in L^2([0,\pi]); \no \\
& \hspace*{5.6cm} g(\pi)=e^{it}g(0), \, g'(\pi)=e^{it}g'(0)\}.
\no
\end{align}
Again, $H(t)$, $t\in [0,2\pi]$, is self-adjoint if and only if $V$ is
real-valued. The spectrum of $H(t)$ is then given by
\begin{equation}
\sigma(H(t))=\{E_k(t)\}_{k\in\bbN_0}=\{z\in\bbC
\,|\, \Delta_+(z)=\cos (t)\}, \quad t\in [0,2\pi],   \lb{3.29}
\end{equation}
and the spectrum of $H$ is given by
\begin{equation}
\sigma(H)=\bigcup_{0\leq t\leq
\pi}\sigma(H(t)). \lb{3.30}
\end{equation}

Finally, we establish the connection between the original
Schr\"odinger operator $H$ in $L^2(\bbR)$ and the family of operators
$\{H(t)\}_{t\in[0,2\pi]}$ in $L^2([0,\pi])$ using the notion of direct
integrals and the Gel'fand transform \cite{Ge50}. To this end we consider
the direct integral of Hilbert spaces
\begin{equation}
\cK=\f{1}{2\pi}\int^{\oplus}_{[0,2\pi]} dt\, L^2([0,2\pi])  \lb{3.31}
\end{equation}
with constant fibers $L^2([0,2\pi])$. Elements $F\in\cK$ are represented
by
\begin{align}
\begin{split}
& F=\big\{F(\cdot,t)\in L^2([0,\pi];dx)\big\}_{t\in[0,2\pi]} \, , \\
& \|F\|^2_{\cK} = \f{1}{2\pi}
\int_{[0,2\pi]} dt\, \|F(\cdot,t)\|^2_{L^2([0,\pi])}
= \f{1}{2\pi}
\int_{[0,2\pi]} dt\int_{[0,\pi]} dx\, |F(x,t)|^2  \lb{3.32}
\end{split}
\end{align}
with scalar product in $\cK$ defined by
\begin{align}
\begin{split}
(F,G)_{\cK} &= \f{1}{2\pi}
\int_{[0,2\pi]} dt\, (F(\cdot,t),G(\cdot,t))_{L^2([0,\pi])} \\
&= \f{1}{2\pi}
\int_{[0,2\pi]} dt \int_{[0,\pi]} dx\,
\ol {F(x,t)} G(x,t), \quad F, G \in \cK.  \lb{3.33}
\end{split}
\end{align}
Frequently, $\cK$ is written as the vector-valued Hilbert space
\begin{equation}
\cK = L^2\big([0,2\pi]; dt/(2\pi); L^2([0,\pi];dx)\big).  \lb{3.33a}
\end{equation}

The Gel'fand transform \cite{Ge50} is then defined by
\begin{equation}
\cG\colon \begin{cases} L^2(\bbR) \to \cK  \\
\hspace*{.75cm} f \mapsto (\cG f)(x,t)=F(x,t)=\slimes_{N \uparrow \infty}
\sum_{n=-N}^N  f(x+n \pi) e^{-int},   \\[1mm]
\end{cases}    \lb{3.34}  
\end{equation}
where $\slimes$ denotes the limit in $\cK$. By inspection,
$\cG$ is a unitary operator. The inverse transform is given by
\begin{equation}
{\cG}^{-1}=\begin{cases} \cK \to L^2(\bbR) \\
\hspace*{.05cm}
F \mapsto ({\cG}^{-1}F)(x+n \pi) = \f{1}{2\pi}\int_{[0,2\pi]}
dt\, F(x,t) e^{int}, \quad n\in\bbZ.  \\[1mm]
\end{cases}  \label{3.35}
\end{equation}
One then infers that
\begin{equation}
\cG \varphi(H) {\cG}^{-1} = \f{1}{2\pi}\int^{\oplus}_{[0,2\pi]} dt\,
\varphi(H(t))
\lb{3.36}
\end{equation}
for a large class of functions $\varphi$ on $\bbC$ (including polynomials
and bounded continuous functions on $\bbC$). In particular, \eqref{3.36}
applies to the resolvent $(H-z I)^{-1}$, $z\in\rho(H)$, of $H$ and one
obtains
\begin{align}
\begin{split}
& G(z,x+m\pi,y+n\pi)=\f{1}{2\pi} \int_{[0,2\pi]} dt \, G(t,z,x,y)
e^{it(m-n)}, \\
& \hspace*{3.05cm} z\in\rho(H), \quad x,y \in [0,\pi],
\quad m,n \in \bbZ, \lb{3.47}
\end{split}
\end{align}
where
\begin{equation}
G(t,z,x,y)=(H(t)-z I)^{-1}(x,y), \quad z\in \rho(H(t)), \; t\in [0,2\pi],
\;\; x,y\in [0,\pi].  \lb{3.48}
\end{equation}

As discussed by McGarvey \cite{Mc62}, \cite{Mc65}, \cite{Mc65a}, a bounded
decomposable operator $B$ commuting with translations by $\pi$ is a
spectral operator of scalar type if, roughly speaking, the associated
family $\{B(t)\}_{t\in [0,2\pi]}$ consists of spectral operators which are
uniformly spectral with respect to $t\in [0,2\pi]$. Since one can show that 
$H$ is a spectral operator of scalar type if and only if
$(H- z I)^{-1}$, $z\in\rho(H)$, is one whose spectral resolution satisfies
$E_{(H-z I)^{-1}}(\{0\})=0$, this seems to offer a reasonable strategy to
characterize those Hill operators which are spectral operators of scalar
type. However, to the best of our knowledge, this strategy has never been
pursued successfully. In this paper we proceed somewhat differently and
focus our attention directly on analyzing spectral projections
$P(\sigma)$ of $H$ in terms of the corresponding spectral projections
of $H(t)$, $t\in [0,2\pi]$.

\section{Principal Results} \lb{s4}

For every compactly supported element $g\in L^2(\R)$ and every regular
spectral arc $\sigma\subset\sigma(H)$ of $H$ we set, following \cite{Tk64},
\begin{align}\label{defproj}
&(P(\sigma) g)(x)  \no \\
&\quad=\frac{1}{2\pi}\int_{\sigma}
                     \frac{d\lambda}{\sqrt{1-\Delta_+(\lambda)^2}}
\Big[\big[\phi(\lambda,\pi)\theta(\lambda,x)
-\Delta_-(\lambda)\phi(\lambda,x)\big] F_\theta(\lambda;g) \no \\
& \hspace*{4.2cm} -\big[\theta'(\lambda,\pi)\phi(\lambda,x)
+\Delta_-(\lambda)\theta(\lambda,x)\big]F_\phi(\lambda;g)\Big],
\end{align}
where
\begin{align} \label{l}
\begin{split}
& F_\theta(\lambda;g)=\int_{\R} dy\, \theta(\lambda,y)g(y), \quad
F_\phi(\lambda;g)=\int_{\R} dy\, \phi(\lambda,y)g(y), \\
& \hspace*{2.55cm} \lambda \in \sigma(H), \quad g\in L^2(\bbR),
\, \text{ $\supp (g)$ compact,}
\end{split}
\end{align}
and the square root in \eqref{defproj} is understood to have a positive
value on $\sigma(H)$.

Using the Floquet solutions \eqref{3.20}, we can write
\eqref{defproj} in the form
\begin{equation}\label{l361}
(P(\sigma)g)(x) =\frac{1}{4\pi}
\int_{\sigma} d\lambda
\frac{\phi(\lambda,\pi)}{\sqrt{1-{\Delta}_+(\lambda)^2}}
\big[\psi_+(\lambda,x)F_-(\lambda;g)+
\psi_-(\lambda,x)F_+(\lambda;g)\big],
\end{equation}
where
\begin{equation} \label{l36}
F_{\pm}(\lambda;g)=\int_{{\mathbb R}} dy\,
\psi_{\pm}(\lambda,y)g(y),  \quad \lambda \in \sigma(H), \quad g\in L^2(\bbR),
\, \text{ $\supp (g)$ compact.}
\end{equation}

\begin{theorem}\label{t4.1}
If $\sigma$ is a regular spectral arc of $H$, then there exists a constant 
$C(\sigma)>0$ such that 
\begin{equation}\label{l362}
\int_{\sigma} |d\lambda|\frac{|\phi(\lambda,\pi)|}{\sqrt{1-\Delta_+(\lambda)^2}}\bigg(
\frac{w_+(\lambda)}{w_-(\lambda)}|F_{-}(\lambda;g)|^2+\frac{w_-(\lambda)}{w_+(\lambda)}|F_+(\lambda;g)|^2\bigg)\leq C(\sigma)^2\|g\|^2_{L^2(\R)}, 
\end{equation} 
where $|d\lambda|$ denotes the arc length measure along $\sigma$ and 
\begin{equation}\label{l363}
w_{\pm}(\lambda)=\|\psi_{\pm}(\lambda,\cdot)\|_{L^2([0,\pi])},\quad
 \lambda\in\sigma(H). 
\end{equation}
The operator  $P(\sigma)$ extends to
a bounded linear projection operator in $L^2(\bbR)$ satisfying
\begin{equation}\label{proj}
P(\sigma_1)P(\sigma_2)=P(\sigma_1\cap\sigma_2) \,
\text{ for all regular spectral arcs $\sigma_k$, $k=1,2$, of $H$.}
\end{equation}\label{sproj}
The subspace $\ran(P(\sigma))$ is invariant with respect to $H$ and
\begin{equation}
\sigma\big(H|_{\ran(P(\sigma))}\big)=\sigma.   \lb{range} 
\end{equation}
Finally, if for some $g\in L^2(\bbR)$, $P(\sigma)g=0$ for all regular spectral arcs $\sigma$ of $H$, then $g=0$.
\end{theorem}

Later we will prove that if $H$ is a spectral operator of scalar type,
then the family of projections $P(\sigma)$, $\sigma$ a regular spectral
arc of $H$, extends to a spectral resolution $E_H(\cdot)$ of $H$ such that
\begin{equation}
\|E_H(\sigma)\|_{\cB(L^2(\bbR))} \leq C, \quad \sigma \subseteq \sigma(H)
\, \text{ measurable}
\end{equation}
for some finite positive constant $C$ independent of
$\sigma\subseteq\sigma(H)$. 
(We denote by $\cB(\cH)$ the Banach space of bounded, linear operators on 
the complex, separable Hilbert space $\cH$.)

Given the notion of projections $P(\sigma)$ for regular spectral arcs
$\sigma$ of $H$ as defined in \eqref{defproj}, we can now define the notion
of spectral singularities for general Hill operators $H$.

\begin{definition} \lb{d4.1A} The point $\lambda_0\in\sigma(H)$ is called
a {\it spectral singularity of $H$} if for all sufficiently small
$\delta>0$ the relation
\begin{equation}
\sup_{\sigma\in \Sigma(\lambda_0,\delta)} \|P(\sigma)\|_{\cB(L^2(\bbR))} = \infty
\end{equation}
holds, where $ \Sigma(\lambda_0,\delta)$
is the set of all regular spectral arcs of $H$ located in the set
\begin{equation}
\sigma(H)\cap\{z\in\bbC \,|\, 0<|z-\lambda_0|<\delta\}.
\end{equation}
\end{definition}
According to Theorem \ref{t4.1} (cf.\ also Definition \ref{d3.4}), all spectral 
singularities of $H$ are contained in the discrete set
\begin{equation}
\cS_+=\{\lambda\in\sigma(H) \,|\, \Delta_+^{\bullet}(\lambda)=0\}.
\end{equation}

We note that the Riesz projection associated with a contour
surrounding a spectral singularity $\lambda_0$ is not necessarily unbounded
as cancellations may occur when several spectral arcs of $H$ end at
(or cross through) $\lambda_0$. The crucial point in Definition \ref{d4.1A}
is the blowup of the norm of the projection $P(\sigma)$ for some arc $\sigma$ 
that converges to $\lambda_0$.

\begin{lemma} \lb{l4.1a} If the Hill operator $H$ in \eqref{3.2} is a
spectral operator of scalar type,
$E_H$ is its spectral resolution, and $\sigma$ is a regular spectral arc
of $H$, then $E_H(\sigma)=P(\sigma)$.
\end{lemma}

Next, we define the Paley--Wiener class ${\mathbb{PW}_\pi}$ as the set of
all entire functions of exponential type not exceeding $\pi$ satisfying
\begin{equation}
\|f\|^2_{\mathbb{PW}_\pi}=\int_{\R} dx\, |f(x)|^2  <\infty.
\end{equation}
We note that
\begin{equation}
|f(\zeta)|\leq C \|f\|_{L^2(\bbR)}e^{\pi|\Im(\zeta)|}, \quad
\zeta\in\bbC, \; f\in\mathbb{PW}_\pi,
\end{equation}
for some $C>0$ independent of $f\in \mathbb{PW}_\pi$.

The following theorem in \cite{ST96} (see also \cite{ST96a}, \cite{ST97})
provides a spectral parametrization of the operator $H$ in \eqref{3.2}.

\begin{theorem} [\cite{ST96}, \cite{ST97}]  \label{t4.2}
A triple of entire functions $\{s=s(\zeta), u_+=u_+(\zeta),
u_-=u_-(\zeta)\}$, with $\zeta\in\bbC$, coincides with a  triple
$\{\phi(\zeta^2, \pi),
\Delta_+(\zeta^2), \Delta_-(\zeta^2)\}$ of some Hill  operator if and
only if the following conditions are satisfied:
\begin{align}
& (i) \quad \,\,
s(\zeta)=\frac{\sin(\pi \zeta)}{\zeta}-\pi\langle V \rangle
\frac{\cos(\pi \zeta)}{2\zeta^2}+\frac{g(\zeta)}{\zeta^2},
\quad \zeta\in\bbC, \; g\in{\mathbb{PW}_\pi}.\label{l601} \\
& (ii) \quad \,
u_+(\zeta)=\cos(\pi
\zeta)+\pi \langle V \rangle \frac{\sin(\pi\zeta)}{2\zeta}
-\pi^2 \langle V \rangle^2 \frac{\cos(\pi \zeta)}{8\zeta^2}+
\frac{f(\zeta)}{\zeta^2},    \label{l61} \\
& \hspace*{7.2cm}  \zeta\in\bbC, \;
f\in {\mathbb{PW}_\pi},   \no  \\
& \hspace*{.9cm} \text{with } \, \{f(n)\}_{n\in\bbZ}\in\ell^1(\bbZ).  \label{l32}  \\
& (iii) \quad
\frac{u_+(\zeta)^2-1-u_-(\zeta)^2}{s(\zeta)} \,
\text{ is an entire function in $\zeta$.}\label{l62}  \\
& (iv) \quad  \zeta u_-(\cdot)\in{\mathbb{PW}_\pi}.  \\
& (v) \quad \,
\text{The Gel'fand--Levitan equation for the transformation kernel $K(x,y)$,}
\no \\
&  \hspace*{.9cm} K(x,y)+ F(x,y)+\int_0^x ds \, K(x,s) F(s,y) =0, \quad
0 \leq y \leq x \leq \pi, \\
&  \hspace*{.9cm}  F(x,y)= [\Phi(x-y)+\Phi(x+y)]/2, \quad
0 \leq y \leq x \leq \pi,  \\
&  \hspace*{.9cm} \Phi(x)= \sum_{\ell\in\bbZ} \bigg[\res_{\zeta=m_{\ell}}
\bigg(\f{u_-(\zeta)+u_+(\zeta)}{\zeta s(\zeta)} \cos(\zeta x)\bigg)
-\f{1}{\pi} \cos(\ell x)\bigg], \quad  x\in [0,2\pi], \\
& \hspace{.9cm} \text{is uniquely solvable in $L^2([0,x])$ for all\
$x\in[0,\pi]$ } \, with \no \\
&  \hspace*{.9cm} V(x)=2\f{d}{dx} K(x,x) \, \text{ for a.e.\ $x\in [0,\pi]$.}
\end{align}
\end{theorem}

\noindent Here, $\langle V \rangle $ denotes the mean value of $V$
\begin{equation}
\langle V \rangle = \f{1}{\pi} \int_0^\pi dx\, V(x)
\end{equation}
and $\{m_\ell\}_{\ell\in\bbZ}$ denotes all zeros of $\zeta s(\zeta)$.

In the following theorems (and in the remainder of this paper), $C$ denotes 
a finite positive constant whose value varies from place to place.

Next we turn to the principal results of this paper. The fundamental
question to be answered is the following: When is $H$ a
spectral operator of scalar type (cf.\ \cite[p.\ 1938 and 2242]{DS88a} and
Appendix \ref{A})? The answer, in terms of intrinsically Floquet theoretic
terms participating in the spectral parametrization of $H$ in
Theorem \ref{t4.2}, reads as follows:

\begin{theorem}\label{t4.3}
A Hill operator $H$ is a spectral operator of scalar type if and only if
the following conditions $(i)$ and $(ii)$ are satisfied: \\
$(i)$ The function
\begin{equation}\label{l17a}
\frac{\Delta_+(z)^2-1-\Delta_-(z)^2}
{\phi(z,\pi)\Delta_+^{\bullet}(z)}\quad
\end{equation}
is analytic in an open neighborhood of $\sigma(H)$.\\
$(ii)$ The inequalities 
\begin{equation}\label{l19}
\left|
\frac{\phi(\lambda,\pi)}{\Delta_+^{\bullet}(\lambda)}\right|\leq C,\quad
\bigg| \frac{\Delta_-(\lambda)}{(\sqrt{|{\lambda}|}+1)
\Delta_+^{\bullet}(\lambda)}\bigg|\leq C, \quad
\lambda\in\sigma(H),\end{equation}
are satisfied with $C$ a finite positive constant independent of
$\lambda\in\sigma(H)$.
\end{theorem}
If both conditions \eqref{l17a} and \eqref{l19} are satisfied, and a point
$\lambda_0\in\sigma(H)$ is such that $\Delta_+^\bullet(\lambda_0)=0$, then
\begin{equation}\label{l52}
\Delta_+(\lambda_0)^2-1=\Delta_-(\lambda_0)=\Delta_+^\bullet(\lambda_0)
=0,\quad \Delta_+^{\bullet\bullet}(\lambda_0)\neq0,
\end{equation}
implying that the spectrum of a Hill operator, which is a spectral operator
of scalar type, is formed by a system of countably many, simple,
nonintersecting, analytic arcs. $($The latter may degenerate into finitely
many simple analytic arcs and a simple analytic semi-infinite arc, all
of which are nonintersecting$)$.


According to Theorem \ref{t4.2}, the operator $H$ and, in particular, the
function $\theta'(z,\pi)$ is uniquely determined by the triple
$\{\phi(z,\pi),\ \Delta_+(z),\ \Delta_-(z)\}$. If all points
$\{\mu_k\}_{k\in\bbN}$ of the Dirichlet spectrum are simple zeros of
$\phi(\cdot,\pi)$, then the function
$\Delta_-(z)$ may be uniquely recovered from its values
$\Delta_-(\mu_k)=\pm\sqrt{\Delta_+(\mu_k)^2-1}$, and such a connection
restricts the freedom in fixing $\Delta_-(z)$ to a choice of a sign for
every $k$. Next, we state an alternative version of Theorem
\ref{t4.3} in terms of the functions appearing in the
definition \eqref{defproj} of the spectral projections $P(\sigma)$.

\begin{theorem}\label{t4.4}
A Hill operator $H$ is a spectral operator of scalar type if and only if
the estimates
\begin{equation}\label{l18}
\bigg|
\frac{\phi(\lambda,\pi)}{\Delta_+^{\bullet}(\lambda)}\bigg|\leq C,\quad
\left|
\frac{\theta'(\lambda,\pi)}{(|\lambda|+1)\Delta_+^{\bullet}(\lambda)}
\right|\leq C,\quad \bigg|
\frac{\Delta_-(\lambda)}{(\sqrt{|{\lambda}|}+1)
\Delta_+^{\bullet}(\lambda)}\bigg|\leq C
\end{equation}
hold for all  $\lambda\in\sigma(H)$, with $C$ a finite positive
constant independent of $\lambda\in\sigma(H)$.

If the conditions \eqref{l18} are satisfied, then the functions
\begin{equation}\label{l16}
\frac{\phi(z,\pi)}{\Delta_+^{\bullet}(z)},\quad
\frac{\theta'(z,\pi)}{\Delta_+^{\bullet}(z)},\quad
\frac{\Delta_-(z)}{\Delta_+^{\bullet}(z)}
\end{equation}
are analytic in an open neighborhood of  $\sigma(H)$.
\end{theorem}

The following criterion involves the spectrum $\sigma(H)$ of $H$, the
Dirichlet spectrum $\{\mu_k\}_{k\in\bbN}$, the periodic
spectrum $\{E_k(0)\}_{k\in {\mathbb\N_0}}$, and the antiperiodic spectrum
$\{E_k(\pi)\}_{k\in {\mathbb\N_0}}$, and is connected with the algebraic
and geometric multiplicities of the eigenvalues in the sets
$\sigma(H(t))$. It follows from \eqref{l61}, see
\cite[Sect.\ 1.5]{Ma86} and also Lemma \ref{l5.1} below, that the union
of the periodic and antiperiodic spectra is formed by the numbers
\begin{align}
& \lambda_k^\pm=\left(k+\frac{\langle
V\rangle}{2k}+\frac{f_k^\pm}{k}\right)^2,
\; k\in\bbN, \quad
\sum_{k\in{\mathbb\N}}|f_k^\pm|^2<\infty,  \\
& \{E_k(0), E_k(\pi)\}_{k\in\bbN_0} =
\{\lambda_0^+, \lambda_k^{\pm}\}_{k\in\bbN}.
\end{align}
\begin{theorem}\label{t4.5}
A Hill operator $H$ is a spectral operator of scalar type if and only if
the following conditions $(i)$--$(iii)$ are satisfied: \\
$(i)$  Every multiple point of either the periodic or antiperiodic
spectrum of $H$ is a point of its Dirichlet spectrum. \\
$(ii)$ For all $t\in[0,2\pi]$ and all $E_k(t)\in\sigma(H(t))$,
each root function of the operator $H(t)$ associated with $E_k(t)$ is an
eigenfunction of $H(t)$. In particular, the geometric and algebraic
multiplicity of each eigenvalue $E_k(t)$ of $H(t)$ coincide. \\
$(iii)$ Let
\begin{equation}\label{l610}
\mathcal{Q}=\{k\in\N\,|\, d_k={\rm dist}(\delta_k,\sigma(H))>0\},
\end{equation}
then
\begin{equation} \label{l60}
\displaystyle
\sup_{k\in\cQ}\frac{|\lambda_k^+-\lambda_k^-|}
{{\rm dist}(\delta_k,\sigma(H))}<\infty,\quad
\sup_{k\in\cQ}\frac{|\mu_k-\lambda_k^-|}
{{\rm dist}(\delta_k,\sigma(H))}<\infty,\quad
\sup_{k\in\cQ}\frac{|\mu_k-\lambda_k^+|}
{{\rm dist}(\delta_k,\sigma(H))}<\infty.
\end{equation}
\end{theorem}

Here a root function of $H(t)$ associated with the eigenvalue $E_k(t)$
denotes any element $f$ satisfying $(H(t)-E_k(t))^mf=0$ for some $m\in\N$ 
(i.e., any element in the algebraic eigenspace of $H(t)$ corresponding to $E_k(t)$).

\begin{remark}\label{r4.1}
It is well-known (see, e.g., \cite[Sect.\ 8.3]{CL85}, \cite[Ch.\ 6]{Ea73},
\cite[Sect.\ 10.8]{In56}, \cite[Sect.\ 3.4]{Ma86},
\cite[Sect.\ XIII.16]{RS78}, \cite{Ti50}, \cite[Ch.\ XXI]{Ti58}) that if
$V\in L^2_{\loc}(\bbR)$ is real-valued, then $H$ and $H(t)$,
$t\in [0,2\pi]$, are self-adjoint operators, all relevant spectra are real,
and the interlacing conditions
\begin{equation}
\lambda_0^+\leq E_0(t)\leq \lambda_1^- \leq\mu_1\leq \lambda_1^+\leq
E_1(t) \leq \lambda_2^- \leq \mu_2\leq \lambda_2^+\leq E_2(t)\leq
\lambda_3^- \leq \mu_3\leq \cdots
\end{equation}
are satisfied. The spectrum $\sigma(H)$ of $H$ is then formed by the system
of bands
\begin{equation}
\sigma(H)=\bigcup_{k=0}^\infty\;[\lambda_k^+,\lambda_{k+1}^-]
\end{equation}
separated by spectral gaps
\begin{equation}
(\lambda_k^-,\lambda_k^+),  \quad k\in\bbN.  \lb{3.33A}
\end{equation}
For each $k\in\bbN$, the closure of the spectral gap \eqref{3.33A},
$[\lambda_k^-,\lambda_k^+]$, contains exactly one point $\mu_k$ of the
Dirichlet spectrum and one non-degenerate critical point $\delta_k$
of $\Delta_+$. The conditions of Theorem
\ref{t4.5} in the self-adjoint setting are quite transparent: Indeed, if
$\lambda_0$ is a multiple point of the periodic/antiperiodic spectra,
then for some $k\geq 1$, the interval $[\lambda_k^-,\lambda_k^+]$
collapses to the point $\lambda_0=\lambda_k^-=\mu_k=\lambda_k^+$, the
point $\mu_k$ of the Dirichlet spectrum is trapped in the collapsed gap
and condition $(i)$ is satisfied.

Being self-adjoint, no operator $H(t)$ has a root function which
is not its eigenfunction, implying property $(ii)$.

The interlacing conditions yield
\begin{align}
& |\mu_k-\lambda_k^\pm|\leq\lambda_k^+-\lambda_k^-
\underset{k\uparrow\infty}{=}\pi(\gamma_k^2-1)^{1/2}[1+o(1)], \\
& {\rm
dist}(\delta_k,\sigma(H))
=\min\{\lambda_k^+-\delta_k,\delta_k-\lambda_k^-\}
\underset{k\uparrow\infty}{=}\frac{\pi}{2} (\gamma_k^2-1)^{1/2}[1+o(1)]
\end{align}
which proves \eqref{l60} and hence property $(iii)$.

We also note that Marchenko and Ostrovskii \cite{MO75} gave a complete
description of all  sequences $\{\lambda_0^+, \lambda_k^\pm\}_{k\in\bbN}$
which are spectra of self-adjoint Hill operators and of all their Hill
discriminants $\Delta_+$. Such a description for non-self-adjoint Hill operators 
was given in \cite{Tk92} and \cite{Tk96}.
\end{remark}

\section{Spectral Properties of the Operators $H(t)$} \lb{s5}

In this section we take a closer look at the spectral properties of
$H(t)$, $t\in [0,\pi]$ defined in \eqref{3.28}. 
For simplicity, we assume
\begin{equation}
\langle V \rangle =0   \lb{mean}
\end{equation}
for the remainder of this paper. For notational convenience we will also
identify $E^{\pm}_0(t)=E_0(t)$ in the following.

\begin{lemma}\label{l5.1}
The spectrum $\sigma(H(t))=\{E_k(t)\}_{k\in\bbN_0}$ of $H(t)$,
$t\in[0,\pi]$, can be represented in the form
\begin{equation}
\sigma(H(t))=\{E_0(t), E_n^\pm(t)\}_{n\in\bbN},
\end{equation}
whose asymptotic expansion $($since $\langle V\rangle =0$$)$ is of the type
\begin{equation}
E_n^\pm(t)\underset{n\uparrow\infty}{=}(p_n^\pm(t))^2,\quad
p_n^\pm(t)\underset{n\uparrow\infty}{=}2n\pm\frac{t}{\pi}
+\frac{g_n^\pm(t)}{n},\label{l33C}
\end{equation}
where
\begin{equation}\label{l33B}
\sum_{n=1}^\infty|g_n^\pm(t)|^2\leq C<\infty,
\end{equation}
and $C>0$ is independent of $t\in [0,\pi]$.
\end{lemma}
\begin{proof} We set $z=\zeta^2$, $u_+(\zeta)=\Delta_+(z)$ and represent the equation 
\begin{equation}  \label{l31}
  u_+(\zeta)=\cos (t), \quad t\in[0,\pi],
\end{equation}
in the form 
\begin{equation} \label{l26}
2\sin\bigg(\frac{\pi \zeta-t}{2}\bigg)
\sin\bigg(\frac{\pi \zeta+t}{2}\bigg)
=\frac{f(\zeta)}{\zeta^2}.
\end{equation}

For every $n\in\bbN$, $n>1$, let $\delta_n=n^{-1}$ and 
$0<\delta_n\leq t\leq\pi-\delta_n$. Setting  
$\pi \zeta=2\pi n+t+\xi$ and assuming $|\xi|\leq n^{-1}$, one obtains
\begin{equation}
{\rm dist}(t+(\xi/2),\{\pi n\}_{n\in\N})\geq (2n)^{-1},\quad |\sin(t+(\xi/2))|\geq Cn^{-1}
\end{equation}
for some constant $C>0$ and represents \eqref{l26} in the form 
\begin{equation}
\xi-\frac{f_n^+(\xi,t)}{n}=0
\end{equation}
with $f_n^+$ given by
\begin{equation}
f_n^+(\xi,t)=n\;\frac{f(2 n+t/\pi+\xi/\pi)}{(2 n+t/\pi+\xi/\pi)^2}\;
\displaystyle\frac{\xi}{2\sin(\xi/2)\sin(t+(\xi/2))}.\label{l271}
\end{equation}
To estimate $f_n^+$ we now use well-known 
subharmonic function arguments (cf., \cite[Lecture 20, p.\ 150]{Le96}). Indeed, for 
$n\in\bbZ$ and $\tau\in\C$ such that $ |\tau|\leq 2$, let
\begin{align}
Q_n & =\{z\in\bbC \,|\, z=2n+x+iy,|x|\leq4,|y|\leq4\},     \\
D_n(\tau)&=\{z\in\bbC \,|\, |z-2n-\tau|\leq2\}.
\end{align}
For every $q\in\N$ the function $|f(\cdot)|^q$  is  subharmonic and hence the inequalities
\begin{align}
& (4\pi)^{-1} \int\!\!\!\int_{Q_n}dx \, dy \, |f(x+iy)|^q  \geq 
(4\pi)^{-1}\int\!\!\!\int_{D_n(\tau)} dx \,dy \, |f(x+iy)|^q  \no \\
& \quad  = (4\pi)^{-1}\int_0^2 r \, dr \int_0^{2\pi} d\theta \, |f(2n+\tau+re^{i\theta})|^q 
\geq|f(2n+\tau)|^q    \lb{4.11}
\end{align}
hold. Next, we set $q=2$ in \eqref{4.11}, and since $f\in\mathbb{PW}_\pi$, we obtain
\begin{equation}\label{l274}
\sum\limits_{n\in\Z}|f(2n+\tau)|^2\leq (2 \pi)^{-1}\int_{-\infty}^\infty
dx\; \!\!\!\int_{|y|\leq4} dy \, |f(x+iy)|^2 < \infty,
\end{equation}
where the right-hand side of the latter inequality is independent of the variable $\tau$ appearing in the left-hand side of \eqref{l274}. As result, for every $\varepsilon>0$, there exists $N\in\bbN$, $N>1$, such that
\begin{equation}
\sum\limits_{n=N}^\infty\ |f(2n+\tau)|^2\leq\varepsilon^2.\label{l272}
\end{equation}
It follows now that if we set
\begin{equation}
r_n^+=\max_{|\xi|\leq n^{-1}, \, \delta_n\leq t\leq\pi-\delta_n}|f_n^+(\xi,t)|, 
\end{equation}
then
\begin{equation}
\sum\limits_{n=N}^\infty\ (r_n^+)^2 \leq \varepsilon^2.
\end{equation}
Now we apply Rouch\'e's theorem to the function $\xi-{f_n^+(\xi,t)}/n$ in the disc 
$\{\xi\in\bbC\,|\,|\xi|\leq2r_n^+/n\}$ and find that there exists a solution 
$p_n^+(t)=2n+(t+\xi_n^+(t))\pi^{-1}$ of \eqref{l26}
such that 
\begin{equation}\label{l27}
|\xi_n^+(t)|\leq2\frac{r_n^+}{n}, \quad \delta_n\leq t\leq\pi-\delta_n,\quad 
\sum_{n=1}^\infty\ (r_n^+)^2<\infty.
\end{equation} 
Similar arguments prove that the substitution $\pi \zeta=2\pi n-t+\xi$
produces another solution $p_n^-(t)=2n-(t+\xi_n^-(t))\pi^{-1}$ of \eqref{l26} such that 
\begin{equation}
|\xi_n^-(t)|\leq2\frac{r_n^-}{n}, \quad \delta_n\leq t\leq\pi-\delta_n,\quad 
\sum_{n=1}^\infty\ (r_n^-)^2<\infty,\label{l28}
\end{equation}
which proves \eqref{l33C} with $g_n^\pm(t)=\pi^{-1}n \, \xi^\pm(t)$. 

Next, let $0\leq t\leq\delta_n$. If $\pi \zeta=2\pi n+\xi$, then \eqref{l26}  takes on the form
\begin{equation}
\xi^2-t^2=\frac{f_n(\xi,t)}{n^2},   \label{l275}
\end{equation}
where
\begin{equation} \label{l29}
f_n(\xi,t)=n^2\frac{f(2 n+\xi/\pi)}{(2 n+\xi/\pi)^2}\;
\displaystyle\frac{\xi-t}{2\sin((\xi-t)/2)}
\frac{\xi+t}{\sin((\xi+t)/2)}.
\end{equation}
We note that since $f'\in {\mathbb{PW}_\pi}$, the estimate \eqref{l272} 
is satisfied with $f'$ instead of 
$f$. In addition, according to condition \eqref{l32}, we have $\{f(2n)\}_{n\in\bbZ}\in\ell^1(\bbZ)$. We use \eqref{4.11} with $q=1$ and the representation
\begin{equation}
f(\zeta)=f(2n)+\int_{2n}^{\zeta} d\tau \, f'(\tau) 
\end{equation}
to arrive at the estimate
\begin{equation}
\sum_{n=N}^\infty\ |f(2n+\tau)|\leq\varepsilon,\quad |\tau|\leq 2.
\end{equation}
Therefore, the sequence
\begin{equation}
r_n=\max_{|\xi|\leq n^{-1}, \, 0\leq t\leq\delta_n}|f_n(\xi,t)|
\end{equation}
is such that 
\begin{equation}
\sum_{n=N}^\infty\ r_n \leq \varepsilon.
\end{equation}
If $|\xi|=2\delta_n$, then $|\xi^2-t^2|\geq|\xi^2|-|t^2|\geq3\delta_n^2$, 
and by Rouch\'e's theorem, for all $n\geq N$, every disc $\{\xi\in\bbC \,|\, 
|\xi|\leq 2\delta_n\}$ contains two solutions $\xi_n^\pm(t)$  of \eqref{l26}. If for a
given  $n\geq N$ we have $t\leq2n^{-1}\sqrt{r_n}$, then \eqref{l275}
implies $|\xi_n^{\pm}(t)|^2\leq5n^{-2}r_n$, which means that there exist two solutions 
\begin{equation} \label{l30}
p_n^{\pm}(t)=2 n\pm\frac{t}{\pi}
+\frac{g_n^{\pm}(t)}{n},\quad g_n^{\pm}(t)=\frac{\mp t+\xi_n^\pm(t)}{\pi},\quad |g_n^\pm(t)|\leq C n^{-1}\sqrt{r_n}
\end{equation}
of \eqref{l26}, and hence \eqref{l33B} is satisfied.
On the other hand, if $2n^{-1}\sqrt{r_n}\leq t\leq n^{-1}$, then 
\begin{equation}
\xi=\pm t\left(1+\frac{f_n(\xi,t)}{n^2t^2}\right)
^{1/2}=\pm t+\frac{g_n^\pm(\xi,t)}{n},
\end{equation}
where $|g_n^\pm(\xi,t)|\leq C\sqrt{|f_n(\xi,t)|}\leq C\sqrt{r_n}$. The same arguments apply to 
the case $\pi-\delta_n\leq t\leq\pi$ and we
again obtain \eqref{l30}, completing the proof of Lemma \ref{l5.1} for
all $t\in[0,\pi]$ .
\end{proof}

\begin{corollary} \label{c5.2}
There exists a finite positive constant $C$ such that for every sequence
$\{c_n^\pm\}_{n=0}^\infty
\in\ell^2(\bbN_0)$ and every $t\in[0,\pi]$, the series
\begin{equation}
\Theta^\pm (x,t)=\sum_{n=0}^\infty\
c_n^\pm\theta(E_n^\pm(t),x),\quad
\Phi^\pm (x,t) =\sum_{n=0}^\infty\ c_n^\pm\sqrt{|E_n^\pm (t)|+1}
\phi(E_n^\pm (t),x)
\end{equation}
converge in $L^2([0,\pi])$ and their sums satisfy the estimate
\begin{equation}
\left\|\Theta^\pm (\cdot,t)
\right\|^2_{L^2([0,\pi])}+
\left\|\Phi^\pm (\cdot,t)\right\|^2_{L^2([0,\pi])}
\leq C\sum_{n=0}^\infty\ |c_n^\pm|^2.
\end{equation}
\end{corollary}

\begin{corollary}\label{c5.3}
There exists a finite positive constant $C$ such that for every element
$f\in L^2([0,\pi])$  and every $t\in[0,\pi]$, the estimate
\begin{align}
& \sum_{n=0}^\infty
\left|
\int_0^\pi dx\, f(x)\theta(E_n^\pm(t),x) \right|^2+
\sum_{n=0}^\infty
\left|\int_0^\pi dx\, f(x)\sqrt{|E_n^\pm(t)|+1}\phi(E_n^\pm(t),x)
\right|^2   \no \\
& \quad \leq C\|f\|^2_{L^2([0,\pi])}
\end{align}
holds.
\end{corollary}
\begin{proof} Indeed, using the transformation operators
\cite[Ch.\ 1]{Ma86} for $L$ and Lemma \ref{l5.1} one obtains
\begin{align}\label{l33D}
\theta(E_n^\pm (t),x)&=\cos\bigg(\bigg(2n\pm\frac{t}{\pi}\bigg)x\bigg)+
\frac{\theta_n^\pm (x,t)}{n}, \\
\sqrt{E_n^\pm (t)}\phi(E_n^\pm (t),x)&=\sin\bigg(\bigg(2n\pm\frac{t}{\pi}
\bigg)x\bigg)+ \frac{\phi_n^\pm (x,t)}{n},\quad n\geq1,
\end{align}
where
\begin{equation}\label{l33E}
\sum_{n=0}^\infty\big(\|\theta_n^\pm (\cdot,t)\|_{L^2([0,\pi])}^2
+ \|\phi_n^\pm (\cdot,t)\|_{L^2([0,\pi]}^2\big)\leq C,\quad 0\leq
t\leq\pi,
\end{equation}
with $C$ independent of $t$. The statements of both Corollaries follow from
the corresponding properties of the trigonometric system $\{\cos
(2nx),\sin (2nx)\}_{n\in\bbN_0}$ in the space $L^2([0,\pi])$.
\end{proof}

A standard application of the Lagrange formula yields the following result
(cf.\ \eqref{3.21}).

\begin{lemma}\label{l5.4}
If $\lambda(t)\in\sigma(H(t)),\;0\leq t\leq \pi$, and
$\Delta_+^{\bullet}(\lambda(t))\ne0,\; \phi(\lambda(t),\pi)\ne0$, then
$\lambda(t)$ is a simple eigenvalue of $H(t)$ and $H(2\pi-t)$, the
corresponding eigenvectors are given by $\psi_+(\lambda(t),\cdot)$ and
$\psi_-(\lambda(t),\cdot)$, and the identity
\begin{equation}\label{l25}
\int_0^\pi dx\,
\psi_+(\lambda(t),x)\psi_-(\lambda(t),x)
=-2\,\frac{\Delta_+^{\bullet}(\lambda(t))
}{\phi(\lambda(t),\pi)}=2\,\frac{\sqrt{1-\Delta_+(\lambda(t))^2}}
{\lambda'(t)\phi(\lambda(t),\pi)}
\end{equation}
holds with the positive value of the square root understood.
\end{lemma}
The last equality in \eqref{l25} follows from differentiating
$\Delta_+(\lambda(t))=\cos(t)$ with respect to $t$.

For $t\neq0,\pi$ the boundary conditions in \eqref{3.28} are regular
(cf.\ \cite[Ch.\ V]{Na68}) and if all zeros $E_k(t)$ of the function
$\Delta_+(z)-\cos (t)$ are simple, then the spectral resolution for the
corresponding operator $H(t)$ has the form
\begin{align}
f(x)&=-\frac{1}{2}\sum_{k\in\bbN}\frac{1}{
\Delta_+^{\bullet}(E_k(t))} \no \\
& \quad \times \Bigg(\theta(E_k(t),x)\Bigg|
\begin{matrix}
-\Delta_-(E_k(t))-i\sqrt{1-\Delta_+(E_k(t))^2} & -\phi(E_k(t),\pi) \\[1mm]
\wti F_\theta(E_k(t);f) & \wti F_\phi(E_k(t);f)
\end{matrix} \Bigg|  \no \\[2mm]
& \quad \quad -\phi(E_k(t),x)\Bigg|\begin{matrix}
\theta'(E_k(t),\pi) & -\Delta_-(E_k(t))
+i\sqrt{1-\Delta_+(E_k(t))^2} \\[1mm]
\wti F_\theta(E_k(t);f) &
\wti F_\phi(E_k(t);f) \end{matrix}\Bigg| \, \Bigg), \no \\
& \hspace*{8cm} f\in L^2([0,\pi]),  \lb{5.22}
\end{align}
where $\wti F_\theta(E_k(t);f)$ and $\wti F_\phi(E_k(t);f)$ are defined by
\begin{align} \label{la}
\begin{split}
& \wti F_\theta(\lambda(t);f)=\int_{0}^{\pi} dy\, \theta(\lambda(t),y)f(y),
\quad \wti F_\phi(\lambda(t);f)=\int_{0}^{\pi} dy\, \phi(\lambda(t),y)f(y),
\\
& \hspace*{4.45cm} \lambda(t) \in \sigma(H(t)), \; t\in [0,2\pi], \quad f\in
L^2([0,\pi]).
\end{split}
\end{align}

The Floquet form of the same expansion, if $\phi(E_k(t),\pi)\neq 0$, is then of the following type
\begin{align} \label{l37}
\begin{split}
f(x)&=-\frac{1}{2}\sum_{k\in\bbN}\frac{\phi(E_k(t),\pi)}
{\Delta_+^{\bullet}(E_k(t))}
\psi_+(E_k(t),x)\int_0^\pi dy\, \psi_-(E_k(t),y)f(y), \\
& \hspace*{5.2cm} t\in [0,\pi], \quad
f\in L^2([0,\pi]).
\end{split}
\end{align}

Here $\{\psi_+(E_k(t),\cdot), \psi_-(E_k(t),\cdot)\}_{k\in\bbN}$ are a
biorthogonal system of eigenfunctions of $H(t)$ and $H(t)^*$, that is,
\begin{align}
\begin{split}
& \big(\ol{\psi_-(E_k(t),\cdot)},\psi_+(E_\ell(t),\cdot)\big)_{L^2([0,\pi])}
= \int_0^{\pi} dx\, \psi_-(E_k(t),x), \psi_+(E_\ell(t),x)  \\
& \quad = -\delta_{k,\ell}\f{2\Delta_+^{\bullet}(E_k(t))}{\phi(E_k(t),\pi)},
\quad  k, \ell \in \bbN, \quad t\in [0,2\pi]. \lb{4.37}
\end{split}
\end{align}

\section{Proofs of Theorem \ref{t4.1} and Lemma \ref{l4.1a}} \lb{s6}

We start with some preliminary considerations.

For a given regular spectral arc $\sigma\subset\sigma(H)$ of $H$, there exists an interval
$I=[\alpha,\beta]\subseteq[0,\pi]$ and an analytic function
$\lambda(\cdot)$  such that
\begin{equation}
\sigma=\{z\in\bbC \,|\, z=\lambda(t)
,\, {\Delta}_+(\lambda(t))=\cos (t), \, t\in[\alpha,\beta]\}.
\end{equation}
In addition, let $I^*=[2\pi-\beta,2\pi-\alpha]$.

Next, we 
use the functions $
w_{\pm}(\lambda)$ defined by \eqref{l363} and introduce the $2\times2$ matrix-valued functions (cf.\ \eqref{l363}) 
\begin{align}
W(\lambda)&=\begin{pmatrix}
[w_-(\lambda)/w_+(\lambda)]^{1/2} & 0\\
0& [w_+(\lambda)/w_-(\lambda)]^{1/2}
\end{pmatrix}, \quad \lambda\in\sigma(H), \\
\Omega(z)&=\begin{pmatrix} 1 & m_+(z) \\
1 & m_-(z) \end{pmatrix}, \quad z\in\bbC\backslash\{\mu_k\}_{k\in\bbN},
\end{align}
and the $\bbC^2$-vector functions
\begin{align}
{\mathbf Y}(z,x)&=\begin{pmatrix} \theta(z,x) \\ \phi(z,x)
\end{pmatrix},\; z\in\bbC, \quad {\mathbf \Psi}(z,x)=\begin{pmatrix}
\psi_+(z,x) \\ \psi_-(z,x) \end{pmatrix}, \; z\in\bbC\backslash\{\mu_k\}_{k\in\bbN}, \\
{\mathbf F}(\lambda;g)&=\begin{pmatrix}
F_\theta(\lambda;g)\\F_\phi(\lambda;g) \end{pmatrix}, \;  
\lambda\in\sigma(H),
\end{align}
where $g\in L^2(\R)$ is a compactly supported function and $F_\theta(\cdot;g)$ and 
$F_\phi(\cdot;g)$ are defined by \eqref{l}.

We use the standard scalar product in $\bbC^2$,
\begin{equation}
(\mathbf Y, \mathbf Z)_{\bbC^2}=\ol{y_1} z_1 + \ol{y_2} z_2, \quad 
\mathbf Y = \begin{pmatrix} y_1 \\ y_2 \end{pmatrix}, \; 
\mathbf Z = \begin{pmatrix} z_1\\ z_2 \end{pmatrix} \in \bbC^2,
\end{equation}
and set
\begin{equation}
J=\begin{pmatrix} 0&1\\1&0
 \end{pmatrix},\quad
P_1=\begin{pmatrix}1&0 \\0&0
 \end{pmatrix},\quad
P_2=\begin{pmatrix}0&0 \\0&1
 \end{pmatrix}.
 \end{equation} 
Given a measurable set $\sigma\subseteq\sigma(H)$, we define the Hilbert space
 $L^2(\sigma)^2$ of measurable $\bbC^2$-vector elements 
\begin{equation}
{\mathbf F}(\lambda)= \begin{pmatrix}
F_1 (\lambda) \\ F_2 (\lambda) \end{pmatrix}, 
\quad \lambda\in\sigma,
\end{equation}
with the finite norm
\begin{equation}\label{l508}
\|{\mathbf F}\|_{L^2(\sigma)^2}=
\bigg(\int_\sigma |d\lambda| \frac{|\phi(\lambda,\pi)|}
{\sqrt{1-(\Delta_+(\lambda))^2}}\|W(\lambda)\Omega(\lambda)
{\mathbf F}(\lambda)\|^2_{\C^2}\bigg)^{1/2}
\end{equation}
and scalar product
\begin{equation}
({\mathbf F},{\mathbf G})_{L^2(\sigma)^2}=\int_{\sigma} |d\lambda| 
\f{|\phi(\lambda,\pi)|}{\sqrt{1-(\Delta_+(\lambda))^2}}
(W(\lambda)\Omega(\lambda){\mathbf F}(\lambda), 
W(\lambda)\Omega(\lambda){\mathbf G}(\lambda))_{\bbC^2}.  \lb{scalarprod}
\end{equation}
We note that $|d \lambda|$ abbreviates the arc length measure along the arc 
$\sigma$ and hence the integrals in \eqref{l508} and \eqref{scalarprod} are independent of the orientation of the spectral arc $\sigma$.

\begin{lemma}
The relation
\begin{equation}\label{l503}
B_\sigma({\mathbf F,\mathbf H})
=\int_\sigma d\lambda\frac{\phi(\lambda,\pi)}
{\sqrt{1-(\Delta_+(\lambda))^2}}\big(\overline{\Omega(\lambda)
{\mathbf F}(\lambda)}, J\Omega(\lambda)
{\mathbf H}(\lambda)\big)_{\bbC^2}
\end{equation}
defines a bounded bilinear form in $L^2(\sigma)^2$ and
the estimate
\begin{equation}
|B_\sigma({\mathbf F,\mathbf H})|\leq\|{\mathbf F}\|_{L^2(\sigma)^2}
\|{\mathbf H}\|_{L^2(\sigma)^2}
\end{equation}
holds.
\end{lemma}
\begin{proof}
It suffices to note that 
\begin{align}
& W(\lambda)^{-1}J=JW(\lambda), \\
& \big(\ol{\Omega(\lambda)
{\mathbf F}(\lambda)}, J\Omega(\lambda) {\mathbf H}(\lambda)\big)_{\bbC^2}=
\big(\overline{W(\lambda)\Omega(\lambda){\mathbf F}(\lambda)}, 
JW(\lambda)\Omega(\lambda){\mathbf H}(\lambda)\big)_{\bbC^2},
\end{align}
and then to apply the Schwarz inequality.
\end{proof}

\begin{lemma}\label{l510}
Let $\sigma$ be a regular spectral arc of $H$ and introduce 
\begin{equation}\label{l511}
K_{\sigma}(H)^2= 2\pi \max_{\lambda\in\sigma, \, x\in[0,\pi]} 
\big[|\theta(\lambda,x)|^2+(|\lambda|+1) |\phi(\lambda,x)|^2\big].
\end{equation}
Then for every compactly supported function $g\in L^2(\R)$ with associated Gel'fand transform $G(x,t)$, the estimate 
\begin{equation} 
\|{\mathbf F}(\cdot;g)\|_{L^2(\sigma)^2}^2 \leq 
C(\sigma)^2 
\int_{I\cup I^*}dt\, \big\|M(\lambda(t)){\wti{\mathbf F}}
(\lambda(t);G(.,t))\big\|^2 _{\C^2} 
\label{l501}
\end{equation}
holds. Here $I = [\alpha,\beta] \subseteq [0,\pi]$, $I^* = [2\pi-\beta, 2\pi-\alpha]$, 
\begin{align}
C(\sigma)^2 & =4 K_{\sigma}(H)^2    \\
& \quad \times \sup_{\lambda\in\sigma}\Bigg[\bigg|
\frac{\phi(\lambda,\pi)}{\Delta_+^{\bullet}(\lambda)}\bigg|^2+
\left|
\frac{\theta'(\lambda,\pi)}{(|\lambda|+1)\Delta_+^{\bullet}(\lambda)}
\right|^2  
+2 \frac{[1- \Delta_+(\lambda)^2+|\Delta_-(\lambda)|^2]}{({|\lambda|}+1)|\Delta_+^{\bullet}(\lambda)|^2}\Bigg],   \no 
\end{align}
the matrix function $M(\cdot)$ is defined by 
\begin{equation}
M(\lambda)=\begin{pmatrix} 1&0\\0&(|\lambda|+1)^{1/2}  \end{pmatrix}, 
\quad \lambda \in \sigma(H),
\end{equation} 
and the vector function $\wti{\mathbf F}(\lambda;f)$ is of the form 
\begin{equation}
{\mathbf{\wti F}}(\lambda;f)=\int_0^\pi dy\,{\bf{Y}}(\lambda,y)  f(y)=
\begin{pmatrix}  {\wti F}_{\theta}(\lambda;f) \\ 
{\wti F}_{\phi}(\lambda;f) \end{pmatrix}, \quad \lambda \in \sigma(H), \; 
f\in L^2([0,\pi])
\end{equation}
with $\wti F_\theta$ and $\wti F_\phi$ defined in \eqref{la}.
\end{lemma}
\begin{proof} For most of this proof we agree to suppress the explicit $t$-dependence of  $\lambda=\lambda(t)$, $t\in I$, and for notational simplicity replace it by $\lambda$.
Since ${\bf\Psi}(z,x)=\Omega(z){\mathbf Y}(z,x)$, one obtains (we recall that 
$\Delta_+(\lambda(t))=\cos(t)$, $t\in I$)
\begin{align}
\Omega(\lambda){\mathbf{F}}(\lambda;g)&=\int_{{\mathbb R}} dy\, {\bf\Psi}(\lambda,y) g(y)= \sum_{n\in\Z}\int_0^\pi dy \begin{pmatrix}
e^{int}\psi_+(\lambda,y) \\ e^{-int}\psi_-(\lambda,y)
\end{pmatrix} g(y+n\pi) \no \\
&=P_1\Omega(\lambda) \mathbf{\wti F}(\lambda;G(\cdot,2\pi-t))+
P_2\Omega(\lambda) \mathbf{\wti F}(\lambda;G(\cdot,t)).   \label{l502}
\end{align}
Thus, one obtains, 
\begin{align}
&\|W(\lambda)\Omega(\lambda)
{\mathbf F}(\lambda;g)\|^2_{\C^2} \no \\
& \quad =\frac{w_-(\lambda)}{w_+(\lambda)}
\bigg|\int_0^\pi dy\, \psi_+(\lambda,y) G(y,2\pi-t)\bigg|^2
+\frac{w_+(\lambda)}{w_-(\lambda)} 
\bigg|\int_0^\pi dy\, \psi_-(\lambda,y) G(y,t)\bigg|^2  \no \\
& \quad = \frac{1}{w_-(\lambda)w_+(\lambda)}  \no \\
& \qquad \times \Bigg\{w_-(\lambda)^2
\bigg|\int_0^\pi dy\, \psi_+(\lambda,y) G(y,2\pi-t)\bigg|^2  \no \\
& \hspace*{1.4cm}  +w_+(\lambda)^2
\bigg|\int_0^\pi dy\, \psi_-(\lambda,y) G(y,t)\bigg|^2\Bigg\}  \no \\
& \quad = \frac{1}{w_-(\lambda)w_+(\lambda)}  \no \\
& \qquad \times \Bigg\{w_-(\lambda)^2  
\bigg|\int_0^\pi dy\, \bigg(\theta(\lambda,y) + 
\f{e^{it}-\theta(\lambda,\pi)}{(|\lambda|+1)^{1/2}\phi(\lambda,\pi)}(|\lambda|+1)^{1/2} 
\phi(\lambda,y)\bigg)  \no \\
& \hspace*{9cm} \times G(y,2\pi-t)\bigg|^2  \no \\
& \qquad +w_+(\lambda)^2
\bigg|\int_0^\pi dy\, \bigg(\theta(\lambda,y) + 
\f{e^{-it}-\theta(\lambda,\pi)}{(|\lambda|+1)^{1/2}\phi(\lambda,\pi)}(|\lambda|+1)^{1/2} 
\phi(\lambda,y)\bigg) G(y,t)\bigg|^2\Bigg\}  \no \\
& \quad = \frac{1}{w_-(\lambda)w_+(\lambda)} \Bigg\{ w_-(\lambda)^2 
\bigg| \widetilde F_{\theta}(\lambda,G(\cdot,2\pi-t))  \no \\
& \hspace*{3cm} 
+ \f{e^{it}-\theta(\lambda,\pi)}{(|\lambda|+1)^{1/2}\phi(\lambda,\pi)}(|\lambda|+1)^{1/2} 
\widetilde F_{\phi}(\lambda,G(\cdot,2\pi-t)))\bigg|^2 \no \\
& \qquad + w_+(\lambda)^2 
\bigg| \widetilde F_{\theta}(\lambda,G(\cdot,t)) 
+ \f{e^{-it}-\theta(\lambda,\pi)}{(|\lambda|+1)^{1/2}\phi(\lambda,\pi)}(|\lambda|+1)^{1/2} 
\widetilde F_{\phi}(\lambda,G(\cdot,t)))\bigg|^2\Bigg\} \no \\
& \quad \leq  \frac{2}{w_-(\lambda)w_+(\lambda)}  \no \\
& \qquad \times \Bigg\{ w_-(\lambda)^2 \bigg(1 
+ \f{|e^{it}-\theta(\lambda,\pi)|}{(|\lambda|+1)^{1/2}|\phi(\lambda,\pi)|}\bigg)^2 
\big\|M(\lambda)\widetilde{\mathbf F}(\lambda;G(\cdot,2\pi-t))\big\|^2_{\C^2}  \no \\
& \hspace*{1.35cm} + w_+(\lambda)^2 \bigg(1 
+ \f{|e^{-it}-\theta(\lambda,\pi)|}{(|\lambda|+1)^{1/2}|\phi(\lambda,\pi)|}\bigg)^2 
     \big\|M(\lambda)\widetilde{\mathbf F}(\lambda;G(\cdot,t))\big\|^2_{\C^2} \Bigg\} \no \\ 
& \quad \leq \f{2K_{\sigma}(H)^2}{w_-(\lambda)w_+(\lambda)} \bigg(1 
+ \f{|e^{it}-\theta(\lambda,\pi)|}{(|\lambda|+1)^{1/2}|\phi(\lambda,\pi)|}\bigg)^2 
\bigg(1 
+ \f{|e^{-it}-\theta(\lambda,\pi)|}{(|\lambda|+1)^{1/2}|\phi(\lambda,\pi)|}\bigg)^2 
\no \\
& \qquad \times 
 \Big[\big\|M(\lambda)\widetilde{\mathbf F}(\lambda;G(\cdot,2\pi-t))\big\|^2_{\C^2} 
+  \big\|M(\lambda)\widetilde{\mathbf F}(\lambda;G(\cdot,t))\big\|^2_{\C^2} \Big], 
     \quad \lambda \in \sigma,   \lb{5.21}
\end{align}
where we used
\begin{align}
w_{\pm}(\lambda)^2 &= \|\psi_{\pm}(\cdot,\lambda\|^2_{L^2([0,\pi])}  \no \\
&=\int_0^{\pi} dy\, \bigg|\theta(\lambda,y) + 
\f{e^{\pm it}-\theta(\lambda,\pi)}{(|\lambda|+1)^{1/2}\phi(\lambda,\pi)}(|\lambda|+1)^{1/2} 
\phi(\lambda,y)\bigg|^2  \no \\
& \leq 2 \int_0^{\pi} dy\, \Bigg[|\theta(\lambda,y)|^2 + 
\bigg|\f{e^{\pm it}-\theta(\lambda,\pi)}{(|\lambda|+1)^{1/2}\phi(\lambda,\pi)}\bigg|^2 
(|\lambda|+1) |\phi(\lambda,y)|^2\Bigg]  \no \\
& \leq 2\pi \max_{\lambda\in\sigma, \, x\in[0,\pi]}    
\big[|\theta(\lambda,x)|^2+(|\lambda|+1) |\phi(\lambda,x)|^2\big]  \no \\
& \quad \times \bigg(1+
\f{|e^{\pm it}-\theta(\lambda,\pi)|}{(|\lambda|+1)^{1/2}|\phi(\lambda,\pi)|}\bigg)^2   \no \\
& = K_{\sigma}(H)^2 \bigg(1+
\f{|e^{\pm it}-\theta(\lambda,\pi)|}{(|\lambda|+1)^{1/2}|\phi(\lambda,\pi)|}\bigg)^2, 
\quad \lambda\in\sigma, 
\end{align}
to arrive at the last inequality in \eqref{5.21}. 

Next, we recall that 
\begin{equation}
m_{\pm}(\lambda)=\f{e^{\pm it} -\theta(\lambda,\pi)}{\phi(\lambda,\pi)}, \quad 
t\in(0,\pi), \; \lambda\in\sigma(H),   \lb{5.23}
\end{equation}
and hence \eqref{3.23} implies
\begin{align}
& \bigg(1+
\f{|e^{it}-\theta(\lambda,\pi)|}{(|\lambda|+1)^{1/2}|\phi(\lambda,\pi)|}\bigg)^2
\bigg(1+
\f{|e^{-it}-\theta(\lambda,\pi)|}{(|\lambda|+1)^{1/2}|\phi(\lambda,\pi)|}\bigg)^2 \no \\
& \quad \leq \f{4}{|\phi(\lambda,\pi)|^2} \Bigg[|\phi(\lambda,\pi)|^2 
+\Bigg(\f{|e^{it}-\theta(\lambda,\pi)|^2+|e^{-it}-\theta(\lambda,\pi)|^2}{|\lambda|+1}\Bigg) + \f{|\theta'(\lambda,\pi)|^2}{(|\lambda|+1)^2}\Bigg], \no \\
& \hspace*{9.5cm} \lambda\in\sigma(H).   \lb{5.24}
\end{align}
Insertion of \eqref{5.24} into \eqref{5.21} then yields
\begin{align}
&\|W(\lambda)\Omega(\lambda)
{\mathbf F}(\lambda;g)\|^2_{\C^2} \no \\
& \quad \leq \frac{8K_{\sigma}(H)^2}{w_-(\lambda)w_+(\lambda)|\phi(\lambda,\pi)|^2}  \no \\
& \qquad \times \Bigg[|\phi(\lambda,\pi)|^2 
+\Bigg(\f{|e^{it}-\theta(\lambda,\pi)|^2+|e^{-it}-\theta(\lambda,\pi)|^2}{|\lambda|+1}\Bigg) + \f{|\theta'(\lambda,\pi)|^2}{(|\lambda|+1)^2}\Bigg]   \no \\
& \qquad \times \Big[\big\|M(\lambda)
\widetilde{\mathbf F}(\lambda;G(\cdot,2\pi-t))\big\|^2_{\C^2}+
     \big\|M(\lambda)\widetilde{\mathbf F}(\lambda;G(\cdot,t))\big\|^2_{\C^2}
\Big], \quad \lambda \in \sigma,   \lb{5.25}
\end{align}
Applying the Schwarz inequality to \eqref{l25} implies 
\begin{equation}
\f{1}{w_+(\lambda)w_-(\lambda)}\leq \f{|\phi(\lambda,\pi)|}{2 |\Delta^\bullet(\lambda)|}, 
\quad \lambda \in \sigma.   \lb{5.26}
\end{equation}
In addition, we note that by \eqref{3.20} and \eqref{5.23} one computes
\begin{equation}
|e^{it}-\theta(\lambda,\pi)|^2 + |e^{-it}-\theta(\lambda,\pi)|^2 
= 2\big[1-\Delta_+(\lambda)^2 + |\Delta_-(\lambda)|^2\big], \quad \lambda \in \sigma.  \lb{5.27}
\end{equation}
Moreover, since by \eqref{l25}
\begin{equation}
|d\lambda| \f{|\phi(\lambda,\pi||}{\sqrt{1-\Delta_+(\lambda)^2}} = 
dt \, \f{|\phi(\lambda(t),\pi)|}{|\Delta_+^{\bullet}(\lambda(t))|}, \quad t\in I,   \lb{5.28}
\end{equation}
\eqref{5.25}--\eqref{5.28} imply 
\begin{align}
& \int_{\sigma} |d\lambda|  \f{|\phi(\lambda,\pi||}{\sqrt{1-\Delta_+(\lambda)^2}} 
\|W(\lambda)\Omega(\lambda) {\mathbf F}(\lambda;g)\|^2_{\C^2}  \no \\
& \quad \leq 4 K_{\sigma}(H)^2  \sup_{\lambda\in\sigma} 
\Bigg[\f{\phi(\lambda,\pi)|^2}{|\Delta_+^{\bullet}(\lambda)|^2} 
+ \f{|e^{it}-\theta(\lambda,\pi)|^2 
+ |e^{-it}-\theta(\lambda,\pi)|^2}{(|\lambda|+1)|\Delta_+^{\bullet}(\lambda)|^2} \no \\
& \hspace*{6.8cm} 
+\f{|\theta'(\lambda,\pi)|^2}{(|\lambda|+1)^2 |\Delta_+^{\bullet}(\lambda)|^2}\Bigg] \no \\
& \qquad \times \int_{I\cup I^*} dt\, 
\big\|M(\lambda) \mathbf{ \widetilde F}(\lambda,G(\cdot,t))\big\|^2_{\bbC^2}  \no \\
& \quad = 4 K_{\sigma}(H)^2  \sup_{\lambda\in\sigma} 
\Bigg[\f{\phi(\lambda,\pi)|^2}{|\Delta_+^{\bullet}(\lambda)|^2} 
+ 2 \f{[1-\Delta_+(\lambda)^2 + |\Delta_-(\lambda)|^2]}{(|\lambda|+1)
|\Delta_+^{\bullet}(\lambda)|^2} \no \\
& \hspace*{5.7cm} 
+\f{|\theta'(\lambda,\pi)|^2}{(|\lambda|+1)^2 |\Delta_+^{\bullet}(\lambda)|^2}\Bigg] \no \\
& \qquad \times \int_{I\cup I^*} dt\, 
\big\|M(\lambda) \mathbf{ \widetilde F}(\lambda,G(\cdot,t))\big\|^2_{\bbC^2}, 
\quad \lambda\in\sigma, 
\end{align}
 implying \eqref{l501}.
\end{proof}

Using analogous arguments one finds that if $\sigma$ is a regular spectral arc of $H$, then
\begin{equation}
 \|{\mathbf Y}(\cdot,x)\|_{L^2(\sigma)^2}^2
\leq C(\sigma)^2\sup_{\mu\in\sigma}
\Big(|\theta(\mu,x)|^2+(|\mu|+1)|\phi(\mu,x)|^2\Big)<\infty
\end{equation}
with $C(\sigma)$ independent of $x\in\bbR$.
\begin{proof}[Proof of Theorem \ref{t4.1}] 
First we note that according to  \eqref{l501},
\begin{equation} 
\|{\mathbf F}(\cdot;g)\|_{L^2(\sigma)^2} \leq 2\pi C(\sigma)K(H)
\|g\|_{L^2(\R)}
\end{equation}
holds. Hence \eqref{l362} is satisfied and the mapping 

\begin{equation}
\mathbf T\colon\begin{cases} L^2(\R) \rightarrow L^2(\sigma)^2 \\
\quad \;\;\;\; g \mapsto {\mathbf F}(\lambda;g) \end{cases}   \lb{mapT}
\end{equation}
defines a bounded linear operator.
 
To show that the range of $\mathbf T$ coincides with $L^2(\sigma)^2$, we assume that 
${\mathbf H}\in L^2(\sigma)^2$ is an arbitrary element and set
\begin{align}
\begin{split}
G_{\mathbf H}(x,t)=-\frac{\phi(\lambda,\pi)}{2\Delta^\bullet(\lambda)}&
\Big[\chi(t)
\big(\overline{P_2\Omega(\lambda){\mathbf H}(\lambda)}, J\Omega(\lambda)
{\mathbf Y}(\lambda,x)\big)_{\bbC^2}  \\
&\, +\chi^*(t)\big(\overline{P_1\Omega(\lambda){\mathbf H}(\lambda)}, J\Omega(\lambda){\mathbf Y}(\lambda,x)\big)_{\bbC^2}\Big],    \label{l504}
\end{split}
\end{align}
where $\lambda=\lambda(t)=\lambda(2\pi -t)\in\sigma$, and 
$\chi(t)$ and $\chi^*(t)$ are the characteristic functions of the intervals $I$ and $I^*$,
respectively. Then
\begin{align}
|G_{\mathbf H}(x,t)| & \leq\bigg|\frac{\phi(\lambda,\pi)}{2\Delta^\bullet(\lambda)}\bigg|
\|W(\lambda)\Omega(\lambda){\mathbf H}(\lambda)\|_{\C^2}  \\
& \quad \times\bigg[\chi(t)
\left(\frac{w_-(\lambda)}{w_+(\lambda)}\right)^{1/2}|\psi_+(\lambda,x)|+
\chi^*(t)\left(\frac{w_+(\lambda)}{w_-(\lambda)}\right)^{1/2}|\psi_-(\lambda,x)|
\bigg],   \no 
\end{align}
and similar to the proof of \eqref{l501} one obtains
\begin{equation}
\|G_{\mathbf H}\|_{\mathcal K}\leq C(\sigma)
\|{\mathbf H}\|_{L^2(\sigma)^2},
\end{equation}
where $C(\sigma)$ is a parameter independent of $\mathbf H \in L^2(\sigma)^2$.
We denote by $h\in L^2(\R)$ the inverse Gel'fand transform of $G_{\mathbf H}$.

Next, we abbreviate $e_1= (1 \; 0)^\top$ and $e_2= (0 \; 1)^\top$,   
where $\top$ denotes transposition. Then, 
\begin{equation}
P_2J\Omega(\lambda)
{\mathbf Y}(\lambda,x)=\psi_+(\lambda,x)e_2, \quad 
P_2{\bf \Psi}(\lambda,x)=\psi_-(\lambda,x)e_2, 
\end{equation}
and for $t\in I$ we have
\begin{align}
& P_2\Omega(\lambda) \mathbf{\wti F}(\lambda;G_{\mathbf H}(\cdot,t))=
\int_0^\pi dx \, G_{\mathbf H}(x,t) (P_2{\bf \Psi})(\lambda,x) \\
& \quad =-\frac{\phi(\lambda,\pi)}{2\Delta^\bullet(\lambda)}
\int_0^\pi dx\; \psi_+(\lambda,x) \psi_-(\lambda,x)
\big(\ol{\Omega(\lambda){\bf H}(\lambda)},e_2\big)_{\bbC^2}e_2=
\big(\ol{\Omega(\lambda){\bf H}(\lambda)},e_2\big)_{\bbC^2}e_2. \no
\end{align}
Similarly, 
\begin{equation}
P_1\Omega(\lambda)\mathbf{\wti F}(\lambda;G_{\mathbf H}(\cdot,2\pi-t))=
\big(\ol{\Omega(\lambda){\bf H}(\lambda)},e_1\big)_{\bbC^2}e_1.
\end{equation}
Using \eqref{l502} one finds that $\Omega(\lambda){\mathbf F}(\lambda;h)=
\Omega(\lambda){\bf H}(\lambda)$ and ${\mathbf F}(\lambda;h)={\bf H}(\lambda)$.

A comparison of \eqref{l361} and \eqref{l503} then shows that the operator 
$P(\sigma)$, initially defined on the subspace of compactly supported functions
$g\in L^2(\R)$, can be extended to all $g\in L^2(\R)$ by the relation
\begin{equation}
(P(\sigma)g)(x)=\frac{1}{4\pi}B_\sigma({\mathbf F}(\cdot;g),{\mathbf Y}(\cdot,x)), 
\quad g\in L^2(\bbR).
\end{equation}
If $h\in L^2(\R)$, then
\begin{equation}
|({\overline h},P(\sigma)g)_{L^2(\bbR)}|
=\frac{1}{4\pi}|B_\sigma({\mathbf F}(\cdot;g),{\mathbf F}(\cdot;h))|\leq C(\sigma)\|g\|_{L^2(\R)}\|h\|_{L^2(\R)}.
\end{equation}
Thus, $P(\sigma)$ is a bounded operator on $L^2(\R)$. 

To prove \eqref{proj} we assume that
$\sigma_1$ and $\sigma_2$ are two regular spectral arcs of $H$:
\begin{equation}
\sigma_j=\{z\in\bbC \,|\, z=\lambda_j(t),\;\Delta_+(\lambda_j(t))
=\cos (t), \; t\in[\alpha_j,\beta_j]\subseteq[0,\pi]\}, \; j=1,2.
\end{equation}
Then for $g\in L^2(\bbR)$ with compact support,
\begin{equation}
F_+(\lambda_1(t);P(\sigma_2)g)=\Phi_+(t)+\Psi_+(t),\quad
t\in[\alpha_1,\beta_1],
\end{equation}
where
\begin{align}
\Phi_+(t)&=\frac{1}{4\pi}\;
\int_{{\mathbb R}} dy\, \psi_+(\lambda_1(t),y)
\int_{\sigma_2}
\frac{d\lambda\, \phi(\lambda,\pi)}{\sqrt{1-\Delta_+(\lambda)^2}}
\, \psi_+(\lambda,y)F_-(\lambda;g)  \\
\intertext{and}
\Psi_+(t)&=\frac{1}{4\pi}\;
\int_{{\mathbb R}} dy\, \psi_+(\lambda_1(t),y)
\int_{\sigma_2}
\frac{d\lambda\, \phi(\lambda,\pi)}{\sqrt{1-\Delta_+(\lambda)^2}}
\, \psi_-(\lambda,y)F_+(\lambda;g).
\end{align}
Substituting $\lambda=\lambda_2(s)$ in the inner integrals, using the
distributional relations
\begin{equation}\label{l41b}
\frac{1}{2\pi}\sum_{n\in\Z}\;e^{i(t-s)n}=\delta(t-s),\quad
\frac{1}{2\pi}\sum_{n\in\Z}\;e^{i(t+s)n}=0,\quad t,s\in(0,\pi),
\end{equation}
and the identity \eqref{l24}, we obtain
\begin{align}
\begin{split}
\Phi_+(t)
&=\int_{\alpha_2}^{\beta_2}
\frac{ds\, \phi(\lambda_2(s),\pi)\lambda_2'(s) }{\sqrt{1-
\Delta_+(\lambda_2(s))^2}}F_-(\lambda_2(s);g) \\
& \quad\times\frac{1}{4\pi}\sum_{n\in\Z} e^{i(t+s)n}
\int_0^\pi dx\, \psi_+(\lambda_1(t),x)\psi_+(\lambda_2(s),x)=0, 
\end{split}
\intertext{and}
\Psi_+(t)
&=\int_{\alpha_2}^{\beta_2}
\frac{ds\, \phi(\lambda_2(s),\pi)\lambda_2'(s)}{\sqrt{1-
\Delta_+(\lambda_2(s))^2}}F_+(\lambda_2(s);g) \no \\
& \quad \times\frac{1}{4\pi}\sum_{n\in\Z}
e^{i(t-s)n}
\int_0^\pi dx\, \psi_+(\lambda_1(t),x)\psi_-(\lambda_2(s),x) \no\\
& \hspace*{-.5cm} =\frac{1}{2}\chi_{\sigma_2}(\lambda_2(t))
\frac{\phi(\lambda_2(t),\pi)\lambda_2'(t)}{\sqrt{1-
\Delta_+(\lambda_2(t))^2}}F_+(\lambda_2(t);g))
\int_0^\pi dx\, \psi_+(\lambda_1(t),x)\psi_-(\lambda_2(t),x)  \no\\
&\hspace*{-.5cm} =\frac{1}{2}\chi_{\sigma_1\cap\sigma_2}(\lambda_1(t))
F_+(\lambda_1(t);g)). 
\end{align}
In the last step we used the biorthogonality of
$\{\psi_+(E_k(t),\cdot), \psi_-(E_k(t),\cdot)\}_{k\in\bbN}$ in
\eqref{5.25}. Thus, by \eqref{l25},
\begin{equation}
F_+(\lambda_1(t);P(\sigma_2)g)=\chi_{\sigma_1\cap\sigma_2}(\lambda_1(t))
F_+(\lambda_1(t);g)),\label{l41a}
\end{equation}
and in exactly the same way one obtains
\begin{equation}
F_-(\lambda_1(t);P(\sigma_2)g)=\chi_{\sigma_1\cap\sigma_2}(\lambda_1(t))
F_-(\lambda_1(t);g)).\label{l42}
\end{equation}
Hence,
\begin{align}
&(P(\sigma_1)P(\sigma_2)g)(x)=\frac{1}{4\pi}
\int_{\sigma_1}
\frac{d\lambda\,\phi(\lambda,\pi)}{\sqrt{1-{\Delta}_+(\lambda)^2}}
\big[\psi_+(\lambda,x)F_-(\lambda;P(\sigma_2)g)  \no \\
& \hspace*{6.7cm}
+ \psi_-(\lambda,y)F_+ (\lambda;P(\sigma_2)g))\big] \no \\
& \quad =\frac{1}{4\pi}
\int_{\sigma_1}
\frac{d\lambda\,\phi(\lambda,\pi)}{\sqrt{1-{\Delta}_+(\lambda)^2}}
\chi_{\sigma_2}(\lambda)
\big[\psi_+(\lambda,x)F_-(\lambda;g)+
\psi_-(\lambda,y)F_+(\lambda;g))\big] \no \\
& \quad =(P(\sigma_1\cap\sigma_2)g)(x), \quad x\in\bbR.
\end{align}

It is evident that the space $L^2(\sigma)^2$ is invariant with respect to the operator 
$\Lambda$ of multiplication by the independent variable. This implies that the space 
$\ran(P(\sigma))$ is contained in the domain of $H$ and hence it is invariant with  respect to $H$. Thus, \eqref{range} holds.

To prove the last statement of Theorem \ref{t4.1} we first note that the sequence
of roots of the function $\Delta^\bullet_+(\lambda)$ is asymptotic to
$\{k^2\}_{k=1}^\infty$ (we recall the assumption \eqref{mean}, $\langle V\rangle=0$), and therefore, the set 
\begin{equation}
{\mathcal T}=\{t\in[0,\pi]\,|\,\text{there exists}\; \lambda\in\sigma\; \text{such that } 
\Delta_+(\lambda)=\cos(t), \, \Delta^\bullet_+(\lambda)=0\}
\end{equation}
is either finite or countable, with the only possible accumulation points at $0$ and $\pi$. 

Let $I$ be a closed interval in the set $[0,\pi]\backslash{\mathcal T}$ and let $t\in I$. Then the spectrum of each operator $H(t)$ consists of simple eigenvalues $E_k(t)$, 
$k=1,2,\dots$, and each set $\sigma_k=\{\lambda\in\bbC\,|\,\lambda=E_k(t), \, t\in I\}$ is a regular spectral arc of $H$. If $P(\sigma_k)g=0$ for all such arcs, then 
$B_{\sigma_k}({\mathbf F}(\cdot;g),{\mathbf Y}(\cdot,x))=0$, $x\in\R$. Applying the operator
$H^m$ to the latter identity one obtains $B_{\sigma_k}(\Lambda^m{\mathbf F}(\cdot;g),{\mathbf Y}(\cdot,x))=0$, $x\in\R$, and since the system $\{\lambda^m\}_{m=0}^\infty$ is complete in the Hilbert space $L^2(\sigma_k)^2$, one has  
$\Omega(\lambda){\mathbf F}(\lambda;g)=0,\;\lambda\in\sigma_k$. Therefore,
\begin{align}
& F_-(E_k(t);g)=\wti{F}_{-}(E_k(t);G(\cdot,t))=
F_+(E_k(t);g)=\wti{F}_{+}(E_k(t);G(\cdot,2\pi-t))=0,   \no \\
& \hspace*{10.2cm}   t\in I, 
\end{align}
where
\begin{align}
\begin{split}
& \wti F_{\pm} (\lambda(t);f) = \int_0^{\pi} dy \, \psi_\pm (\lambda(t),y) f(y),   \\
& \lambda(t) \in \sigma(H(t)), \; t\in[0,2\pi], \quad f \in L^2([0,\pi]). 
\end{split}
\end{align}
But then the completeness of the eigensystems of the operators $H(t)$
and $H(2\pi-t)$ implies $G(x,t)=G(x,2\pi-t)=0$ for almost all $(x,t)\in [0,\pi] \times I$,  
and hence $g=0$.
\end{proof}

\begin{proof}[Proof of Lemma \ref{l4.1a}]
Let $\sigma $ be a regular spectral arc of $H$ given by \eqref{3.16b}. For
every
$ g\in L^2(\bbR)$ we set
\begin{equation}
(Q(z,\sigma)g)(x) =\frac{1}{4\pi}
\int_{\sigma}
\frac{d\lambda\,\phi(\lambda,\pi)}{\sqrt{1-{\Delta}_+(\lambda)^2}}
\big[\psi_+(\lambda,x)F_-(\lambda;g)+
\psi_-(\lambda,x)F_+(\lambda;g)\big]\frac{1}{\lambda-z}
\end{equation}
and check that $Q(z,\sigma)$ is an
analytic operator-valued function in $z\in\mathbb{C}\backslash\sigma$.
A straightforward computation shows that it coincides with the resolvent
\begin{equation}
R(z)=(H-zI)^{-1}, \quad z\in\rho(H),
\end{equation}
on the subspace $\ran(P(\sigma))$ and hence
$R(z)P(\sigma) =Q(z,\sigma)$. In addition, according to \eqref{l41a} and
\eqref{l42},  one has 
\begin{equation}
F_+(\lambda;P(\sigma)g)=F_+(\lambda;g),\quad
F_-(\lambda;P(\sigma)g)=F_-(\lambda;g),\quad \lambda\in\sigma,
\end{equation}
and the identity
\begin{equation}
Q(z,\sigma)P(\sigma)=Q(z,\sigma)\label{l43}
\end{equation}
holds.

To describe the analytic behavior of $R(z)$ with $z$ approaching the spectral
arc $\sigma$,  we use the representation \eqref{3.47}. If $z$ is sufficiently close to
$\sigma$, then
\begin{align}
\begin{split}
& G(z,x+m\pi,y+n\pi)=
\frac{1}{2\pi}\int_{[\alpha,\beta]\cup[2\pi-\beta,2\pi-\alpha]}dt \\
& \quad \times \left(
\frac{\psi_+(\lambda(t),x)\psi_-(\lambda(t),y)}{\lambda(t)-z}\;
\frac{\phi(\lambda(t),\pi)}{2\Delta_+^{\bullet}(\lambda(t))}+\Gamma(z,x,y,t)
\right)e^{it(m-n)},
\end{split}
\end{align}
where the remainder $\Gamma(z,x,y,t)$ is analytic at interior points 
$z\in\sigma$ for all $x,y\in[0,\pi]$ and $t\in[\alpha,\beta]\cup[2\pi-\beta,2\pi-\alpha]$. Equations
\eqref{l502} then show that the difference
$R(z)-Q(z,\sigma)$ is also analytic at the same points. Using \eqref{l43}
we obtain
\begin{equation}
R(z)(I-P(\sigma))=R(z)-R(z)P(\sigma)^2=R(z)-Q(z,\sigma),
\end{equation}
proving that the operator-valued function $R(z)(I-P(\sigma))$ is analytic
at interior points of $\sigma$.

Next, assume that $H$ is a spectral operator of scalar type and $E_H$
is its spectral resolution. Since  $R(z)$ and $E_H(\sigma)$ commute for
all regular arcs $\sigma$, the operator-valued functions
$R(z)E_H(\sigma)$ and $R(z)(I-E_H(\sigma))$ are analytic at interior
points of $\sigma(H)\backslash\sigma$ and $\sigma$, respectively.
Therefore, the equations
\begin{align}
\Gamma_1(z)&=R(z)E_H(\sigma)(I-P(\sigma))=E_H(\sigma)R(z)(I-P(\sigma)),
\\
\Gamma_2(z)&=R(z)(I-E_H(\sigma))P(\sigma)=(I-E_H(\sigma))R(z)P(\sigma)
\end{align}
define operator-valued functions analytic in $\mathbb{C}\backslash\{
\lambda(\alpha),\lambda(\beta)\}$. It follows from \eqref{3.20} and
\eqref{3.47} that
\begin{equation}
\|\Gamma_j(z)\|_{\cB(L^2(\bbR))} \leq C(|z-\lambda(\alpha)||z-\lambda(\beta)|)^{-\gamma}
\end{equation}
in some neighborhoods of $\lambda(\alpha)$ and $\lambda(\beta)$ with
$\gamma>0,\; C>0$ independent of $z$. In addition,
\begin{equation}
\lim_{M\uparrow\infty}\bigg[\sup_{|\Im z|\geq M}
\|\Gamma_j(z)\|_{\cB(L^2(\bbR))} \bigg]=0,
\end{equation}
and hence $\Gamma_j$, $j=1,2$, are rational operator-valued functions with
possible poles at $\lambda(\alpha)$ and $\lambda(\beta)$. If one of these
points is a pole of $\Gamma_j$, then it is an eigenvalue of $H$ which
contradicts Theorem \ref{t3.2}. Hence $\Gamma_j$ are entire functions
and by the Liouville theorem they are both equal to the zero operator
implying
\begin{equation}
E_H(\sigma)(I-P(\sigma))=(I-E_H(\sigma))P(\sigma)=0.
\end{equation}
Finally,
\begin{equation}
E_H(\sigma)-P(\sigma)=E_H(\sigma)(I-P(\sigma))-(I-E_H(\sigma))P(\sigma)
=0,
\end{equation}
completing the proof of Lemma \ref{l4.1a}.
\end{proof}

\section{Necessary Conditions for a Hill Operator to be \\ a Spectral
Operator of Scalar Type}
\lb{s7}

Let $H$ be a spectral operator of scalar type and let $E_H$ be the
corresponding resolution of the identity. Then, according to Lemma
\ref{l4.1a}, for every regular spectral arc $\sigma$ of $H$ we have
$E_H(\sigma)=P(\sigma)$ and hence there exists
a constant $C>0$ such that
\begin{equation}
\sup_{\sigma\in\sigma(H)}\|P(\sigma)\|_{\cB(L^2(\bbR))} \leq C.  \lb{6.1}
\end{equation}

The proof of the necessity of the conditions of Theorem \ref{t4.4} is
contained in the following two lemmas.

\begin{lemma}\label{l7.1}
If $H$  is a spectral operator of scalar type and
$\Delta_+^{\bullet}(\lambda_0)=0,\;
\lambda_0\in\sigma(H)$, then $\lambda_0$ is a
simple root of  $\Delta_+^{\bullet}$, that is,
$\Delta_+^{\bullet\bullet}(\lambda_0)\neq 0$. Moreover,
\begin{equation}\label{l48}
\phi(\lambda_0,\pi)=\theta'(\lambda_0,\pi)=\Delta_-(\lambda_0)=0,
\end{equation}
and
\begin{equation}\label{l49}
{\Delta}_+(\lambda_0)^2=1.
\end{equation}
\end{lemma}
\begin{proof} Let $\sigma$ be a spectral arc of $H$ in $\sigma(H)$ with
$\lambda_0$ one of its endpoints, let $[a,c]\subseteq[0,\pi]$ be a closed
interval and let $\lambda(\cdot)$ be the function on $[a,c]$ such that
\begin{equation}\label{l34}
\sigma=\{z\in\bbC \,|\,
z=\lambda(t),\, {\Delta}_+(\lambda(t))=\cos (t),\, t\in[a,c]\}.
\end{equation}

We start with some preparatory constructions. To this end we first assume
that $\lambda_0=\lambda(a)$. Since $\phi(\lambda_0,x)\not\equiv0$, we can
choose an element $Q \in L^2([0,\pi])$ such that
\begin{equation}
\int_0^\pi dx\, \phi(\lambda_0,x) Q(x) \ne 0,
\end{equation}
and find $b\in[a,c]$ such that for every interval $I=[\alpha,\beta]
\subseteq[a,b]$, the condition
\begin{equation}
\int_0^\pi dx\, \phi(\lambda(t),x) Q(x) \ne0,\quad t\in I,
\end{equation}
is satisfied. Furthermore, we set $I^*=[2\pi-\beta,2\pi-\alpha]$, 
and for every $F_0\in L^2([0,\pi])$ satisfying
\begin{equation}\label{l22}
\int_0^\pi dx\, \phi(\lambda_0,x)F_0(x) =0,
\end{equation}
we define
\begin{equation}\label{l0}
F(x,t)= \begin{cases}
F_0(x)+\rho(t) Q(x),&t\in I\cup I^*,\\
0,&{\rm otherwise}
\end{cases}
\end{equation}
with
\begin{equation}
\displaystyle
\rho(t)=-\;\frac{\int_0^\pi dx\, \phi(\lambda(t),x)F_0(x)}
{\int_0^\pi dx\,\phi(\lambda(t),x) Q(x)},\quad t\in I\cup I^*.
\end{equation}
Then $\rho(a+0)=0$,
\begin{equation}\label{l1}
\int_0^\pi dx\,\phi(\lambda(t),x)F(x,t)=0,\quad t\in[0,2\pi],
\end{equation}
the function
\begin{equation}
\Theta(t)=\int_0^\pi dx\, \theta(\lambda(t),x)F(x,t),\quad
t\in[0,2\pi],
\end{equation}
is continuous in $I\cup I^*$, and there exists the finite limit
\begin{equation}
{\Theta_0=\lim_{\alpha\downarrow a}\Theta(\alpha)=
\int_0^\pi dx\, \theta(\lambda_0,x)F_0(x).}\label{l47}
\end{equation}
In addition,
\begin{equation}
\int_0^\pi dx\,\int_0^{2\pi} dt\, |F(x,t)|^2 \leq C|\beta-\alpha|
\end{equation}
with $C>0$ independent of $\alpha$ and $\beta$. If $f$ is the inverse
Gel'fand transform of $F(x,t)$, then
\begin{equation}\label{l3}
\|f\|^2_{L^2(\R)}\leq C|\beta-\alpha|,\quad
|(\ol{f},P(\sigma_{\alpha, \beta})f)_{L^2(\bbR)}|\leq C|\beta-\alpha|,
\end{equation}
where $\sigma_{\alpha, \beta}\subseteq\sigma$ is the spectral arc of $H$
in $\sigma$ with endpoints $\lambda(\alpha)$ and $\lambda(\beta)$, and
$\ol{f}$ denotes the complex conjugate of $f$. According to \eqref{l1}
and \eqref{defproj} we have
\begin{equation}
F_\pm(\lambda(t);f)=\int_{0}^\pi dy\, \psi_\pm(\lambda(t),y)F(y,t)
=\Theta(t).
\end{equation}
Thus,
\begin{equation}\label{l4}
(\ol{f},P(\sigma_{\alpha, \beta})f)_{L^2(\bbR)}=-
\frac{1}{2\pi}\;\int_\alpha^\beta dt\,
\frac{\phi(\lambda(t),\pi)}{\Delta_+^{\bullet}(\lambda(t))}\;\Theta (t)^2.
\end  {equation}

Similarly, we can choose an element $S \in L^2([0,\pi])$
such that
\begin{equation}
\int_0^\pi dx\, \theta(\lambda_0,x)S(x)\ne0,
\end{equation}
and find $b\in[a,c]$ such that for every interval $I=[\alpha,\beta]\subseteq [a,b]$, 
the condition
\begin{equation}
\int_0^\pi dx\, \theta(\lambda(t),x)S(x) \ne0,\quad t\in I,
\end{equation}
is satisfied. Furthermore, for every element $G_0 \in L^2([0,\pi])$
satisfying
\begin{equation}\label{l8}
\int_0^\pi dx\, \theta(\lambda_0,x)G_0(x)=0,
\end{equation}
we set
\begin{equation}
G(x,t)=
\begin{cases}
G_0(x)+\widetilde{\rho}(t)S(x),&t\in I\cup I^*,\\
0,&{\rm otherwise}
\end{cases}
\end{equation}
with
\begin{equation}
\displaystyle
\widetilde{\rho}(t)=-\;\frac{\int_0^\pi dx\, \theta(\lambda(t),x)G_0(x)}
{\int_0^\pi dx\, \theta(\lambda(t),x)S(x)},\quad t\in I\cup I^*.
\end{equation}
Then $\tilde \rho(a+0)=0$, 
\begin{equation}\label{l2}
\int_0^\pi dx\, \theta(\lambda(t),x)G(x,t)=0,\quad t\in[0,2\pi],
\end{equation}
and the function
\begin{equation}
\Phi(t)=\int_0^\pi dx\, \phi(\lambda(t),x)G(x,t),\quad t\in[0,2\pi],
\end{equation}
is continuous in $I\cup I^*$. In addition,
\begin{equation}
\int_0^\pi dx\, \int_0^{2\pi} dt\, |G(x,t)|^2 \leq C|\beta-\alpha|
\end{equation}
with $C>0$ independent of $\alpha$ and $\beta$. If $g$ is the inverse
Gel'fand transform of $G(x,t)$ then
\begin{equation}\label{l5}
\|g\|^2_{L^2(\R)}\leq C|\beta-\alpha|,\quad
|(\ol{g},P(\sigma_{\alpha, \beta})g)_{L^2(\bbR)}|\leq C|\beta-\alpha|. 
\end{equation}
According to \eqref{l2} and \eqref{3.20} we infer  
\begin{equation}
F_\pm(\lambda(t);g)=\int_{0}^\pi dx \, \psi_\pm(\lambda(t),x)G(x,t)
= m_\pm(\lambda(t)) \Phi(t).  \lb{7.23}
\end{equation}
Thus,
\begin{equation}\label{l11}
(\ol{g},P(\sigma_{\alpha, \beta})g)_{L^2(\bbR)}
= \frac{1}{2\pi}\;\int_\alpha^\beta dt\,
\frac{\theta'(\lambda(t)),\pi)}{\Delta_+^{\bullet}(\lambda(t))}\;
\Phi(t)^2.
\end{equation}

Since $\Delta_+^{\bullet}(\lambda_0)=0$, we have
\begin{equation}
{\Delta}_+(\lambda)={\Delta}_+(\lambda_0)+c_0(\lambda-\lambda_0)^k[1+o(1)]
\lb{7.26}
\end{equation}
with $c_0\neq0$ and $k\geq2$.

After these preparations we now turn to the proof of Lemma \ref{l7.1}.

If we assume that $\Delta_+(\lambda_0)^2\ne1$, then $a\ne0,\;\sin (a) \ne0$.
Together with the expansion
\begin{equation}
\cos (t)=\cos (a) -(t-a)\sin (a) [1+o(1)],
\end{equation}
valid for $t$ sufficiently close to $a$, this yields
\begin{equation}\label{l9}
\lambda (t)=\lambda_0+c_1(t-a)^{1/k}[1+o(1)],\quad
\Delta_+^{\bullet}(\lambda(t)) =c_2(t-a)^{1-(1/k)}[1+o(1)]
\end{equation}
with some non-vanishing constants  $c_1, c_2 \in\bbC$.

In addition to $\Delta_+(\lambda_0)^2\ne1$, we now assume  
$\phi(\pi,\lambda_0)\ne0$. For all $F_0\in L^2([0,\pi])$ satisfying
\eqref{l22}, we define
$F(x,t)$ by \eqref{l0} and note that
if for all such $F_0$ the constant $\Theta_0$
defined by \eqref{l47} is zero, then the Hahn--Banach theorem implies
$\theta(\lambda_0,x)=C\phi(\lambda_0,x)$, which represents a contradiction.
Hence we can assume $F_0$ to be such that $\Theta_0\neq0$.
The inverse Gel'fand transform $f$ of
$F(x,t)$ then satisfies \eqref{l3}. Using \eqref{l4} we find
\begin{equation}
\bigg|\int_\alpha^\beta dt\, (\phi(\lambda_0,\pi)\Theta_0^2+\varepsilon(t))
{(t-a)^{-1+(1/k)}} \bigg|\leq C|\beta-\alpha|, 
\end{equation}
where $\varepsilon(t)\downarrow 0$ uniformly with respect to $t\in[\alpha,\beta]$
as $\beta\downarrow a$. Since $\phi(\lambda_0,\pi)\neq0$, we obtain
\begin{equation}\label{l6}
k((\beta-a)^{1/k}-(\alpha-a)^{1/k})=\bigg|\int_\alpha^\beta dt\,
{(t-a)^{-1+(1/k)}}
\bigg|\leq C|\beta-\alpha|.
\end{equation}
Since $k\geq2$ and $\beta$ are close to $a$, this is impossible.
Hence, $\phi(\lambda_0,\pi)=0$.

If, as before, $\Delta_+(\lambda_0)^2\ne1$ and, in addition,
$\theta'(\lambda_0,\pi) \ne0$, then \eqref{l5} yields
\begin{equation}\label{l13}
\bigg|
\int_\alpha^\beta dt\,(\theta'(\lambda_0,\pi)\Phi_0^2+o(1))
{(t-a)^{-1+(1/k)}} \bigg|\leq C|\beta-\alpha|, 
\end{equation}
where
\begin{equation}
\Phi_0=\lim_{\alpha\downarrow a}\Phi(\alpha)
=\int_0^\pi dx\, \phi(\lambda_0,x) G_0(x)
\end{equation}
for every element $G_0\in L^2([0,\pi])$ satisfying \eqref{l8}. Using again
the Hahn--Banach arguments we conclude that $\Phi_0\neq0$, which again
leads to a contradiction. Therefore, $\theta'(\lambda_0,\pi)=0$. Thus,
\begin{equation}\label{l7}
\frac{\phi(\lambda(t),\pi)}{\Delta_+^{\bullet}(\lambda(t))}
=O\big(|t-a|^{-1+(2/k)}\big),
\quad
\frac{\theta'(\lambda(t),\pi)}{\Delta_+^{\bullet}(\lambda(t))}
=O\big(|t-a|^{-1+(2/k)}\big).
\end{equation}

Let $R \in L^2([0,\pi])$ such that
\begin{equation}
X(t)=\int_0^\pi dx\, \theta(\lambda(t),x)R(x) \ne0,\quad
\Psi(t)=\int_0^\pi dx\, \phi(\lambda(t),x)R(x) \ne0,
\end{equation}
for all $t\in[a,\beta]$ with $\beta$ sufficiently close to $a$, and let
\begin{equation}
R(x,t)=\begin{cases}
R(x), &t\in I\cup I^*,\\
0, &{\rm otherwise}. \end{cases}
\end{equation}
If $r$ is the inverse Gel'fand transform of $R(x,t)$, then
\begin{equation}
F_\pm(\lambda(t);r)=\begin{cases}
X(t)+m_\pm(\lambda(t))\Psi(t), &t\in I\cup I^*,\\
0, &{\rm otherwise},
\end{cases}
\end{equation}
and
\begin{align}
\begin{split}
(\ol{r},P(\sigma_{\alpha, \beta})r)_{L^2(\bbR)}
&=\frac{1}{2\pi}\;\int_\alpha^\beta dt\,
\left[\frac{\phi(\lambda(t),\pi)}{\Delta_+^{\bullet}(\lambda(t))}\,X(t)^2
                         -
\frac{\theta'(\lambda(t),\pi)}{\Delta_+^{\bullet}(\lambda(t))}\,
\Psi(t)^2\right]   \\
& \quad -\frac{1}{\pi}\;\int_\alpha^\beta dt\,
\frac{\Delta_-(\lambda(t))}{\Delta_+^{\bullet}(\lambda(t))}\,X(t)
\Psi(t).
\end{split}
\end{align}
Because of \eqref{l7}, the absolute value of the former integral is not
greater than $C((\beta-a)^{2/k}-(\alpha-a)^{2/k})$. On the other hand, 
$\|r\|_{L^2(\R)}^2=\|R\|_\mathcal K^2\leq C|\beta-\alpha|$, and hence,
\begin{equation}\label{l12}
\bigg|
\displaystyle\int_\alpha^\beta dt\,
\frac{\Delta_-(\lambda(t))}{\Delta_+^{\bullet}(\lambda(t))}\,X(t)\Psi(t)
\bigg|\leq C((\beta-a)^{2/k}-(\alpha-a)^{2/k}).
\end{equation}
Since $X(a)\Psi(a)\neq0$ and \eqref{l9} is valid, we conclude that
$\Delta_-(\lambda_0)=0$. Together with $\theta'(\pi,\lambda_0)=
\phi(\pi,\lambda_0)=0$ this implies $\Delta_+(\lambda_0)^2=1$ contradicting the
assumption $\Delta_+(\lambda_0)^2\neq1$.

Repeating the above arguments for the case
$\lambda_0=\lambda(c)$ we find that the condition
$\Delta_+^{\bullet}(\lambda)=0$ for $\lambda\in\sigma(H)$ is only
satisfied for
\begin{equation}
\Delta_+(\lambda_0)^2=1.
\end{equation}
Next, we will show that if this is the case, then $k=2$ in \eqref{7.26}. For simplicity 
we restrict ourselves to the case $a=0$, $\lambda_0=\lambda(0)$
and $\Delta_+(\lambda_0)=1$.

First we note that instead of \eqref{l9} we have
\begin{equation}\label{l10}
\lambda (t) =\lambda_0+c_1t^{2/k}[1+o(1)], \quad
\Delta_+^{\bullet}(\lambda(t))= c_2t^{2-(2/k)}[1+o(1)]
\end{equation}
for some $c_1, c_2 \in \bbC\backslash\{0\}$. Assume for a moment that $k\geq3$. Then if $\phi(\lambda_0,\pi)\neq0$ then according to \eqref{l3} and
\eqref{l4} we now obtain, in place of \eqref{l6},
\begin{equation}
\bigg|
\phi(\lambda_0,\pi) \, \Theta_0^2\int_\alpha^\beta dt\,
{t^{-2+(2/k)}}
\bigg|\leq C|\beta-\alpha|  \lb{7.41}
\end{equation}
for some constant $\Theta_0\neq0$. Hence, 
\begin{equation}
\big|\phi(\lambda_0,\pi)
\big(\beta^{-1+(2/k)}-\alpha^{-1+(2/k)}\big)\big|
\leq C|\beta-\alpha|, 
\end{equation}
which is impossible for $k\geq3$. Therefore, as before,  $\phi(\lambda_0,\pi)=0$ and
$\phi(\lambda,\pi)= A(\lambda-\lambda_0)^p[1+o(1)]$ for some integer
$p\in\bbN$ and $A\neq0$. In the same manner we obtain
$\theta'(\lambda,\pi)=B(\lambda-\lambda_0)^q[1+o(1)]$  and
$\Delta_-(\lambda)=C(\lambda-\lambda_0)^r[1+o(1)]$ with integers $q\geq
1$, $r\geq 1$. Using again \eqref{l3}, \eqref{l4}, \eqref{l5}, \eqref{l11}
and\eqref{l12} we deduce $p\geq k, q\geq k, r\geq k$. The identity
\begin{equation}\label{l21}
\Delta_+(\lambda)^2-1-\Delta_-(\lambda)^2
=-\phi(\lambda,\pi)\theta'(\lambda,\pi)
\end{equation}
then leads to the relation
\begin{equation}\label{l64}
D(\lambda-\lambda_0)^k
=C^2(\lambda-\lambda_0)^{2r}[1+o(1)]-AB(\lambda-\lambda_0)^{p+q}[1+o(1)]
\end{equation}
for some $D\neq0$, which is impossible for $k\geq3$. Thus, $k=2$ in
\eqref{7.26}.

Assuming $\phi(\lambda_0,\pi)\neq0$ and inserting $k=2$ in \eqref{7.41} yields 
\begin{equation}
|\phi(\lambda_0,\pi) \, \Theta_0^2 \, \ln(\beta/\alpha)|
= \bigg|\phi(\lambda_0,\pi) \, \Theta^2_0 \, 
\ln\bigg(1+\f{\beta-\alpha}{\alpha}\bigg)\bigg|
\leq C|\beta-\alpha|.  \lb{7.45}
\end{equation}
Since $\Theta_0\neq 0$
and $\alpha$ can be chosen arbitrarily close to
zero, \eqref{7.45} immediately leads to a contradiction. Therefore,
\begin{equation}
\phi(\lambda_0,\pi)=0.
\end{equation}
Analogous considerations yield
\begin{equation}
\theta'(\lambda_0,\pi)=0 \, \text{ and hence } \, \Delta_-(\lambda_0)=0,
\end{equation}
completing the proof of Lemma \ref{l7.1}.
\end{proof}

\begin{remark}\label{r7.2}
If for some $\lambda_0\in\bbC$,
\begin{equation}
\phi(\lambda_0,\pi)=\theta'(\lambda_0,\pi)=\Delta_-(\lambda_0)=0,
\end{equation}
then the Dirichlet and Neumann operators $H^D$ and $H^N$ (cf.\
\eqref{3.18}, \eqref{3.19a}) and, in fact, all operators $H^{\alpha}$,
$\alpha\in\bbR$ (cf.\ \eqref{3.19b}), have the eigenvalue $\lambda_0$. This
is the familiar situation in the self-adjoint case when
$\Delta^{\bullet}(\lambda_0)=0$ and $\Delta(\lambda_0)^2=1$. In this case
a spectral gap closes at $\lambda_0$ and an eigenvalue for all separated
self-adjoint boundary conditions associated with $L$ restricted to
$[0,\pi]$ is trapped at the position $\lambda_0$. Consequently, no
self-adjoint Hill operator can have spectral singularities.
Of course, this is clear from the general fact that all self-adjoint
(respectively, normal) operators in a complex separable Hilbert space are
of course spectral operators of scalar type.
\end{remark}

\begin{lemma}\label{l7.3}
If $H$  is a spectral operator of scalar type
then the functions
\begin{equation}\label{l16a}
\frac{\phi(z,\pi)}{\Delta_+^{\bullet}(z)},\quad
\frac{\theta'(z,\pi)}{\Delta_+^{\bullet}(z)},\quad
\frac{\Delta_-(z)}{\Delta_+^{\bullet}(z)}
\end{equation}
are analytic in an open neighborhood of
$\sigma(H)$ and the functions
\begin{equation}\label{l14}
\frac{\phi(\lambda,\pi)}{\Delta_+^{\bullet}(\lambda)}\;,\quad
\frac{\theta'(\lambda,\pi)}{\lambda\Delta_+^{\bullet}(\lambda)}\;,\quad
\frac{\Delta_-(\lambda)}{\sqrt{\lambda}\Delta_+^{\bullet}(\lambda)}
\end{equation}
are bounded on the set
\begin{equation}
\sigma_{{\rm ext}, R}(H)=\sigma(H)\cap \{z\in\bbC \,|\, |z|\geq R\}
\end{equation}
with $R>0$ sufficiently large.
\end{lemma}
\begin{proof} The analyticity of the functions \eqref{l16a} in an open
neighborhood of $\sigma(H)$ is an immediate consequence of Lemma
\ref{l7.1}. To prove the remaining part of Lemma \ref{l7.3} we now fix a
constant $R>0$ sufficiently large for all expressions of the form
$1+o(1)$ to be close to $1$ uniformly with respect to the variables
$t,n,x$ on which they depend.

Let $\sigma\subset\sigma(H)$ be a closed spectral arc of $H$. Without loss
of generality we can assume that the latter is a part of a spectral 
arc
\begin{equation}
\sigma_n^+=\{z\in\bbC \,|\, z=E_n^+(t),\;\Delta_+(E_n^+(t))=\cos (t),
\; t\in[0,\pi]\}, \; n\in\bbN,
\end{equation}
where $n \geq N_0$ is sufficiently large and
\begin{equation}
E_n^+(t)=(p_n^+(t))^2,\quad p_n^+(t)=2 n+\frac{t}{\pi}
+\frac{g_n^+(t)}{n}, \quad |g_n^+(t)|\leq C
\end{equation}
with $C>0$ independent of $t\in [0,\pi]$ and $n \geq N_0$.

We denote by $T=T(\sigma)\subseteq [0,\pi]$ the pre-image of $\sigma$
with respect to the  function $E_n^+(t)$ and set $T^*=\{t^*\in
[\pi,2\pi] \,|\, t^*=2\pi -t, t\in T\}\subseteq [\pi,2\pi]$. Since
\begin{align}
\theta(z,x)&=\cos (\sqrt{z}x) +\int_0^x dy\, K_\theta(x,y)
\cos (\sqrt{z}y),\\
\sqrt{z}\phi(z,x)&=\sin (\sqrt{z}x) +\int_0^x dy\, K_\phi(x,y)
\sin (\sqrt{z}y),
\end{align}
where $K_\theta(x,y)$ and $K_\phi(x,y)$ are integral kernels of transformation operators,
(cf.\ \cite[Ch.\ 1]{Ma86}), we obtain
\begin{align}
\theta(E_n^+(t),x)&=\cos\left(\left(2n+\frac{t}{\pi}\right)x\right)+\frac
{\varepsilon_{\theta,n}^+(t)}{n}, \\
\sqrt{E_n^+(t)}\phi(E_n^+(t),x)&=\sin\left(\left(2n+\frac{t}{\pi}\right)
x\right)+\frac{\varepsilon_{\phi,n}^+(t)}{n},
\end{align}
with
\begin{equation}
|\varepsilon_{\theta,n}^+(t)|+|\varepsilon_{\phi,n}^+(t)|\leq C
\end{equation}
and $C>0$ independent of $t\in [0,\pi]$ and $n\geq N_0$.

Next, set
\begin{equation}
\alpha_n(t)=
\begin{cases}
\displaystyle
-\frac{\big(\ol{\theta(E_n^+(t),\cdot)},
\phi(E_n^+(t),\cdot)\big)_{L^2([0,\pi])}}
{\sqrt{E_n^+(t)}
\big(\ol{\phi(E_n^+(t),\cdot)},\phi(E_n^+(t),\cdot)\big)_{L^2([0,\pi])}},
& t\in T\cup T^*,\\ 0, &  \text{otherwise.}
\end{cases}
\end{equation}
The above asymptotic representations show that
an inequality $|\alpha_n(t)|\leq C$  is satisfied with $C>0$ independent
of $t\in [0,\pi]$ and $n\geq N_0$. For the function
\begin{equation}
F_n(x,t)=\begin{cases}
\theta(E_n^+(t),x)+\alpha_n(t)\sqrt{E_n^+(t)}\phi(E_n^+(t),x),
&t\in T\cup T^*,\\
0,& \text{otherwise,}
\end{cases}
\end{equation}
we obtain
\begin{equation}
\int_0^\pi dx\, F_n(x,t)\phi(E_n^+(t),x)=0, \quad t\in[0,2\pi],
\end{equation}
and
\begin{equation}
\int_0^{\pi}\int_0^{2\pi} dx\, dt\, |F_n(x,t)|^2  \leq C|T|,
\end{equation}
where $|\cdot|$ abbreviates Lebesgue measure.

If $f_n\in L^2(\mathbb{R})$ is the inverse Gel'fand transform of
$F_n(x,t)$, then $\|f_n\|_{L^2(\bbR)}^2\leq C |T|$ and
\begin{equation}
(\ol{f_n},P(\sigma) f_n)_{L^2(\bbR)}
=-\frac{1}{2\pi}\int_{T} dt\,
\frac{\phi(E_n^+(t),2\pi)}{\Delta_+^{\bullet}(E_n^+(t))}\;
\left(\int_0^\pi dx\, \theta(E_n^+(t),x)^2 +o(1)\right)^2.
\end{equation}
According to \eqref{6.1}, there exists a finite
positive constant $C$  such that for every set $\sigma\subset\sigma(H)$
and every element $g\in L^2(\mathbb{R})$ the inequality
\begin{equation}
|(\ol{g},P(\sigma)g)_{L^2(\bbR)}|\leq C\|g\|_{L^2(\bbR)}^2
\end{equation}
holds. Thus,
\begin{equation}
\bigg|\int_{T} dt\,
\frac{\phi(E_n^+(t),\pi)}{{\Delta}_+^{\bullet}(E_n^+(t))}
[1+o(1)] \bigg|\leq C\|f_n\|_{L^2(\bbR)}^2\leq C |T|.
\end{equation}
Since $\sigma$ is an arbitrary closed part of the spectral arc
$\sigma_n$ of $H$, we find
\begin{equation}\label{l141}
\left|\frac{\phi(\lambda,\pi)}{\Delta_+^{\bullet}(\lambda)}\right|
\leq C,\quad \lambda\in\sigma_{{\rm ext},R} (H)
\end{equation}
proving that the first function in \eqref{l14} is bounded on the set
$\sigma(H)$.

Next, we replace $F_n(x,t)$ by the function
\begin{equation}
G_n(x,t)=\begin{cases}
\beta_n (t)\theta(E_n^+(t),x)+
\sqrt{E_n^+(t)}\phi(E_n^+(t),x), &t\in T\cup T^*, \\
0, &  \text{otherwise}
\end{cases}
\end{equation}
with
\begin{equation}
\beta_n (t)=-\frac{\sqrt{E_n^+(t)}
\big(\ol{\theta(E_n^+(t),\cdot)},
\phi(E_n^+(t),\cdot)\big)_{L^2([0,\pi])}}
{\big(\ol{\theta(E_n^+(t),\cdot)},
\theta(E_n^+(t),\cdot)\big)_{L^2([0,\pi])}},
\quad t\in[0,2\pi].
\end{equation}
Then
\begin{equation}
\int_0^\pi dx\, G_n(x,t)\theta(E_n^+(t),x)=0
\end{equation}
and
\begin{equation}
\int_0^{\pi}\int_0^{2\pi} dx dt\, |G_n(x,t)|^2  \leq C |T|.
\end{equation}
Using \eqref{l361} we obtain
\begin{align}
\begin{split}
& (\ol{g_n},P(\sigma)g_n)_{L^2(\bbR)} \\
& \quad = \frac{1}{2\pi}\int_{T} dt\,
\frac{\theta'(E_n^+(t),\pi)}{E_n^+(t)\Delta_+^{\bullet}(E_n^+(t))}
\left(E_n^+(t)\int_0^\pi dx\, \phi(E_n^+(t),x)^2+o(1) \right)^2,
\end{split}
\end{align}
where $g_n\in L^2(\mathbb{R})$ is the inverse Gel'fand transform of
$G_n(x,t)$. This implies the estimate
\begin{equation}\label{l142}
\left|\frac{\theta'(\lambda,\pi)}
{\lambda\Delta_+^{\bullet}(\lambda)}
\right| \leq C,\quad \lambda\in\sigma_{{\rm ext}, R} (H).
\end{equation}

Finally, set
\begin{equation}
R_n(x,t)=\begin{cases}
\theta(E_n^+(t),x)+\sqrt{E_n^+(t)}\phi(E_n^+(t),x),&t\in T\cup T^*,\\
0,&  \text{otherwise,}
\end{cases}
\end{equation}
and denote by $r_n(x)$ the inverse Gel'fand transform of $
R_n(x,t)$. Then 
\begin{align}
& (\ol{r_n},P(\sigma) r_n)_{L^2(\bbR)}
=\f{1}{2\pi}\int_T dt\, \frac{\phi(E_n^+(t),\pi)}
{\Delta_+^{\bullet}(E_n^+(t))}
\big(\ol{\theta(E_n^+(t),\cdot)},R_n(\cdot,t)\big)_{L^2([0,\pi])}^2
\no \\
& \quad -\f{1}{2\pi}\int_T dt\,
\frac{\theta'(E_n^+(t),\pi)}{\Delta_+^{\bullet} (E_n^+(t))}
\big(\ol{\phi(E_n^+(t),\cdot)},R_n(\cdot,t)\big)_{L^2([0,\pi])}^2 \\
& \quad
-\f{1}{\pi}\int_T dt\, \frac{\Delta_-(E_n^+(t))}{\Delta_+^{\bullet}
(E_n^+(t))}\big(\ol{\theta(E_n^+(t),\cdot)},
R_n(\cdot,t)\big)_{L^2([0,\pi])}
\big(\ol{\phi(E_n^+(t),\cdot)},R_n(\cdot,t)\big)_{L^2([0,\pi])}.  \no
\end{align}
Since
\begin{align}
\big(\ol{\theta(E_n^+(t),\cdot)},R_n(\cdot,t)\big)_{L^2([0,\pi])}
&=\frac{\pi}{2}+o(1), \\
\sqrt{E_n^+(t)}\big(\ol{\phi(E_n^+(t),\cdot)},
R_n(\cdot,t)\big)_{L^2([0,\pi])}
&=\frac{\pi}{2}+o(1),
\end{align}
and the estimates \eqref{l141} and \eqref{l142} are already proved, we have
\begin{equation}
\left|\int_T dt\, \frac{\Delta_-(E_n^+(t))}
{\sqrt{E_n^+(t)}\Delta_+^{\bullet}(E_n^+(t))}
[1+o(1)]\right|\leq C |T|,
\end{equation}
and
\begin{equation}
\bigg|
\frac{\Delta_-(\lambda)}{\sqrt{\lambda}\Delta_+^{\bullet}(\lambda)}
\bigg|\leq C<\infty,
\quad\lambda\in\sigma_{{\rm ext}, R} (H),
\end{equation}
completing the proof of Lemma \ref{l7.3}.
\end{proof}

The necessity of the conditions in Theorem \ref{t4.3} now follows from
the necessity of the conditions in Theorem \ref{t4.4} and the identity
\eqref{l21}.

\begin{proof}[Proof of necessity of the conditions in Theorem \ref{t4.5}]
To prove necessity of condition $(i)$ of Theorem \ref{t4.5} we note
that if $\lambda_0$ is a multiple point of the periodic or antiperiodic
spectra then $\lambda_0\in\sigma(H)$ and $\Delta_+^\bullet(\lambda_0)=0$
and the claim follows from the analyticity of the first fraction in
\eqref{l16}.

Now, let $\lambda_0$ be an eigenvalue of an operator $H(t), \;
t\in[0,\pi]$. Then
$\Delta_+(\lambda_0)=\cos(t)+c(\lambda-\lambda_0)^m[1+o(1)]$ with some
$c\neq0$ and an integer $m\geq1$. Let
\begin{equation}
\mathcal{E}(\lambda_0,t)=\{f\in\dom(H(t)^m) \,|\,
(H(t)-\lambda_0I)^mf=0\}   \lb{6.82}
\end{equation}
and
\begin{equation}
\ker(H(t)-\lambda_0 I)=\{f\in\dom(H(t))\,|\,(H(t)-\lambda_0I)f=0\}
\lb{6.83}
\end{equation}
the corresponding algebraic (root) and geometric eigenspaces,
respectively. (Cf.\ \cite{GW95} for a detailed discussion of algebraic
multiplicities in terms of the behavior of the discriminant $\Delta_+$.)
If $m=1$ then both subspaces coincide. In the case
$m\geq2$ we have $\Delta_+^{\bullet}(\lambda_0)=0$ and, according to
\eqref{l52}, either $t=0$ or $t=\pi$, and in both cases $m=2$,
$\mathcal{M}-\Delta_+(\lambda_0)I_2=0$ and
$2=\dim(\ker(H(t)-\lambda_0))\leq\dim(\mathcal{E}(\lambda_0,t))$=2.
Hence, $\ker(H(t)-\lambda_0)=\mathcal{E}(\lambda_0,t)$, proving claim
$(ii)$ of Theorem \ref{t4.5}.

To prove that $(iii)$ is also nessesary, we first note that
for every $\delta>0$ there exists a constant $C=C(\delta)>0$ such that
\begin{equation}
C^{-1}\leq|\sin\pi \zeta|e^{-\pi|\Im(\zeta)|}\leq C,
\quad \zeta\notin\bigcup_{k\in\Z}D_k,\quad D_k=
\{\zeta' \in\bbC \,|\,|\zeta'-k|\leq\delta\}.
\end{equation}
In the following we use the abbreviations,
\begin{equation}
\ell_k^{\pm}=\sqrt{\lambda_k^{\pm}}, \quad m_k=\sqrt{\mu_k}, \quad
d_k=\sqrt{\delta_k},   \lb{6.85}
\end{equation}
where, as before,  $\{\lambda_k^{\pm}\}_{k\in\N}$ are the periodic/antiperiodic spectra, 
$\{\mu_k^{\pm}\}_{k\in\N}$ denotes the Dirichlet spectrum, and 
$\{\delta_k^{\pm}\}_{k\in\N}$
is the set of critical points of $\Delta_+(\mu)$. Moreover, we use the notation 
(cf.\ Theorem \ref{t4.2})
\begin{equation}
u_+(\zeta)=\Delta_+(\zeta^2), \quad u_-(\zeta)=\Delta_-(\zeta^2),  \quad 
s(\zeta)=\phi(\zeta^2,\pi).   \lb{6.86}
\end{equation}
Next, we fix a small $\delta$ and use representations \eqref{l601}-\eqref{l62}
to arrive at the estimates
\begin{equation}\label{l63}
C^{-1}\leq \left|\frac{\zeta s(\zeta)}{\zeta-m_k}\right|\leq C,\quad
C^{-1}\leq \left|\frac{u_+^{\bullet}(\zeta)}{\zeta-d_k}\right|\leq C,
\quad \zeta\notin D_k,
\end{equation}
for all $k\geq k_0$ with sufficiently large $k_0$. In addition, we will
use the expansions
\begin{equation}\label{l52a}
u_+(\zeta)=\gamma_k+\frac{1}{2}u_+^{\bullet\bullet}(d_k)(\zeta-d_k)^2
[1+\varepsilon_k(\zeta)],\quad \zeta\in D_k,
\end{equation}
where
\begin{equation}
\lim_{k\uparrow\infty}\max_{\zeta\in
D_k}|\varepsilon_k(\zeta)|=0.
\end{equation}
Substituting $\zeta=\ell_k^\pm$ we obtain
\begin{equation}\label{l53}
\ell_k^\pm-d_k=\pm\frac{1}{\pi}(\gamma_k^2-1)^{1/2}[1+o(1)], \quad k\geq
k_0
\end{equation}
and
\begin{equation}\label{l54}
\ell_k^+-\ell_k^-=\frac{2}{\pi}(\gamma_k^2-1)^{1/2}[1+o(1)], \quad k\geq
k_0.
\end{equation}
Therefore, there exists a constant $C>0$ such that
\begin{equation}
C^{-1}|\ell_k^+-\ell_k^-|\leq|\ell_k^\pm-d_k|\leq
C|\ell_k^+-\ell_k^-|.
\end{equation}
Since
\begin{align}\label{l55}
\begin{split}
u_-(m_k)^2& =u_+(m_k)^2-1  \\
& = \gamma_ku_+^{\bullet\bullet}(d_k) \big\{(m_k-d_k)^2
[1+o(1)]+(\ell_k^\pm-d_k)^2[\pi^2+o(1)]\big\},
\end{split}
\end{align}
we have
\begin{equation}\label{l56}
|u_-(m_k)|^2\leq C(|m_k-d_k|^2+|\ell_k^\pm-d_k|^2).  
\end{equation}
Next, let $\lambda=\zeta^2, \zeta \in D_k\cap \sigma'(H)$ with
\begin{equation}\label{l561}
\sigma'(H)=\{\zeta\in\bbC \,|\, u_+(\zeta)=\cos (t), \, 0\leq t\leq \pi\}.
\end{equation}
 Then for every $k\in\mathcal{Q}$ 
defined in \eqref{l610}, the inequalities \eqref{l19} and \eqref{l63} imply
\begin{equation}\label{l57}
\left|\frac{\zeta-m_k}{\zeta-d_k}\right|
=\left|\frac{\zeta-m_k}{2\zeta\phi(\zeta^2)}\right|
\left|\frac{u_+^\bullet(\zeta)}{\zeta-d_k}\right|
\left|\frac{\phi(\lambda,\pi)}{\Delta_+^{\bullet}(\lambda)}\right|\leq C,
\quad \zeta\in D_k\cap\sigma'(H), 
\end{equation}
proving an estimate
\begin{equation}\label{l58}
\left|\frac{d_k-m_k}{\zeta-d_k}\right|\leq C, \quad
\zeta\in D_k\cap\sigma'(H).
\end{equation}
On the other hand, $|u_-(\zeta)-u_-(m_k)|\leq C|\zeta-m_k|$, and hence
\begin{equation}
\left|\frac{u_-(m_k)}{\zeta-d_k}\right|\leq
\left|\frac{\Delta_-(\lambda)}{2\lambda^{1/2}\Delta_+^\bullet(\lambda)}
\frac{u_+^\bullet(\zeta)}{\zeta-d_k}\right|+C\left|
\frac{\zeta-m_k}{\zeta-d_k}\right|.
\end{equation}
Combining the latter inequality with \eqref{l19} we obtain
\begin{equation}\label{l59}
\left|\frac{u_-(m_k)}{\zeta-d_k}\right|\leq C,\quad
\zeta\in D_k\cap\sigma'(H).
\end{equation}
Therefore, \eqref{l55} and \eqref{l58} yield
\begin{equation}\label{l59a}
\left|\frac{\ell_k^\pm-d_k}{\zeta-d_k}\right|+
\left|\frac{(\gamma_k^2-1)^{1/2}}{\zeta-d_k}\right|
\leq C,\quad
\zeta\in D_k\cap\sigma'(H).
\end{equation}
Inequalities \eqref{l60} follow from \eqref{l58} and \eqref{l54}.
\end{proof}

\section{Sufficient Conditions for a Hill Operator to be \\ a Spectral
Operator of Scalar Type}  \lb{s8}


Similar to the proof of Theorem \ref{t4.1}, we define the Hilbert space $L^2(\sigma)^2$, $\sigma\subseteq\sigma(H)$, of measurable $\bbC^2$-vector elements with the finite norm \eqref{l508}. The following expansion theorem
contains an exact description of the Fourier--Floquet transform generated by $H$
and will turn out to be a cornerstone in the proof of sufficiency of the conditions in Theorem \ref{t4.4}.

\begin{theorem}\label{t7.1}
Suppose the conditions \eqref{l18} of Theorem \ref{t4.4} are satisfied. Then
for every element $g\in L^2(\R)$ 
there exists the $L^2(\sigma(H))^2$-limit
\begin{equation}
{\mathbf F}(\lambda;g) =\slim_{R\uparrow\infty} \int_{-R}^R dy\, 
{\mathbf Y}(\lambda,x) g(y)  \lb{9.21}
\end{equation}
and the $L^2(\R)$-representation
\begin{align}
g(x)&=\frac{1}{2\pi} \slim_{R\uparrow \infty}
\int_{\sigma(H)\cap\{\lambda\in\bbC\,|\, |\lambda|\leq R\}}
\frac{d\lambda}{\sqrt{1-\Delta_+(\lambda)^2}} \no \\
& \quad
\times \Big\{\big[\phi(\lambda,\pi)\theta(\lambda,x)
-\Delta_-(\lambda)\phi(\lambda,x)\big] F_\theta(\lambda;g)\lb{9.210}  \\
& \hspace*{.95cm} -\big[\theta'(\lambda,\pi)\phi(\lambda,x)
+\Delta_-(\lambda)\theta(\lambda,x)\big]F_\phi(\lambda;g)\Big\}   \no
\end{align}
is valid.
\end{theorem}
\begin{proof}
Assume that \eqref{l18} is valid and $\Delta_+^{\bullet}(\lambda_0)=0$
for some $\lambda_0\in\sigma(H)$. Then \eqref{l48} follows, \eqref{l21}
implies \eqref{l49}, and we find that the function
\begin{equation}
\omega(\lambda)=\frac{\Delta_+(\lambda)^2-1}{\Delta_+^{\bullet}(\lambda)^2}
\end{equation}
is analytic in a neighborhood of $\lambda_0$. Moreover, for
$\lambda\to\lambda_0$ we have
\begin{equation}
\Delta_+(\lambda)=\pm1+c(\lambda-\lambda_0)^p[1+o(1)],\quad
\Delta_+^{\bullet}(\lambda)=cp(\lambda-\lambda_0)^{p-1}[1+o(1)]
\end{equation}
with some integer $p\geq2$ and $c\neq0$, and
$\omega(\lambda)=\pm2c^{-1}p^{-2}(\lambda-
\lambda_0)^{-p+2}[1+o(1)]$. The analyticity of $\omega(\lambda)$ at $\lambda_0$ is
possible only if $p=2$, which implies that the point $\lambda_0$ is a
simple zero of $\Delta_+^{\bullet}(\lambda)$. Therefore, no spectral arcs
of $H$ cross at interior points, at most two spectral arcs of $H$ can meet
at a point $\lambda_0\in\bbC$ with $\Delta_+(\lambda_0)^2=1$, and the
spectrum of $H$ is of the form
\begin{equation}
\sigma(H)=\bigcup_{n=1}^\infty\Lambda_n
\end{equation}
with the spectral arcs $\Lambda_n$ of $H$ being given by
\begin{align}
\begin{split}
& \Lambda_n=\{z \in\bbC\,|\, z=E_{n-1}(t),\, \Delta_+(E_{n-1}(t))=\cos(t), \, t\in[0,\pi], \\ 
& \hspace*{4.85cm}  \Delta_+^{\bullet}(E_{n-1}(t))\neq0, \, t\in(0,\pi)\}.
\end{split}
\end{align}
In addition, all functions in
\eqref{l16}  are analytic in a neighborhood of $\sigma(H)$ and the function
\begin{equation}
\frac{\Delta_+(\lambda)^2-1}{(|\lambda|+1)\Delta_+^{\bullet}(\lambda)^2}
\end{equation}
is bounded on $\sigma(H)$.

If $\Delta_+(\lambda)^2-1=\Delta_+^{\bullet}(\lambda)=0$ for some $\lambda
\in\sigma(H)$, then the spectral arcs meeting at the point $\lambda$ are not regular. Nevertheless, according to conditions \eqref{l18}, the constant $C=C(\sigma)$ in \eqref{l501} is finite and is independent of the spectral arc  $\sigma=\Lambda_n\subseteq\sigma(H)$.

If $g\in L^2(\R)$, then using Lemma \ref{l510} one obtains ${\mathbf F}(\cdot;g)\in L^2(\Lambda_n)^2$
for all $n\in\N_0$, and 
\begin{align}
& \|{\mathbf F}(\lambda;g)\|_{L^2(\sigma(H))^2}^2=
\sum_{n=1}^\infty\|{\mathbf F}(\lambda;g)\|_{L^2(\Lambda_n)^2}^2   
\label{l506}  \\ 
& \quad \leq
C^2\displaystyle\int_0^{2\pi} dt \sum_{n=0}^\infty 
\Big[\big|{\wti F}_\theta(E_n(t);G(.,t))\big|^2
+\displaystyle\big|\sqrt{|E_n(t)|+1} {\wti F}_\phi(E_n(t);G(.,t))\big|^2\Big].  \no 
\end{align}
Corollary \ref{c5.3} yields the inequality
\begin{equation}\label{l509}
\|{\mathbf F}(\lambda;g)\|_{L^2(\sigma(H))^2}^2\leq
C^2\int_{0}^{2\pi}dt\|G(\cdot,t)\|^2_{L^2([0,\pi])}=C^2\|g\|_{L^2(\R)}^2,
\end{equation}
proving the existence of the limit in \eqref{9.21}. Therefore, the mapping ${\mathbf T}g={\mathbf F}(\lambda;g)$ defines a bounded linear operator from ${L^2(\R)}$ into 
$L^2(\sigma(H))^2$. 

To prove that ${\mathbf T}$ is a surjective map, we now assume that 
${\mathbf F}\in L^2(\sigma(H))^2$ and set 
${\mathbf F}_k(\lambda)=\chi_{k}(\lambda){\mathbf F}(\lambda)$, $k\in\N$, where
$\chi_{k}$ is the characteristic function of $\Lambda_k$. We have seen in the proof of Theorem \ref{t4.1} that there exists a sequence of functions $\{v_k(x)\}_{k=1}^\infty \subset L^2(\R)$ such that 
\begin{equation}
P(\Lambda_k)v_k=v_k,\quad{\mathbf F}_k(\lambda;v_k)={\mathbf F}_k(\lambda),\;\lambda\in\Lambda_k.
\end{equation}
Next, we set 
\begin{equation}
w_n=\sum_{k=1}^n\;v_k=\sum_{k=1}^n\;P(\Lambda_k)v_k
\end{equation}
and let $h\in L^2{(\R)}$ be an arbitrary element. Then, as we have already proved, 
${\mathbf F}(\cdot;h)\in L^2(\sigma(H))^2$ and 
\begin{align}
\begin{split}
|(\overline{h},w_n)_{L^2(\bbR)}| & \leq\sum_{k=1}^n\;\frac{1}{4\pi}
B_{\Lambda_k}|({\mathbf F}_k,{\mathbf F}(\cdot;h))|  \\
& \leq C
\bigg(\sum_{k=1}^n\;\|{\mathbf F}_k\|_{L^2(\Lambda_k)^2}^2\bigg)^{1/2}\|h\|_{L^2(\R)}.
\end{split}
\end{align}
Therefore,  
\begin{equation}
\|w_n\|_{L^2(\R)}\leq C\bigg(\sum_{k=1}^n\;\|{\mathbf F}_k\|_{L^2(\Lambda_k)^2}^2\bigg)^{1/2} = C\|{\mathbf F}\|_{L^2(\sigma(H))^2}, 
\end{equation}
and there exists the $L^2(\R)$-limit $w=\slim_{n\uparrow\infty}w_n$. For  
$m\in\bbN$ one has in the norm of $L^2(\sigma(H))^2$,
\begin{equation}
\chi_{m} {\mathbf F}(\cdot;w)=\slim_{n\uparrow\infty}\sum_{k=1}^n\;
\chi_{m} {\mathbf F}(\cdot;v_k)=\slim_{n\uparrow\infty}\sum_{k=1}^n 
\chi_{m} {\mathbf F}_k={\mathbf F}_m. 
\end{equation}
Thus, ${\mathbf F}(\cdot;w)={\mathbf F}$ and $\ran({\mathbf T})=L^2(\sigma(H))^2$. 
Moreover, for every ${\mathbf F} \in L^2(\sigma(H))^2$ 
there exists the $L^2(\R)$-limit
\begin{equation}
w(x)=\slim_{n\uparrow\infty}\frac{1}{4\pi}\sum_{k=1}^n\;
B_{\Lambda_k}({\mathbf F},{\mathbf Y}(\cdot,x))
\end{equation}
such that ${\mathbf F}(\cdot;w)={\mathbf F}$ and the estimate 
\begin{equation}\label{l507}
\|w\|_{L^2(\R)}\leq C\|{\mathbf F}\|_{L^2(\sigma(H))^2}
\end{equation}
holds. In particular, these arguments are applicable to 
${\mathbf F}={\mathbf F}(\cdot;g)$ with an arbitrary element  
$g\in L^2(\R)$ and one concludes that there exists the $L^2(\R)$-limit
\begin{equation}\label{l512}
(P(\sigma(H))g)(x)=\slim_{n\uparrow\infty}\frac{1}{4\pi}\sum_{k=1}^n\;
B_{\Lambda_k}({\mathbf F}(\cdot;g),{\mathbf Y}(\cdot,x)).
\end{equation}
Combining \eqref{l509} and \eqref{l507} one finds that 
$P(\sigma(H))$ is a bounded linear operator on $L^2(\R)$.

For an element $g\in L^2(\R)$ we set $h=g-P(\sigma(H))g$. If 
$\sigma\subset\sigma(H)$ is an arbitrary regular spectral arc of $H$, then
\begin{equation}
P(\sigma)P(\sigma(H))g=\slim_{n\uparrow\infty}\sum_{k=1}^n\;
P(\sigma)P(\Lambda_k)g=P(\sigma)g,
\end{equation}
implying $P(\sigma)h=0$. According to the last statement of Theorem \ref{t4.1} this  implies $h=0$, $g=P(\sigma(H))g$ and $P(\sigma(H))=I$, which proves assertion 
\eqref{9.210}.
\end{proof}

\begin{proof}[Proof of sufficiency of the conditions in Theorem \ref{t4.4}]
Theorem \ref{t7.1} contains a description of the functional model for the operator $H$: The latter is 
equivalent to the operator of multiplication by the independent variable in the space
$L^2(\sigma(H))^2$. At this point we can extend \eqref{l361} to an arbitrary Borel set 
$\sigma\subseteq\sigma(H)$. Indeed, for every such set and an arbitrary element  
$g\in L^2(\R)$, we define 
\begin{equation}\label{l513}
(P(\sigma)g)(x)=\slim_{n\uparrow\infty}\frac{1}{4\pi}\sum_{k=1}^n\;
B_{\Lambda_k}(\chi_{\sigma}(\cdot){\mathbf F}(\cdot;g),{\mathbf Y}(\cdot,x))
\end{equation}
and find that the family 
$\{E_H(\omega)=P(\omega\cap\sigma(H))\}_{\omega\in\cB_{\bbC}}$ is a spectral 
resolution for the operator $H$. (Here $\cB_{\bbC}$ denotes the collection of Borel subsets of $\bbC$.)
\end{proof}

\begin{proof}[Proof of sufficiency of the conditions in Theorem
\ref{t4.3}]If the function in \eqref{l17a}
is analytic in an open neighborhood of $\sigma(H)$ and if
$\Delta_+^{\bullet}(\lambda_0)=0$, then
$\Delta_+(\lambda_0)^2-1 =\phi(\lambda_0,\pi)={\Delta}_-(\lambda_0)=0$,
the identity \eqref{l64} implies that $\lambda_0$ is a simple zero of
$\Delta_+^{\bullet}(\lambda_0)=0$, and the functions \eqref{l16} are
locally bounded on $\sigma(H)$.

Next, we will prove that \eqref{l19} implies the second inequality in \eqref{l18}.

In addition to employing the notation \eqref{6.86}, we also introduce 
\begin{equation}
c_1(\zeta)=\theta'(\zeta^2,\pi),
\end{equation}
and recall the abbreviations $m_k=\sqrt{\mu_k}$, $d_k=\sqrt{\delta_k}$.

It is sufficient to prove
\begin{equation}\label{l20}
\sup_{\zeta\in\sigma'(H)\cap S(\delta,k_0)}\left|\frac{c_1(\zeta)}
{\zeta u_+^{\bullet}(\zeta)}\right|<\infty,
\end{equation} 
where $\sigma'(H)$ is defined by \eqref{l561}, 
and $S(\delta,k_0)$ is the union of discs $D_k=\{\zeta\in\bbC \,|\,
|{\zeta}-{{d_k}}|\leq\delta, k\geq k_0\}$ for some small $\delta>0$ and
sufficiently large $k_0$.

Now we use the following expansions in $D_k$:
\begin{align}
\zeta s(\zeta)& =m_ks^{\bullet}(m_k)(\zeta-m_k)[1+o(1)],  \\
u_+(\zeta)&=\gamma_k+\frac{1}{2}
u_+^{\bullet\bullet}(d_k)(\zeta-d_k)^2[1+o(1)], \\
u_+^{\bullet}(\zeta)&=u_+^{\bullet\bullet}(d_k)(\zeta-d_k)[1+o(1)], \\
c_1(\zeta)&=-\zeta\sin(\pi\zeta)[1+o(1)].
\end{align}
Using  relation \eqref{l21}  we find
\begin{equation}
\f{c_1(m_k)}{m_k}=2\frac{u_+(m_k)u_+^{\bullet}(m_k)
-u_-(m_k)u_-^{\bullet}(m_k)}
{m_ks^{\bullet}(m_k)}
\end{equation}
and thus,
\begin{equation}
\left|\f{c_1(\zeta)}{\zeta}\right|\leq C(|u_+^{\bullet}(m_k)|+|u_-(m_k)|
+|\zeta-m_k|),\quad
\zeta\in D_k\cap\sigma'(H).
\end{equation}
As we have seen in the proof of Theorem \ref{t4.5}, conditions \eqref{l19}
imply \eqref{l57}--\eqref{l59}, which in turn lead to \eqref{l20}.
\end{proof}

\begin{proof}[Proof of sufficiency of the conditions in Theorem
\ref{t4.5}] 
To this end we
assume that $\Delta_+^{\bullet}(\lambda_0)=0$ for some
$\sigma(H(t)),\;t\in [0,\pi]$. Then
\begin{equation}
\Delta_+(\lambda)=\cos(t)+c(\lambda-\lambda_0)^p[1+o(1)],\quad
\Delta_+^{\bullet}(\lambda)=c\,p(\lambda-\lambda_0)^{p-1}[1+o(1)]
\end{equation}
with $c\neq0$ and some integer $p\geq2$. Since according to condition $(ii)$,
$p=\dim(\mathcal{E}(\lambda_0,t))=\dim(\ker(H(t)-\lambda_0 I))\leq
2$ (cf.\ \eqref{6.82} and \eqref{6.83}), we conclude that $p=2$. Thus,
the boundary value problem
\begin{equation}
(H(t)-\lambda_0I)y=0,\quad y(\pi)=e^{it}y(0),\;y'(\pi)=e^{it}y'(0),
\end{equation}
has two linearly independent solutions $y_1$ and $y_2$. In this case
the $\bbC^2$-vectors
\begin{equation}
{\mathbf Y}_1 = \begin{pmatrix} y_1(0) \\ y_1'(0) \end{pmatrix}, \quad
{\mathbf Y}_2 = \begin{pmatrix} y_2(0) \\ y_2'(0) \end{pmatrix}
\end{equation}
are linearly independent solutions of the equation
\begin{equation}
\left(
\mathcal{M}(t)-e^{it}I_2 \right){\mathbf Y}=0
\end{equation}
and $\rho=e^{it}$ is a double zero of the characteristic
equation $\rho^2-2\Delta_+(\lambda_0)\rho+1=0$. This is possible only if
$\Delta_+(\lambda_0)^2-1=0$. The upshot is that either $t=0$, or $t=\pi$
and hence $\lambda_0$ is a double zero of $\Delta_+(\lambda)^2-1$ and a
simple zero of $\Delta_+^\bullet(\lambda)$. We use \eqref{l64} once more
and obtain
$\phi(\lambda_0,\pi)={\Delta}_-(\lambda_0)=
\theta'(\lambda_0,\pi)=0$. Therefore, all functions in \eqref{l16} are
locally bounded on $\sigma(H)$. It remains to prove \eqref{l19} for
$\lambda\in\sigma_R(H)$ with a sufficiently large $R>0$. Moreover, because
of representations \eqref{l601}--\eqref{l62}, it is sufficient to prove
\eqref{l19} for $\lambda=\zeta^2,\;\zeta\in D_k\cap\sigma'(H)$, for all
sufficiently large $k\geq k_0$.

If $\delta_k\in\sigma(H)$ for such $k$, then the functions
\begin{equation}
\frac{\zeta\phi(\zeta^2)}{u_+^{\bullet}(\zeta)},\quad
\frac{u_-(\zeta)}{u_+^{\bullet}(\zeta)}
\end{equation}
are analytic on the open interior of the disc $D_k$ and, according to
\eqref{l63}, are bounded on the boundary $\partial D_k$ by a constant
independent of $k$. It follows from the maximum principle that
they are bounded by the same constant inside $D_k$.

If $\delta_k\notin\sigma(H)$ for $k\geq k_0$, then similar to \eqref{l57}
(cf.\ also the abbreviations in \eqref{6.85}),
\begin{align}
\left|\frac{\phi(\lambda,\pi)}{\Delta_+^{\bullet}(\lambda)}\right|
&=\left|\frac{2\zeta\phi(\zeta^2)}{\zeta-m_k}\right|
\left|\frac{\zeta-d_k}{u_+^\bullet(\zeta)}\right|
\left|\frac{\zeta-m_k}{\zeta-d_k}\right| \no \\
& \leq C
\left(1+\frac{|m_k-\ell_k^\pm|+|\ell_k^+-\ell_k^-|}{|\zeta-d_k|}\right),
\quad \zeta\in D_k\cap\sigma'(H).
\end{align}
At last, we apply \eqref{l56} and obtain
\begin{align}
     \bigg|\frac{\Delta_-(\lambda)}{\sqrt{\lambda}
\Delta_+^{\bullet}(\lambda)}\bigg|
&=\bigg|\frac{u_-(\zeta)}{u_+^\bullet(\zeta)}\bigg|\leq
C\frac{|u_-(m_k)|+|\zeta-m_k|}{|u_+^{\bullet}(\zeta)|} \no \\
& \leq C
\left(1+\frac{|m_k-\ell_k^\pm|+|\ell_k^+-\ell_k^-|}{|\zeta-d_k|}\right),
\quad \zeta\in D_k\cap\sigma'(H).
\end{align}

According to the assumptions \eqref{l60}, the inequalities \eqref{l19}
are satisfied, and by virtue of Theorem \ref{t4.3}, $H$ is a spectral
operator of scalar type.
\end{proof}

\section{Concluding Remarks} \lb{s9}

In our final section we further illustrate the principal results of this
paper in a series of remarks.

\begin{remark} \lb{r4.5} We emphasize again some of the main points and
consequences of Theorems \ref{t4.3}--\ref{t4.5}. \\
If the conditions \eqref{l18} are satisfied, then the following
assertions are valid: \\
\hspace*{4mm} $(\alpha)$ $H$ has no spectral singularities and
$\{P(\sigma)\}_{\sigma\subseteq \sigma(H)}$ defines a spectral
resolu- \\
\hspace*{1cm} tion for $H$ in the sense that
$E_H(\omega)=P(\omega\cap\sigma(H))$ for all Borel sets $\omega\in\bbC$. \\
\hspace*{4mm}  $(\beta)$ The spectrum of $H$ consists of a system of
countably many, simple, noninter- \\
\hspace*{1cm} secting, analytic arcs (the latter may degenerate into
finitely many simple \\
\hspace*{1cm} analytic arcs and a simple analytic semi-infinite arc, all
of which are nonin- \\
\hspace*{1cm} tersecting). \\
\hspace*{4mm} $(\gamma)$ The functions
\begin{equation}\label{l17} 
\frac{\phi(z,\pi)}{\Delta_+^{\bullet}(z)},\quad
\f{\Delta_-(z)}{\Delta_+^{\bullet}(z)}, \quad 
\f{\theta'(z,\pi)}{\Delta_+^{\bullet}(z)}, \quad
\frac{\Delta_+(z)^2-1}{\Delta_+^{\bullet}(z)^2},\quad
\frac{\Delta_+(z)^2-1-\Delta_-(z)^2}
{\phi(z,\pi)\Delta_+^{\bullet}(z)} 
\end{equation}
\hspace*{1cm}  
are analytic in an open neighborhood of $\sigma(H)$. 
\end{remark}

\begin{remark} \lb{r9.8} 
Theorem \ref{t7.1} shows that conditions \eqref{l18} are sufficient for establishing a spectral decomposition of non-self-adjoint Hill operators with complex-valued potentials  in complete analogy to the well-known spectral decomposition of self-adjoint Hill operators (cf.\ \cite{Ti50}, \cite[Ch.\ XXI]{Ti58}). In this sense, Hill operators which are spectral operators of scalar type, demonstrate their great similarity to the case of self-adjoint Hill operators. We also note that equation \eqref{9.210} can be rewritten in the form 
\begin{equation}
g(x)= \slim_{R\uparrow\infty}
\int_{\sigma(H)\cap\{\lambda\in\bbC\,|\, |\lambda|\leq R\}} {d\lambda} \, 
({\mathcal S}(\lambda){\mathbf F}(\lambda;g),{\mathbf Y}(\lambda,x))_{\bbC^2}, 
\end{equation}
where
\begin{equation}
{\mathcal S}
(\lambda)=
\frac{1}
{2\pi\sqrt{1-\Delta_+(\lambda)^2}} 
\begin{pmatrix}
\phi(\lambda,\pi) & -\Delta_-(\lambda) \\
-\Delta_-(\lambda) & -\theta'(\lambda,\pi)
\end{pmatrix}, \quad \lambda\in\sigma(H),
\end{equation}
represents the analog of the $2\times 2$ spectral matrix of $H$ familiar from the self-adjoint context (cf.\ \cite[Sect.\ XIII.5]{DS88},  \cite[Ch.\ VI]{Na68}, \cite[Ch.\ III]{Ti58}).
\end{remark}

\begin{remark} \lb{r9.2} Next, we will show that the two inequalities in \eqref{l19} are independent. We start by constructing a Hill operator for which
the first inequality in \eqref{l19} holds but the second does not.
To this end, let $\gamma_{2k-1}=-1$,
$\gamma_{2k}=1-\pi^2(\rho_{|k|}^2+i\omega_{|k|}^2)/2,\;k=0,\pm1,\dots$, with sufficiently small and rapidly vanishing
numbers $\rho_k>0$ and $\omega_k>0$ such that 
$\omega_k\rho_k^{-1}\underset{k\uparrow  \infty} \to 0$. According to \cite{Tk96} there exists a function $u_+(\zeta)$ satisfying property
$(ii)$ in Theorem \ref{t4.2} such that $\{\gamma_k\}_{k\in\Z}$ is its set of critical values. The corresponding critical points are of the form 
$d_k\underset{k\uparrow\infty}{=}k+o(1)$ and if $r>0$ is a sufficiently small fixed 
constant and $D_k=\{\zeta\in\bbC \,|\, |\zeta-d_k|\leq r\}$, then
\begin{equation}
u_+(\zeta)=\gamma_k-\pi^2(-1)^k(\zeta-d_k)^2
(1+o(1))/2,\quad \zeta\in D_{k}.
\end{equation}
An elementary analysis show that the set 
\begin{equation}
\sigma_{2k}'=\{\zeta\in\bbC \,|\, u_+(\zeta)=\cos (t),\;t\in[0,\pi]\}\cap D_{2k}
\end{equation}
is formed by two analytic arcs which are close to the sets
\begin{equation}
c_{2k}^\pm=
\Big\{\zeta\in\bbC \,|\, \zeta=d_{2k}\pm i\sqrt{\rho_{|k|}^2+i\omega_{|k|}^2-t^2/(\pi)^2}
\Big\}
\end{equation}
and therefore, there exists a point $\nu_{2k}\in\sigma_{2k}'$ such that $|\nu_{2k}
-d_{2k}|\leq2\omega_{k}$. We note that the condition $\omega_k>0$ eliminates the
intersection of spectral gaps. 

Furthermore, we set $s(\zeta)=-{u^\bullet(z)}/({\pi\zeta})$ and
define a function $u_-(\zeta )$ using the interpolation data
$u_-(d_k)^2=u_+(d_k)^2-1$. According to Theorem \ref{t4.2}, with $\rho_k>0$ and $\omega_k>0$ being sufficiently small, there exists a Hill operator $H$
such that 
$\Delta_+(\lambda)=u_+(\sqrt{\lambda}),\; \phi(\lambda,\pi)=
s(\sqrt{\lambda}),\;\Delta_-(\lambda)=u_-(\sqrt{\lambda})$.
It is evident that the first condition in \eqref{l19} is satisfied for every such operator. 
On the other hand, for all sufficiently large $|k|\in\N$, we have 
$|u_-(d_{2k})|=|\gamma_{2k}^2-1|^{1/2}\geq2^{-1}\rho_{|k|}$, and if  
$\lambda_k=\nu^2_{2k}$, then the estimates from below 
\begin{align} 
\left|\frac{\Delta_-(\lambda_k)}{2\lambda_k^{1/2}\Delta_+^{\bullet}(\lambda_k)}\right|&=\left|\frac{u_-(\nu_{2k})}{u_+^{\bullet}(\nu_{2k})}\right|
\geq10^{-1}\left(
\left|\frac{u_-(d_{2k})}{d_{2k}-\nu_{2k}}\right|-
\left|\frac{u_-(\nu_{2k})-u_-(d_{2k})}{d_{2k}-\nu_{2k}}\right|
\right)  \no \\
& \geq C \frac{\rho_{|k|}}{\omega_{|k|}}
\end{align}
hold with $C>0$ independent of $k$. Thus, $H$ is not a spectral operator
of scalar type and the second condition in \eqref{l19} cannot be omitted in  
Theorem \ref{t4.3}.

Next, we consider the Hill operator $H$ associated with the triple $\{\Delta_+(\lambda),\; \phi(\lambda)=s(\sqrt{\lambda}), 0\}$ with the same function $\Delta_+(\lambda)$ as above
and with the entire function $s(\zeta)$ uniquely determined by the interpolation data 
$s(l_k^+)=0$, $k\neq0$, and by its exponential type $\pi$. With the same 
$\lambda_k=\nu_k^2$ as above, one obtains
\begin{equation} 
\left|\frac{\phi(\lambda_k)}{\Delta_+^{\bullet}(\lambda_k)}\right|=\left|\frac{2\nu_{2k}s(\nu_{2k})}{u_+^{\bullet}(\nu_{2k})}\right|
\geq10^{-1}
\left|\frac{l^+_{2k}-\nu_{2k}}{d_{2k}-\nu_{2k}}\right|
\geq C \frac{\rho_{|k|}}{\omega_{|k|}}.
\end{equation}
Thus, the first inequality in \eqref{l19} fails, while the second obviously holds since 
$\Delta_-(\lambda)\equiv0$.
\end{remark}

\begin{remark} \lb{r9.10} It is evident that the first inequality in \eqref{l60}
follows from two remaining ones and thus it could have been omitted. Nevertheless, we decided to keep it in the statement of Theorem \ref{t4.5} because of the following well-known interpretation of $|\lambda_k^+-\lambda_k^-|$ in the self-adjoint situation: 
If $V$ is real-valued, then $\lambda_k^+$ and $\lambda_k^-$
are the end-points of spectral gaps, $|\lambda_k^+-\lambda_k^-|$ are the lengths of 
spectral gaps, the points $\mu_k$ of the Dirichlet spectrum are trapped inside the gaps,  and the second and third inequality in \eqref{l60} become redundant. However, we will show next that, generally speaking, they are indispensable for complex-valued potentials $V$.

Indeed, if all critical values $\gamma_k$ of
$\Delta_+(z)$ are real  and $|\gamma_k|\geq 1$, then the first condition in 
\eqref{l60} is satisfied. We can choose $\gamma_k$ exponentially close to
$(-1)^k$ and afterwards fix real zeroes $m_k$ of
$s(\zeta)=\phi(\zeta^2,\pi)$ outside spectral gaps in such a way that
$|m_k-d_k|=\varepsilon |k|^{-\text{any big integer}}$. The corresponding
non-self-adjoint Hill operator exists according to Theorem \ref{t4.2}, but
neither the second nor the third inequality in \eqref{l60} are satisfied.

The numerators of fractions in \eqref{l60} contain distances between the closest points of periodic, antiperiodic, and Dirichlet boundary problems generated by the expression \eqref{3.1a} in the space $L^2([0,\pi])$. Their denominator also has a spectral meaning: Multiplied by the factor 2, it is 
equivalent to the distance between adjacent spectral arcs and may be considered as the length of the gap between them. These lengths decay as $k$ grows to infinity and 
equation \eqref{l60}, if valid, states that the distances between these three spectra 
vanish at least at the same rate.

It is interesting to note that the smoothness or analyticity properties of $V$ are characterized in terms of $\ell^2$-norms of the numerators in \eqref{l60}, with no connection with the relative rate of their decay (cf.\ \cite{Tk92}, \cite{Tk01}, \cite{DM06}). 
\end{remark}

\begin{remark} \lb{r9.3} If $\lambda_0$  is a point of either the periodic
or antiperiodic spectrum of the operator $H$ and a simple Dirichlet eigenvalue of
\eqref{3.18} such that $\Delta_+^\bullet(\lambda_0)=0$, then the analyticity of the first two functions in \eqref{l17} at $\lambda_0$ implies that $\lambda_0$ is a simple zero of
$\Delta_+^{\bullet}(z)$ and the three remaining functions are analytic at
$\lambda_0$ as well. We will next show that the situation is different
if a Dirichlet eigenvalue $\lambda_0$
has algebraic multiplicity larger than one, that is,
$\phi(\cdot,\pi)$ has a zero at $\lambda_0$ of order larger than one
(cf.\ \cite{GW95}).

According to \cite{Tk02}, an arbitrary complex number $\lambda_0$ may be
simultaneously a point of the Dirichlet spectrum of multiplicity $m>2$ and
a point of the periodic spectrum of multiplicity $2$ of the operator $H$
in \eqref{3.2}. We fix $\lambda_0\in\C$, set $m=3$, and consider a
corresponding Hill operator $H$. Following \cite{ST96}, \cite{Tk96}, we
can assume  that the triple $\{\phi(z,\pi), \Delta_+(z), \Delta_-(z)\}$
parameterizing such an operator has, in addition to properties $(i)$--$(v)$
in Theorem \ref{t4.2}, the following properties $(vi)$--$(ix)$:

$(vi)$ If $z \to\lambda_0$, then
\begin{equation}
\Delta_+(z)=1+a(z-\lambda_0)^2[1+o(1)],\quad \phi(z,\pi)=
b(z-\lambda_0)^3[1+o(1)]
\end{equation}
with $ab\neq0$.

$(vii)$ There exists $R>0$ such that if $\lambda\in\sigma_{{\rm
ext}, R} (H),\; \Delta_+^{\bullet}(\lambda)=0$, then
$\Delta_+^{\bullet\bullet}(\lambda)\neq0$ and
\begin{equation}
\Delta_+^{\bullet}(\lambda)^2-1={\phi(\lambda,\pi)}=\Delta_-(\lambda)=
\lim_{z\to\lambda}\frac{\Delta_+(z)^2-1-\Delta_-(z)^2}
{\phi(z,\pi)}=0.
\end{equation}

$(viii)$ All critical points of $\Delta_+(\cdot)$ in the disc
$\{z\in\bbC \,|\, |z|\leq R\}$, with the exception of $\lambda_0$, do
not belong to the set $\sigma(H)$.

$(ix)$ All zeros of $\phi(\cdot,\pi)$, except $\lambda_0$, are simple.

All conditions of Theorems \ref{t4.3} and \ref{t4.4} are satisfied
for such an operator $H$, except perhaps, the analyticity of the function
\begin{equation}\label{l41}
\frac{\Delta_+(z)^2-1-\Delta_-(z)^2}
{\phi(z,\pi)\Delta_+^{\bullet}(z)}
\end{equation}
at $\lambda_0$. For $z\to\lambda_0$ we have
\begin{align}
&\Delta_+(z)^2-1-\Delta_-(z)^2=
\big[\Delta_+^{\bullet\bullet}(\lambda_0)
-\Delta_-^{\bullet}(\lambda_0)^2\big](z-\lambda_0)^2 \\
& \quad + \big[3^{-1}{\Delta}^{\bullet\bullet\bullet}_+(\lambda_0)-
\Delta_-^{\bullet}(\lambda_0)
\Delta_-^{\bullet \bullet}(\lambda_0)\big](z-\lambda_0)^3+o((z-\lambda_0)^3)
=O((z-\lambda_0)^3),  \no
\end{align}
and hence  
$\Delta_-^{\bullet}(\lambda_0)=
\pm\sqrt{\Delta_+^{\bullet\bullet}(\lambda_0)}\neq0$. In addition to the latter
relation,  the analyticity of the function \eqref{l41} leads to
another relation, $\Delta^{\bullet\bullet\bullet}_+(\lambda_0)
-3\Delta_-^{\bullet}(\lambda_0)
\Delta_-^{\bullet\bullet}(\lambda_0)=0$. In any case, whether the latter
relation is satisfied or not,
we construct a function $v(\lambda)$ using the interpolation data
\begin{align}
\begin{split}
&v(\mu_k) =\Delta_-(\mu_k), \; \mu_k\neq\lambda_0, \;
v(\lambda_0)=0, \; v^{\bullet}(\lambda_0)
=\Delta_-^{\bullet}(\lambda_0), \\
& v^{\bullet\bullet}(\lambda_0) \neq
\Delta_+^{\bullet\bullet\bullet}(\lambda_0)/[3 v^{\bullet}(\lambda_0)]
\end{split}
\end{align}
such that
$\|\sqrt{z}[v(z)-\Delta_-(z)]\|_{{\mathbb{PW}_\pi}}$
is sufficiently  small and find the Hill operator $H$ (cf.\ \eqref{3.2})
corresponding to the triple $\{\phi(z,\pi),\Delta_+(z), v(z)\}$.
Conditions $(ii)$ of Theorem \ref{t4.3} are satisfied for such
an operator, but analyticity of the function \eqref{l41} at
$\lambda_0$ fails. In other words, the analyticity of the function
\eqref{l17a} does not follow from \eqref{l19} and condition $(i)$ in  
Theorem \ref{t4.3} is indispensable.
\end{remark}

\begin{remark} \lb{r9.4} Gasymov showed in \cite{Ga80} (cf.\ also \cite{Ga80a}, \cite{Sh03}) 
 that if
\begin{equation}
V(x)=\sum_{n=1}^\infty\ c_ne^{2inx}, \quad
\{c_n\}_{n\in\bbN}\in \ell^1(\bbN), \;\, x\in\bbR,  \lb{8.9}
\end{equation}
then $\Delta_+(z)=\cos(\pi\sqrt{z})$. Thus, in this case
$\sigma(H) =[0,\infty)$ and the function
\begin{equation}
\frac{\phi(z,\pi)}{\Delta_+^{\bullet}(z)}
=-\frac{2\sqrt{z}\phi(z,\pi)}
{\pi\sin(\pi\sqrt{z})}
\end{equation}
is analytic in an open neighborhood of $\sigma(H)$ if and only if
$\phi(z,\pi)=\sin(\pi\sqrt{z})/\sqrt{z}$. In the latter case
$\Delta_-(z)\equiv 0$ and $V(x)=0$ for a.e.\ $x\in\bbR$. This means that
no smoothness or analyticity conditions imposed on $V$ can guarantee that
the Hill operator \eqref{3.2} is a spectral operator of scalar type. Nevertheless, every 
Hill operator with a complex-valued locally
square-integrable potential is an operator with  a separable spectrum as
defined by Lyubich and Matsaev \cite{LM60}, \cite{LM62}.

For every operator $H$ with a nontrivial potential \eqref{8.9} there
exists at least one integer $n_0\in\N$ such that $\phi(n_0^2,\pi)\neq0$ and the point
$n_0^2$ is then a spectral singularity. For a generic element $f\in L^2(\R)$,  
norm-convergent expansions of the form \eqref{9.210} are impossible, but 
Gasymov derived a regularized spectral expansion containing, in particular,
terms of the type 
\begin{equation}
\bigg(\int_{\R} dy \, \psi_+(n_0^2,y)f(y)\bigg) \psi_+(n_0^2,\cdot)
\end{equation}
not belonging to $L^2(\R)$, and converging for sufficiently smooth and sufficiently fast decaying functions $f$.

Regularized expansions of a similar type for arbitrary potentials $V$
were also proposed by Veliev in \cite{Ve80} (cf.\ also \cite{Ve83}--\cite{VT02}). However, from our point of view the proofs given in papers \cite{Ve80}--\cite{VT02} are not satisfactory since they completely ignore the possible growth 
of functions in \eqref{l18}, leaving disputable the convergence of the series
even for compactly supported continuous functions $f$.
\end{remark}

\begin{remark} \lb{r9.5} Explicit examples of crossings of spectral arcs
$\lambda(t)$ of $H$ at interior points (i.e., for values of $t$ different
from $0$ and $\pi$) were first constructed by Pastur and Tkachenko
\cite{PT91a}. An explicit example of two crossing arcs in terms of the
Weierstass $\wp$-function was found by Gesztesy and Weikard \cite{GW95}.
These crossings represent spectral singularities of the underlying Hill
operator $H$. As shown by Tkachenko \cite{Tk94}, these crossings, and hence
the existence of such spectral singularities, is unstable with respect to
arbitrarily small perturbations of $V$ in $L^2([0,\pi])$. More precisely, for any
$\varepsilon>0$ there exists a periodic potential $V_{\varepsilon}\in
L^2_{\loc}(\bbR)$ of period $\pi$ with
$\|V-V_{\varepsilon}\|_{L^2([0,\pi])}<\varepsilon$, such that the spectrum
of $H_{\varepsilon}=-d^2/dx^2+V_{\varepsilon}$ is a system of
nonintersecting regular analytic arcs. Moreover, each spectral arc of $H$
is mapped in a one-to-one manner on the interval $[-1,1]$ by the Floquet
discriminant $\Delta_+$. Thus, generically (in the
$L^2([0,\pi])$-sense just described), the spectrum of a Hill operator $H$
is formed by a system of nonintersecting analytic arcs. In addition,
an approximating potential may be chosen in such a way that all conditions of 
Theorems \ref{t4.3}--\ref{t4.5} are satisfied. Put differently,
generically, Hill operators $H$ are spectral operators of scalar type.
\end{remark}

\begin{remark} \lb{r9.6} There are two obstacles for the operator $H$ not to be a  spectral operator of scalar type, and both of them may be elucidated using
the direct integral representation 
\begin{equation}
P(\sigma)=\frac{1}{2\pi}
\int^{\oplus}_{\{t\in[0,2\pi] \,|\, \Delta_+(\lambda)=\cos (t),\;\lambda\in\sigma\}} 
dt\, P(t;\sigma), 
\end{equation}
where $P(t;\sigma)$ is the spectral projection generated by $H(t)$.

The first obstacle is the  presence of critical points 
of the Hill discriminant in the spectrum, and to explain their influence we assume that $t_0\in(0,\pi)$ and $\lambda_0\in\sigma(H(t_0))$ are such that $\Delta_+(\lambda_0)^2
-\cos (t_0)=\Delta_+^\bullet(\lambda_0)= \cdots =
\Delta_+^{(m)}(\lambda_0)=0,\;\Delta_+^{(m+1)}(\lambda_0)\neq0$. Then the root subspace (i.e., the algebraic eigenspace) of $H(t_0)$ corresponding to the eigenvalue 
$\lambda_0$ is of dimension $m+1\geq2$. On the other hand, since $\Delta_+(\lambda_0)^2-1=
((\theta(\lambda_0,\pi)-\phi'(\lambda_0,\pi)))/2)^2-
\theta'(\lambda_0,\pi)\phi(\lambda_0,\pi))\neq0$, at least one entry of the matrix $\cM(\lambda_0)-e^{it_0}I$ does not vanish and the corresponding geometric 
eigenspace is one-dimensional. 

If $t\neq t_0$ is close to $t_0$, then the spectrum of $H(t)$ in a neighborhood of 
$\lambda_0$ is formed by $m+1$ simple eigenvalues, and the corresponding eigenvectors span an $(m+1)$-dimensional invariant subspace. The angle between any pair of linearly independent 
eigenvectors tends to zero as $t$ tends to $t_0$ and the projection onto one of them along the other has its norm increasing to infinity. A similar situation arises if 
$t_0=0\, \text{\rm (mod}\, \pi)$ and $m\geq3$. In either case, the family $\{P(\sigma)\}$ with regular spectral arcs $\sigma$ close to $\lambda_0$ is not bounded and we have a blowup at this point.

Another obstacle for $H$ to be a spectral operator of scalar type is connected 
to the behavior of functions in \eqref{l16} at infinity. For every 
$t\neq 0\, \text{\rm (mod}\, \pi)$ the boundary conditions \eqref{3.28} are regular 
\cite{Na68} and the corresponding eigenvalues are asymptotically separated (i.e., 
asymptotically they are simple, and consecutive eigenvalues are a fixed minimal, possibly $t$-dependent  distance apart from each other). This has been discussed in \cite{Ke64}, \cite{Mi62}, and \cite[Ch.\ XIX]{DS88a} for more general boundary conditions than the ones in 
\eqref{3.28}. In addition, the normalized root system of $H(t)$ is a basis in the space 
$L^2([0,\pi])$ in this case. Moreover, as proved in \cite{VT02}, it is a uniform Riesz basis for $t\in[\varepsilon,\pi-\varepsilon]$, which means (in the absence of root 
functions) that the inequalities
\begin{equation}\label{l410}
C_\varepsilon^{-1}\|f\|_{L^2([0,\pi])}^2\leq\sum_{k\in\bbN_0}|({f},\hat{\psi}_\pm(E_k(t),\cdot))|^2\leq C_\varepsilon\|f\|_{L^2([0,\pi])}^2
\end{equation}
are valid for all such $t$'s with some finite constant $C_\varepsilon >0$. As a result, the family  of projections $\{P(\sigma)\}$ with $\sigma\subseteq \{\lambda\in\sigma(H) \,|\,  \Delta_+(\lambda)=\cos (t),\; \varepsilon\leq t\leq \pi-\varepsilon\}$ is bounded,
but the possible dependence of $C_\varepsilon$ on
$\varepsilon$ does not permit us to claim that $H$ is a spectral operator of scalar
type. 

If, however, $t=0\, \text{\rm (mod}\, \pi)$, then the boundary conditions in \eqref{3.28} remain  regular, but the eigenvalues of the operators $H(t)$ with such $t$'s
are grouped in converging pairs $(\lambda_k^+(t),\lambda_k^-(t))$ which may 
amalgamate as $k$ tends to infinity. In the special case where $H$ is self-adjoint, so is every operator $H(t)$, and because of the orthogonality of spectral projections, even in the worst case scenario of coincidence,  
$\lambda_k^+(t)=\lambda_k^-(t)$, the norm of the corresponding spectral projections 
$P(\sigma,t)$ remains equal to $1$. 

The situation is quite different for complex-valued potentials $V$. 
For $t$  close to $0\,\text{\rm (mod} \, \pi)$, the angle between
eigenspaces of $H(t)$ corresponding to $\lambda_k^-(t)$ and $\lambda_k^+(t)$ may
approach zero as $k\uparrow\infty$ with different $k$'s for different $t$'s
and we can find a sequence of regular spectral arcs $\sigma_m$ 
accumulating at infinity such that $\lim_{m\to\infty}\|P(\sigma_m)\|=\infty$. This  may be considered a blowup at the point infinity, and the operator $H$ in 
Remark \ref{r9.2} yields an example of such a behavior. The conditions described in Theorems \ref{t4.3}--\ref{t4.5} prevent both blowup phenomena from happening. 
\end{remark}

\begin{remark}\label{r9.9}
The problem we are discussing in this remark is similar to the following well-known
question: When is the exponential system $\{e^{i\mu_kx}\}_{k\in\Z}$ a Riesz basis in ${L^2([-\pi,\pi])}$? The investigation of such problems for non-orthogonal systems of exponentials $\{\exp(i\mu_k x)\}_{k\in\bbN}$, $x\in[-\pi,\pi]$, which is similar to those presently discussed, has a long history. 
A criterion for the latter system to be a Riesz basis in $L^2([-\pi,\pi])$
was found by Pavlov \cite{Pa79}. Most recently,  Minkin \cite{Mi06} made
essential progress in studying the same property for systems of eigenfunctions of 
two-point boundary problems for higher-order differential operators, and Makin 
\cite{Ma06} found sufficient conditions for the root system of a 
Schr\"odinger operator on the interval $[0,1]$ associated with periodic and antiperiodic boundary conditions to be (or not to be) a Riesz basis. Using well-known asymptotic formulas (see \cite{Ma86}) for the periodic/antiperiodic and Dirichlet spectra of the corresponding operator $H$ restricted to the interval $[0,1]$, it is easy to check
that under the assumptions made on $V$ in \cite{Ma86}, the conditions of Theorem
\ref{t4.4} are either satisfied (or not satisfied), with the corresponding conclusions about the system of eigenfunctions and its property of forming a Riesz basis. Recently, Djakov and Mityagin \cite{DM06} constructed a series of potentials in some weighted spaces of periodic functions for which the corresponding eigensystems are not Riesz bases in  $L^2([-\pi,\pi])$.

Returning to the eigensystems of the operators $H(t)$, we will next show that our
conditions \eqref{l18} are sufficient for
the eigensystem of every such operator to form a Riesz basis in the space $L^2([0,\pi])$ with a constant $C_\varepsilon$ in \eqref{l410} independent of $\varepsilon$.

First of all, we fix $t\in[0,\pi]$ with $t\neq 0\, \text{\rm (mod}\, \pi)$, and note that for every
$k\in\bbN_0$ such that $\phi(E_k(t),\pi)\neq0$, the numbers $w_\pm(E_k(t))$, with
$w_\pm(\lambda)$ being defined by \eqref{l363}, are finite, do not vanish, and
the functions
\begin{equation}
\hat{\psi}_\pm(E_k(t),x)=\frac{{\psi}_\pm(E_k(t),x)}{w_\pm(E_k(t))},  \quad
t\in[0,\pi], \; \text{$t\neq 0 \, {\rm (mod}\,\pi)$,}
\end{equation}
are normalized Floquet solutions of $L\psi=z\psi$. Second, for the function
\begin{equation}
m(\lambda)=-\frac{\phi(\lambda,\pi)w_+(\lambda)w_-(\lambda)}
{2\Delta_+^{\bullet}(\lambda)},
\end{equation}
equation \eqref{l25} yields  $|m(E_k(t))|\geq 1$ while, similar to \eqref{l501},
conditions \eqref{l18} imply $|m(E_k(t))|\leq C$, with $C>0$ independent of $t\in[0,\pi]$.

Furthermore, for each $k\in\N_0$ such that $\phi(E_k(t),\pi)=0$, we define the
normalized Floquet solutions of $L\psi=z\psi$ and the numbers $m(E_k(t))$
by the relations
$$\hat{\psi}_\pm(E_k(t),x))=\lim_{t\to s}
\hat{\psi}_\pm(E_k(s),x),\qquad m(E_k(t))=\lim_{t\to s}m(E_k(s)),$$
and obtain the estimates $1\leq|m(E_k(t))|\leq C$ with the same constant $C$ for 
all $t\in(0,\pi)$ and $k\in\N_0.$

At last, we write
the expansion \eqref{l37} for every $t\in(0,\pi)$ in the form
\begin{equation}
f(x)=\sum_{k\in\bbN_0}\;m(E_k(t))\;(\overline{f},\hat{\psi}_-(E_k(t),\cdot))
\hat{\psi}_+(E_k(t),x),
\end{equation}
and using Corollaries \ref{c5.2} and \ref{c5.3} prove the inequalities \eqref{l410}
for $t\in(0,\pi)$ with $C_\varepsilon$ actually independent of $\varepsilon$.

For $t=0\, \text{\rm (mod}\, \pi)$ the spectrum 
of $H(t)$ degenerates if $\Delta_+^{\bullet}(E_k(t))=0$ and the series in \eqref{l410}
must be modified. Assume, for instance, that $t=0$ and set ${\mathcal K}=
\{k\in\bbN_0 \,|\, |\Delta_+^{\bullet}(E_k(0))=0\}$. Then $E_k(0)=E_{k+1}(0)=\lambda_k$ 
for $k\in{\mathcal K}$ and there exists a small 
neighborhood $U_0$ of $t=0$, the same for all sufficiently large $k\in{\mathcal K}$, with 
solutions $E_k^+(t)$ and $E_{k}^-(t)$ of the equation $\Delta_+(\lambda)=\cos (t)$,
$t\in U_0$, given by \eqref{l33C}. Moreover, for all sufficiently large $k$ 
these solutions are analytic in $U_0$. Since the constant $C$ in \eqref{l33B} is 
independent of $t$, we can differentiate the representations \eqref{l33C} and 
obtain
\begin{equation}
\frac{dE_k^\pm(0)}{dt}=\pm\frac{4k}{\pi}+\frac{g_k^\pm}{k},\quad 
\sum_{k\in\mathcal K}|g_k^\pm|^2<\infty.
\end{equation}
According to Theorem \ref{t4.4}, the functions \eqref{l16}
are analytic in a neighborhood of $\lambda_k$, one has $\phi(\lambda_k,\pi)=
\theta(\lambda_k,\pi)-1=0$, and hence there exist the limits
\begin{align}\no
\begin{split}
 \displaystyle\psi^\pm_k(x)&=\lim_{t\to0}{\psi}_+(E^\pm_k(t),x) \\
& =\theta(\lambda_k,x)+\frac{i\left( \frac{dE_k^\pm(0)}{dt}\right)^{-1}+\theta^\bullet(\lambda_k,\pi)}
{\phi^\bullet(\lambda_k,\pi)}\;\phi(\lambda_k,x), \quad k\in\mathcal K.
\end{split}
\end{align}
These functions are two linearly independent Floquet solutions forming a basis in the eigenspace of $H(0)$ corresponding to its eigenvalue 
$\lambda_k$. Taking into account \eqref{l33D}, \eqref{l33E}, and 
\eqref{l33C}, we obtain the representation
\begin{equation}
\psi^\pm_k(x)=e^{\pm2\pi ix}+G_k(x),\quad \sum_{k\in\mathcal K}
\|G_k\|_{L^2([0,\pi])}^2 \leq C.
\end{equation}
After normalizing the system
\begin{equation}
\{\psi^\pm_k\}_{k\in\mathcal K}\cup 
\{\psi_+(E^\pm_k(0),\cdot)\}_{k\in\N_0\setminus\mathcal K}, 
\end{equation}
it forms a Riesz basis in $L^2([0,\pi])$.

The same arguments are valid for $t=\pi$ and hence inequalities 
\eqref{l410} are satisfied for all $t\in[0,2\pi]$. Thus, if the conditions of at least one of the Theorems \ref{t4.3}--\ref{t4.5} are satisfied, then the system of eigenfunctions of every operator $H(t)$, $t\in [0,2\pi]$, forms a Riesz basis in the space ${L^2([0,\pi])}$. 

Arguments similar to those used in the proof
of Theorem \ref{t4.4} show that the system of eigenfunctions of $H(t)$ with 
$t=0\, \text{\rm (mod}\, \pi)$ is a Riesz basis if and only if 
\begin{equation}
\bigg|\frac{\phi(\lambda_{k}(t),\pi)}
{\Delta_+^{\bullet}(\lambda_{k}(t))}
\bigg|+
\left|
\frac{\theta'(\lambda_{k}(t),\pi)}
{(|\lambda_{k}(t)|+1)\Delta_+^{\bullet}(\lambda_{k}(t))}
\right|+ 
\bigg|
\frac{\Delta_-(\lambda_{k}(t))}
{(\sqrt{|{\lambda_{k}(t)}|}+1)\Delta_+^{\bullet}(\lambda_{k}(t))}
\bigg|\leq C       \lb{l8.12}
\end{equation}
 for such $t$'s and all $k\geq1$. Of course, condition \eqref{l8.12} is weaker than \eqref{l18}, and, as the operator $H$ constructed in Remark \ref{r9.2} demonstrates, the former condition may be satisfied, while the latter may not.
As a result, the eigensystems of $H(t)$  are Riesz bases for all $t\in[0,2\pi]$, while $H$ is not a spectral operator of scalar type, in complete agreement with the unboundness of the 
family $\{C_\varepsilon\}_{\varepsilon>0}$ in \eqref{l410}. 
\end{remark}

\begin{remark} \lb{r9.7}
Finally, we briefly describe some historical developments related to spectral theory
of non-self-adjoint Hill operators: The notion of spectral projections
associated with regular spectral arcs (cf. Definition \ref{d3.4}) was
introduced by Tkachenko \cite{Tk64} in 1964. Spectral projections  were also discussed by Veliev
\cite{Ve80}--\cite{Ve86} who defined spectral singularities as described in 
\ref{d4.1A}. He was first to note  that $\lambda=\lambda_k(t)$ is a point 
of spectral singularity of $H$ if and only if the root subspace of $H(t)$ 
corresponding to $\lambda$ contains a root function which is not an actual eigenfunction of $H(t)$ . However, since he erroneously concluded that
Hill operators lead to non-intersecting analytic arcs (cf.\ \cite[Theorem
5]{Ve83}), the counter examples of Hill operators with crossing spectral
arcs in the interior of these arcs found in \cite{PT91a} (and
subsequently in \cite{GW95}) invalidate some of his results concerning
spectral projections and spectral singularities. A bit earlier spectral expansions for
non-self-adjoint Hill operators were also briefly touched upon by Meiman
\cite{Me77} (but no proofs of his claims were offered). Moreover, he noted
on page 846 that zeros of $\Delta_+^\bullet(\lambda)$  are integrable singularities
for the functions in \eqref{l16} with $\lambda=\lambda(t)$, but this is
incorrect for $t_0\in\{0,\pi,2\pi\}$. Indeed, for potentials of the
type \eqref{8.9} and $t_0\in\{0,\pi,2\pi\}$, one has
$\Delta_+^\bullet(\lambda(t))=(c+o(1))(t-t_0)$, $c\neq0$, and if
$\phi(\lambda(t_0),\pi)\neq 0$, the singularity at $t_0$ is not integrable.
Another claim of that paper on page 846, suitably paraphrased, 
is  the following: ``The crude asymptotic estimates
for $\phi(\lambda,\pi)$, $\Delta_+^\bullet(\lambda)$ and  the Floquet
solutions of $H$ remain valid in the complex case. From this\footnote{and some formulas in Section 6 of \cite{Me77}} it follows that the Fourier integral theory remains essentially valid  also for the expansion in eigenfunctions of a non-self-adjoint
Schr\"odinger operator with a complex-valued periodic potential.'' However,
such a vague statement, without explanations as to what type of convergence
and what type of function space is meant, is unsatisfactory, especially,
taking  into account Theorems \ref{t4.3}--\ref{t4.5}.

For spectral expansion theorems in connection with finite interval problems and limit
circle-type situations we refer to \cite{EFZ05}, \cite{Mi06}, \cite{VT02}, and
the literature cited therein. Spectral resolutions for a special class of
non-self-adjoint operators were also considered by Volk \cite{Vo63}.
Exponentially decaying perturbations of non-self-adjoint Hill operators
were studied by Zheludev \cite{Zh69}. A rather different route was chosen
by Marchenko \cite{Ma63}, \cite[Ch.\ 2]{Ma86} (see also \cite{Fu70}), who
introduced the notion of spectral distributions for non-self-adjoint
problems. For the notion of generalized spectral operators we refer, for
instance, to \cite{CF68}, \cite{Lj66} and the references therein.
\end{remark}

\appendix
\section{Spectral Operators in a Nutshell}
\lb{A}
\renewcommand{\theequation}{A.\arabic{equation}}
\renewcommand{\thetheorem}{A.\arabic{theorem}}
\setcounter{theorem}{0}
\setcounter{equation}{0}

In this appendix we recall a selection of basic facts on spectral
operators as discussed in great detail in volume 3 of Dunford and Schwartz
\cite{DS88a}. Unless explicitly stated otherwise, the material presented
is taken from \cite{DS88a} (see also \cite{Ba54}, \cite{Du58},
\cite{Fo58}, \cite{Sc60}) and the reader is referred to this monograph for
proofs and pertinent references on this subject. For simplicity we
will restrict our discussion to operators in a separable complex Hilbert
space. We note, however, that the material below is developed in a 
general Banach space context in \cite{DS88a}.

To set the stage, we assume the following conventions for the rest of this
section: $\cH$ denotes a separable, complex Hilbert space with scalar
product $(\cdot,\cdot)_{\cH}$ (linear in the second factor), norm
$\|\cdot\|_{\cH}$, and $I_{\cH}$ the identity operator in $\cH$. The
Banach space of bounded linear operators on $\cH$ will be denoted by
$\cB(\cH)$ with norm $\|\cdot\|_{\cB(\cH)}$, and the set of densely
defined, closed linear operators in
$\cH$ will be denoted by $\cC(\cH)$. The domain and range of a linear
operator $S$ are denoted by $\dom(S)$ and $\ran(S)$.

\begin{definition} \lb{d2.1} ${}$ \\
$(i)$ A {\it spectral measure} $E$ in $\cH$ is a homomorphic map of a
$\sigma$-algebra $\cA$ of sets into a Boolean algebra of projection
operators in $\cH$ such that the unit of $\cA$ is mapped to $I_{\cH}$. The
spectral measure $E$ is called {\it bounded} if for some $C>0$,
$\|E(\omega)\|_{\cB(\cH)}\leq C$ for all $\omega \in \cA$.  \\
$(ii)$ If $T\in\cC(\cH)$, then $\sigma \subseteq\sigma(T)$ is called a
{\it spectral set} if $\sigma$ is both open and closed in the topology of
$\sigma(T)$.
\end{definition}

In the concrete applications we have in mind in the bulk of this paper,
$\Omega=\bbC$ and $\Sigma$ typically equals the $\sigma$-algebra
$\cB_{\bbC}$ of Borel subsets of $\bbC$; thus we confine ourselves to
this case for the rest of this appendix.

Next we turn to the special case of bounded spectral operators on $\cH$.

\begin{definition} \lb{d2.2} Let $T\in\cB(\cH)$. \\
$(i)$ A projection-valued spectral measure $E$ on $\cB_{\bbC}$ is called
a {\it resolution of the identity} (or a {\it spectral resolution}) for
$T$ if
\begin{equation}
E(\omega)T=TE(\omega), \quad
\sigma\big(T\big|_{E(\omega)\cH}\big)\subseteq
\ol{\omega}, \quad \omega\in \cB_{\bbC}.  \lb{2.1}
\end{equation}
$(ii)$ A projection-valued spectral measure $E$ in $\cH$ defined on
$\cB_{\bbC}$ is called {\it countably additive} if for all
$f, g \in \cH$, $(f,E(\cdot)g)_{\cH}$ is countably additive on
$\cB_{\bbC}$. \\
$(iii)$ $T$ is called a {\it spectral operator} if it has
a countably additive resolution of the identity defined on $\cB_{\bbC}$.
\end{definition}

\begin{lemma} \lb{l2.3} ${}$ \\
$(i)$ Any countably additive projection-valued spectral measure $E$ on
$\cB_{\bbC}$ is countably additive in the strong operator topology and
bounded. \\
$(ii)$ Let $T\in\cB(\cH)$ be a spectral operator, then
$E(\sigma(T))=I_{\cH}$. \\
$(iii)$ Every bounded spectral operator has a uniquely defined countably
additive resolution of the identity defined on $\cB_{\bbC}$.
\end{lemma}

The uniquely defined spectral resolution for $T$ will frequently be
denoted by $E_T$ in the following.

The important special case of bounded scalar spectral operators is
introduced next.

\begin{definition} \lb{d2.4} ${}$ \\
$(i)$ Let $S\in\cB(\cH)$ be a spectral operator with spectral resolution
$E_S$ defined on $\cB_{\bbC}$. Then $S$ is said to be of {\it scalar
type} (or a {\it scalar spectral operator}) if
\begin{equation}
S=\int_{\bbC} \lambda \, dE_S (\lambda).  \lb{2.2}
\end{equation}
$(ii)$ $N\in\cB(\cH)$ is called {\it quasi-nilpotent} if
$\lim_{n\uparrow \infty} \|N^n\|_{\cB(\cH)}^{1/n}=0$.
\end{definition}

\begin{lemma} \lb{l2.5} ${}$ \\
$(i)$ If $E$ is a countably additive projection-valued spectral measure on
$\cB_{\bbC}$ which vanishes outside a compact subset of $\bbC$, then
\begin{equation}
S = \int_{\supp\,(dE)} \lambda \, dE(\lambda)  \lb{2.3}
\end{equation}
is a bounded spectral operator of scalar type whose spectral resolution
is $E$. \\
$(ii)$ $N\in\cB(\cH)$ is quasi-nilpotent if and only if
$\sigma(N)=\{0\}$.
\end{lemma}

The following is a principal result on bounded spectral operators.

\begin{theorem} [The canonical reduction of bounded spectral
operators] \lb{t2.6} ${}$ \\
Let $T\in\cB(\cH)$. Then $T$ is a spectral operator if and only if
$T=S+N$, where $S\in\cB(\cH)$ is a bounded spectral operator of scalar
type and $N$ is a quasi-nilpotent operator commuting with $S$. This
decomposition is unique and
\begin{equation}
\sigma(T)=\sigma(S).  \lb{2.4}
\end{equation}
Moreover, $T$ and $S$ have the same resolution of the identity.
\end{theorem}

Next we turn to unbounded spectral operators.

\begin{definition} \lb{d2.14}
Let $T\in\cC(\cH)$. Then $T$ is called a {\it spectral operator} if there
exists a regular, countably additive projection-valued spectral measure
$E$ (with respect  to the strong operator topology) defined on
$\cB_{\bbC}$ such that the following conditions hold:
\begin{align*}
& (\alpha) \; \dom(T)\supseteq E(\omega)\cH \, \text{ for
$\omega\in\cB_{\bbC}$ bounded}. \\
& (\beta) \; E(\omega)\dom(T) \subseteq \dom(T), \\
& \quad \;\;\,  TE(\omega)f=E(\omega)Tf, \quad f\in\dom(T), \;
\omega\in\cB_{\bbC}. \\
& (\gamma) \, \text{ Let $\omega\in\cB_{\bbC}$. Then } \,
T\big|_{E(\omega)\cH} \, \text{ on } \,
\dom\big(T\big|_{E(\omega)\cH} \big) = \dom(T)\cap E(\omega)\cH \,
\text{ has spec-} \\
& \quad \;\; \text{ trum contained in $\ol \omega$, }
\, \sigma\big(T\big|_{E(\omega)\cH}\big) \subseteq
\ol{\omega}.
\end{align*}
$E$ is called a {\it resolution of the identity} (or a {\it spectral
resolution}) for $T$.
\end{definition}

\begin{lemma} \lb{l2.15} ${}$ \\
$(i)$ Every spectral operator $T\in\cC(\cH)$ has a uniquely defined,
regular, countably additive resolution of the identity defined on
$\cB_{\bbC}$.
$($It will frequently be denoted by $E_T$ in the following.$)$ \\
$(ii)$ Let $T\in\cC(\cH)$ be a spectral operator with spectral resolution
$E_T$ and $\omega\in\cB_{\bbC}$. Then, $T\big|_{E(\omega)\cH}$ is a
spectral operator with spectral resolution
\begin{equation}
E_{T|_{E(\omega)\cH}}(\sigma) = E_T(\sigma)\big|_{E(\omega)\cH}, \quad
\sigma \in \cB_{\bbC}. \lb{2.11}
\end{equation}
If $\omega$ is bounded, $T\big|_{E(\omega)\cH}$ is a bounded
spectral operator. \\
$(iii)$ Let $T\in\cC(\cH)$ be a spectral operator, $E_T$ its spectral
resolution, and $\omega\in\cB_{\bbC}$ open. Then,
\begin{equation}
\sigma(T)\cap\omega \subseteq \sigma\big(T\big|_{E_T(\omega)}\big)
\subseteq \sigma(T)\cap \ol{\omega}. \lb{2.12}
\end{equation}
$(iv)$ If $P=P^2$, $P\dom(T)\subseteq \dom(T)$, $PTf=TPf$, $f\in\dom(T)$,
then,
\begin{equation}
\sigma(T)\supseteq \sigma\big(T\big|_{P\cH}\big).  \lb{2.13}
\end{equation}
$(v)$ Let $T\in\cC(\cH)$ be a spectral operator with spectral resolution
$E_T$. Then,
\begin{equation}
\sigma(T) = \bigcap_{\{\omega\in\cB_{\bbC} \,|\,E_T(\omega)=I_{\cH}\}}
\ol{\omega}.  \lb{2.14}
\end{equation}
$(vi)$ Let $T\in\cC(\cH)$ be a spectral operator with spectral resolution
$E_T$ satisfying $E_T(\sigma(T))={I_\cH}$. If
$\{\omega_n\}_{n\in\bbN}\subset\cB_{\bbC}$ is an increasing sequence of
bounded Borel sets with $E_T\big(\bigcup_{n\in\bbN}
\omega_n\big)=I_{\cH}$, then
\begin{equation}
\sigma(T)=\ol{\bigcup_{n\in\bbN}
\sigma\big(T\big|_{E_T(\omega_n)\cH}\big)}. \lb{2.14a}
\end{equation}
\end{lemma}

\begin{definition} \lb{d2.17}
Let $S\in\cC(\cH)$. Suppose there exists a regular, countably additive
projection-valued spectral measure $E$ (with respect to the strong
operator topology) defined on $\cB_{\bbC}$ such that
\begin{align}
\begin{split}
& \dom(S)=\bigg\{g\in\cH \,\bigg|\, \slimes_{n\uparrow\infty}
\int_{\{\lambda\in\bbC\,|\, |\lambda|\leq n\}} \lambda \, d(E(\lambda)\,
g) \text{ exists in $\cH$}\bigg\}, \\
& Sf =\slimes_{n\uparrow\infty} \int_{\{\lambda\in\bbC\,|\,|\lambda|\leq
n\}} \lambda \,  d(E(\lambda)\, f), \quad f\in\dom(S).  \lb{2.15}
\end{split}
\end{align}
Then $S$ is called a {\it spectral operator of scalar type} and the
projection-valued measure $E$ is called the {\it resolution of the
identity} for $S$.
\end{definition}

\medskip

\noindent {\bf Acknowledgments.}
We are grateful to Kwang Shin and Rudi Weikard for helpful discussions on
this subject.

\bigskip



\begin{thebibliography}{99}
%
\bi{Ba54} W.\ G.\ Bade, {\it Unbounded spectral operators}, Pac. J. Math.
{\bf 4}, 373--392 (1954).
%
\bi{BG06} V.\ Batchenko and F.\ Gesztesy, {\it On the spectrum of
Schr\"odinger operators with quasi-periodic algebro-geometric KdV
potentials}, J. Analyse Math. {\bf 95}, 333--387 (2005).
%
\bi{Bi86}  B.\ Birnir, {\it Complex Hill's equation and the complex
periodic Korteweg--de Vries equations }, Commun. Pure Appl. Math.
{\bf 39}, 1--49 (1986).
%
\bi{Bi86a} B.\ Birnir, {\it Singularities of the complex Korteweg--de
Vries flows}, Commun. Pure Appl. Math. {\bf 39}, 283--305 (1986).
%
\bi{Bi87} B.\ Birnir, {\it An example of blow-up, for the complex KdV
equation and existence beyond blow-up,} SIAM J. Appl. Math. {\bf 47},
710--725 (1987).
%
\bi{Ch06} T.\ Christiansen, {\it Isophasal, isopolar, and isospectral Schr\"odinger operators and elementary complex analysis}, Amer. J. Math., to appear. 
%
\bi{CL85} E.\ A.\ Coddington and N.~Levinson,{\it Theory of Ordinary
Differential Equations}, Krieger, Malabar, 1985.
%
\bi{CF68} I.\ Colojoar{\v a} and C.\ Foia{\c{s}}, {\it Theory of
Generalized Spectral Operators}, Gordon and Breach, New York, 1968.
%
\bi{DM02} P.\ Djakov and B.\ Mityagin, {\it Smoothness of Schr\"odinger
operator potential in the case of Gevrey type asymptotics of the gaps}, J.
Funct. Anal. {\bf 195}, 89--128 (2002).
%
\bi{DM03} P.\ Djakov and B.\ Mityagin, {\it Spectral gaps of the periodic
Schr\"odinger operator when its potential is an entire function}, Adv.
Appl. Math. {\bf 31}, 562--596 (2003).
%
\bi{DM03a} P.\ Djakov and B.\ Mityagin, {\it Spectral triangles of
Schr\"odinger operators with complex potentials}, Selecta Math. {\bf 9},
495--528 (2003).
%
\bi{DM06} P.\ Djakov and B.\ Mityagin, {\it Instability zones of 1D periodic 
Schr\"odinger and Dirac operators}, Uspehi Math. Nauk, {\bf 61}, 77--183 (2006).
%
\bi{Du58} N.\ Dunford, {\it A survey of the theory of spectral
operators}, Bull. Amer. Math. Soc. {\bf 64}, 217--274 (1958).
%
\bi{DS88}  N.\ Dunford and J.\ T Schwartz, {\it Linear Operators, Part
II: Spectral Theory}, Wiley--Interscience, New York, 1988.
%
\bi{DS88a}  N.\ Dunford and J.\ T Schwartz, {\it Linear Operators, Part
III: Spectral Operators}, Wiley--Interscience, New York, 1988.
%
\bi{Ea73}  M.\ S.\ P.\ Eastham, {\it The Spectral Theory of Periodic
Differential Equations}, Scottish Academic Press, Edinburgh and
London, 1973.
%
\bi{EFZ05} W.\ Eberhard, G.\ Freiling, and A.\ Zettl, {\it Sturm--Liouville
problems with singular non-selfadjoint boundary conditions}, Math. Nachr.,
to appear.
%
\bi{Fo58} S.\ R.\ Foguel, {\it The relation between a spectral operator
and its scalar part}, Pac. J. Math. {\bf 8}, 51--65 (1958).
%
\bi{Fu70} V.\ N.\ Funtakov, {\it Expansions in eigenfunctions of
nonself-adjoint second-order differential equations}, Diff. Eq. {\bf 6},
1528--1535 (1970).
%
\bi{Ga80} M.\ G.\ Gasymov, {\it Spectral analysis of a class of
second-order non-self-adjoint differential operators}, Funct. Anal.
Appl. {\bf 14}, 11--15 (1980).
%
\bi{Ga80a} M.\ G.\ Gasymov, {\it Spectral analysis of a class of
ordinary differential operators with periodic coefficients}, Sov.
Math. Dokl. {\bf 21}, 718--721 (1980).
%
\bibitem {Ge50}   I.\ M.\ Gel'fand,
{\em Expansion in characteristic functions of an equation with periodic
coefficients}, Doklady Akad Nauk SSSR {\bf 73}, 1117-1120 (1950).
(Russian.)
%
\bi{GT06} F.\ Gesztesy and V.\ Tkachenko, {\it When is a non-self-adjoint Hill operator a spectral operator of scalar type?}, C. R. Acad. Sci. Paris, Ser. I, {\bf 343}, 239--242 (2006).
%
\bi{GW95}  F.\ Gesztesy and R.\ Weikard, {\it Floquet theory revisited},
in {\it Differential Equations and Mathematical Physics}, I.\ Knowles
(ed.), International Press, Boston, 1995, pp.\ 67--84.
%
\bi{GW96} F.\ Gesztesy and R.\ Weikard, {\it Picard potentials
and Hill's  equation on  a torus}, Acta Math. \textbf{176}, 73--107
(1996).
%
\bi{GW98} F.\ Gesztesy and R.\ Weikard, {\it A characterization
of all elliptic  algebro-geometric solutions of the AKNS hierarchy}, Acta
Math. {\bf 181}, 63--108 (1998).
%
\bi{GW98a} F.\ Gesztesy and R.\ Weikard, {\it Elliptic
algebro-geometric solutions of the KdV and AKNS hierarchies -- an
analytic approach}, Bull. Amer. Math. Soc. {\bf 35}, 271--317 (1998).
%
\bi{GU83} V.\ Guillemin and A.\ Uribe, {\it Hardy functions and the
inverse spectral method}, Commun. PDE {\bf 8}, 1455--1474 (1983).
%
\bi{In56} E.\ L.\ Ince, {\it Ordinary Differential Equations},
Dover, New York, 1956.
%
\bi{Ke64} G.\ M.\ Kesel'man, {\it On the unconditional convergence of eigenfunction expansions of certain differential operators}, Izv. Vyssh. Uchebn. Zaved. Mat. 
{\bf 39} (2), 82--93 (1964). (Russian.)
%
\bi{Ko97} S.\ Kotani, {\it Generalized Floquet theory for stationary
Schr\"odinger operators in one dimension}, Chaos, Solitons \& Fractals
{\bf 8}, 1817--1854 (1997).
%
\bi{Le96} B.\ Ya.\ Levin, {\it Lectures on Entire Functions}, Transl. 
Math. Monographs, {\bf 150}, Amer. Math. Soc., Providence, RI, 1996.
%
\bi{Lj66} V.\ \`E. Ljance, {\it On a generalization of the concept of
spectral measure}, Amer. Math. Soc. Transl. (2) {\bf 51}, 273--315 (1966).
%
\bi{LM60} U.\ I.\ Lyubi\v c and V.\ I.\ Macaev, {\it On the spectral theory
of linear operators in Banach spaces}, Soviet Math. Dokl.  {\bf 1},
184--186 (1960).
%
\bi{LM62} Ju.\ I.\ Lyubi\v c and V.\ I.\ Macaev, {\it Operators with
separable spectrum}, Mat. Sb. (N.S.)  {\bf 56} (98), 433--468 (1962).
(Russian.)
%
\bi{Ma06} A.\ Makin, {\it On periodic boundary value problem for the 
Sturm--Liouville operator}, preprint, arXiv:math.SP/0601436.
%
\bi {Ma63} V.\ A.\ Marchenko, {\it Expansion in eigenfunctions of
non-self-adjoint singular differential operators of second order}, Amer.
Math. Soc. Transl. (2) {\bf 25}, 77--130 (1963).
%
\bi {Ma86} V.\ A.\ Marchenko, {\it Sturm--Liouville operators and
applications}, Birkh\"auser, Basel, 1986.
%
\bi {MO75} V.\ A.\ Marchenko and I.\ V.\ Ostrovskii, {\it A
characterization of the spectrum of Hill's operator}, Math. USSR Sb.
{\bf 26}, 493--554.
%
\bi{Mc62} D.\ McGarvey, {\it Operators commuting with translations
by one. Part I. Representation theorems}, J. Math. Anal. Appl. {\bf 4},
366-410 (1962).
%
\bi{Mc65} D.\ McGarvey, {\it Operators commuting with translations
by one. Part II. Differential operators with periodic coefficients in
$L_p(-\infty,\infty)$}, J. Math. Anal. Appl. {\bf 11}, 564-596 (1965).
%
\bi{Mc65a} D.\ McGarvey, {\it Operators commuting with translations
by one. Part III. Perturbation results for periodic differential
operators}, J. Math. Anal. Appl. {\bf 12}, 187--234 (1965).
%
\bi{Mc66} D.\ McGarvey, {\it Linear differential systems with periodic
coefficients involving a large parameter}, J. Diff. Eq. {\bf 2},
115--142 (1966).
%
\bi{Me77} N.\ N.\ Meiman, {\it The theory of one-dimensional Schr\"odinger
operators with a periodic potential}, J. Math. Phys. {\bf 18}, 834--848
(1977).
%
\bi{Mi62} V.\ P.\ Miha{\u i}lov, {\it Riesz bases in $L_2(0,1)$}, Sov. Math. Dokl. 
{\bf 3}, 851--855 (1962). 
%
\bi{Mi06} A.\ Minkin, {\it Resolvent growth and Birkhoff-regularity}, J. Math. 
Anal. Appl. {\bf 323}, 387--402 (2006).
%
\bibitem {Na60} M.\ A.\ Naimark, {Investigation of the spectrum and the
expansion in eigenfunctions of a non-selfadjoint differential operator of
the second order on a semi-axis}, Amer. Math. Soc. Transl. (2) {\bf 16},
103--193 (1960).
%
\bibitem {Na68} M.\ A.\ Naimark, {\it Linear Differential operators, Part
II}, Ungar, New York, 1968.
%
\bi{Pa79} B.\ S.\ Pavlov, {\it Basicity of an exponential system and
Muckenhoupt's condition},  Sov. Math. Dokl.  {\bf 20}, 655--659 (1979).
%
\bi{PT88}  L.\ A.\ Pastur and V.\ A.\ Tkachenko, {\it Spectral theory
of Schr\"odinger operators with periodic complex-valued potentials},
Funct. Anal. Appl. {\bf 22}, 156--158 (1988).
%
\bi{PT91} L.\ A.\ Pastur and V.\ A.\ Tkachenko, {\it An inverse
problem for a class of one-dimensional Schr\"odinger operators with a
complex periodic potential}, Math. USSR Izv. {\bf 37}, 611--629 (1991).
%
\bi{PT91a} L.\ A.\ Pastur and V.\ A.\ Tkachenko, {\it Geometry of the
spectrum of the one-dimensional Schr\"odinger equation with a periodic
complex-valued potential}, Math. Notes {\bf 50}, 1045--1050 (1991).
%
\bi{RS78} M.\ Reed and B.\ Simon, {\it Methods of Modern
Mathematical Physics. IV: Analysis of Operators}, Academic Press, New
York, 1978.
%
\bi{Ro63} F.\ S.\ Rofe-Beketov, {\it The spectrum of non-selfadjoint
differential operators with periodic coefficients}, Sov. Math. Dokl.
{\bf 4}, 1563--1566 (1963).
%
\bi{ST96} J.-J.\ Sansuc and V.\ Tkachenko, {\it Spectral
parametrization of non-selfadjoint Hill's operators}, J. Diff. Eq.
{\bf 125}, 366--384 (1996).
%
\bi{ST96a} J.-J.\ Sansuc and V.\ Tkachenko, {\it Spectral properties of
non-selfadjoint Hill's operators with smooth potentials}, in {\it Algebraic
and Geometric Methods in Mathematical Physics}, A.\ Boutel de
Monvel and V.\ Marchenko (eds.), Kluwer, Dordrecht, 1996, pp.\ 371--385.
%
\bi{ST97} J.-J.\ Sansuc and V.\ Tkachenko, {\it Characterization of the
periodic and  antiperiodic spectra of nonselfadjoint Hill's operators,}
in {\it New  Results in Operator Theory and its Applications},
I.\ Gohberg and  Yu.\ Lubich eds.), Operator Theory: Advances and
Applications {\bf 98},  Birkh\"auser, Basel, 1997, pp.\ 216--224.
%
\bi{Sc60} J. Schwartz, {\it Some non-selfadjoint operators}, Commun.
Pure Appl. Math. {\bf 13}, 609--639 (1960).
%
\bi{Se60} M.\ I.\ Serov, {\it Certain properties of the spectrum of a
non-selfadjoint differential operator of the second order}, Sov. Math.
Dokl. {\bf 1}, 190--192 (1960).
%
\bi{Sh03} K.\ C.\ Shin, {\it On half-line spectra for a class of
non-self-adjoint Hill operators}, Math. Nachr. {\bf 261--262}, 171--175
(2003).
%
\bi{Sh04} K.\ C.\ Shin, {\it Trace formulas for non-self-adjoint
Schr\"odinger operators and some applications}, J. Math. Anal. Appl.
{\bf 299}, 19--39 (2004).
%
\bi{Sh04a} K.\ C.\ Shin, {\it On the shape of spectra for non-self-adjoint
periodic Schr\"odinger operators}, J. Phys. A {\bf 37}, 8287--8291 (2004).
%
\bi{Ta51} A.\ E.\ Taylor, {\it Spectral theory of closed distributive
operators}, Acta Math. {\bf 84}, 189--224 (1951).
%
\bi{Ti50} E.\ C.\ Titchmarsh, {\it Eigenfunction problems with periodic
potentials}, Proc. Roy. Soc. London A {\bf 203}, 501--514 (1950).
%
\bi{Ti58} E.\ C.\ Titchmarsh, {\it Eigenfunction Expansions
associated with Second-Order Differential Equations, Part II}, Oxford
University Press, Oxford, 1958.
%
\bi{Tk64} V.\ A.\ Tkachenko, {\it Spectral analysis of the one-dimensional
Schr\"odinger operator with periodic complex-valued potential}, Sov.
Math. Dokl. {\bf 5}, 413--415 (1964).
%
\bi{Tk92} V.\ A.\ Tkachenko, {\it Spectral analysis of a nonselfadjoint
Hill operator}, Sov. Math. Dokl. {\bf 45}, 78--82 (1992).
%
\bi{Tk94} V.\ A.\ Tkachenko, {\it Discriminants and Generic Spectra
of Nonselfadjoint Hill's Operators}, Adv. Sov. Math. {\bf 19}, 41--71
(1994).
%
\bi{Tk96} V.\ A.\ Tkachenko, {\it Spectra of non-selfadjoint
Hill's operators and a class of Riemann surfaces}, Ann. Math. {\bf
143}, 181--231 (1996).
%
\bi{Tk01} V.\ Tkachenko, {\it Characterization of Hill operators with
analytic potentials}, Integral equ. oper. theory {\bf 41}, 360--380
(2001).
%
\bibitem {Tk02} V.\ Tkachenko, {\em Non-selfadjoint Sturm-Liouville
operators with multiple spectra}, in {\it Interpolation Theory,
Systems Theory and Related Topics}, D. Alpay, I. Gohberg, V.
Vinnikov (eds.), Operator Theory: Advances and Applications, {\bf 134},
403--414 (2002).
%
\bi{Ve80} O.\ A.\ Veliev, {\it The one-dimensional Schr\"odinger operator
with a periodic complex-valued potential}, Sov. Math. Dokl. {\bf 21},
291--295 (1980).
%
\bi{Ve83} O.\ A.\ Veliev, {\it Spectrum and spectral singularities of
differential operators with complex-valued periodic coefficients}, Diff.
Eqs. {\bf 19}, 983--989 (1983).
%
\bi{Ve86} O.\ A.\ Veliev, {\it Spectral expansions related to
non-self-adjoint differential operators with periodic coefficients}, Diff.
Eqs. {\bf 22}, 1403--1408 (1986).
%
\bi{VT02} O.\ A.\ Veliev and M.\ Toppamuk Duman, {\it The spectral
expansion for a nonself-adjoint Hill operator with a locally integrable
potential}, J. Math. Anal. Appl. {\bf 265}, 76--90 (2002).
%
\bi{Vo63} V.\ Ja.\ Volk, {\it Spectral resolution for a class of
non-selfadjoint operators}, Sov. Math. Dokl. {\bf 4}, 1279--1281 (1963).
%
\bibitem{We98} R.\ Weikard, {\it On Hill's equation with a singular
complex-valued potential}, Proc. London Math. Soc. {\bf 76}, 603-633
(1998).
%
\bibitem{We98a} R.\ Weikard, {\it On a theorem of Hochstadt}, Math.
Ann. {\bf 311}, 95-105 (1998).
%
\bibitem{Zh69} V.\ A.\ Zheludev, {\it Perturbations of the spectrum of the
Schroedinger operator with a complex periodic potential}, in {\it Topics
in Mathematical Physics {\bf 3}, Spectral Theory}, M. Sh. Birman (ed.),
Consultants Bureau, New York (1969), pp.\ 25--41.
%
\end{thebibliography}
\end{document}